\documentclass[11pt]{amsart}
\usepackage{lmodern} 
\usepackage[T1]{fontenc}
\usepackage{textcomp}
\usepackage{stmaryrd}
\usepackage{enumerate}
\usepackage{amsfonts}
\usepackage{mathrsfs}
\usepackage{bbm}
\usepackage{amssymb}
\usepackage{amsmath,amssymb, amsthm}

\usepackage{color}

\usepackage{geometry}
\usepackage[colorlinks,
            linkcolor=blue,
            anchorcolor=blue,
            citecolor=blue
            ]{hyperref}
\let\newpf\proof \let\proof\relax 
\newenvironment{pf}{\newpf[\proofname]}{\qed\endtrivlist}

\newcommand{\bQ}{\overline{Q}}
\newcommand{\ba}{\overline{A}}

\newcommand{\CD}{\rm CD}

\def\be{\begin{equation}}
\def\ee{\end{equation}}

\def\ba{{\begin{align}}}
\def\ea{{\end{align}}}

\def\bm{\begin{matrix}}
\def\em{\end{matrix}}

\def\0{{\mathbf 0}}

\def\cal{\mathcal}

\newtheorem{Theorem}{Theorem}[section]
\newtheorem{Lemma}[Theorem]{Lemma}
\newtheorem{Proposition}[Theorem]{Proposition}
\newtheorem{Conjecture}[Theorem]{Conjecture}
\newtheorem{Corollary}[Theorem]{Corollary}
\newtheorem{Remark}[Theorem]{Remark}
\newtheorem{Definition}[Theorem]{Definition}

\numberwithin{equation}{section}
\newcommand{\thmref}[1]{Theorem~\ref{#1}}

\def \bn {\hfill \\ \smallskip\noindent}

\theoremstyle{definition}

\newtheorem{definition}{Definition}[section]

\def\proof{\bn {\bf Proof.} }

\renewcommand{\mod}{\operatorname{mod}}

\newcommand{\C}{{\mathbb C}}
\newcommand{\D}{{\mathbb D}}
\newcommand{\E}{{\mathbb E}}

\newcommand{\N}{{\mathbb N}}
\newcommand{\Q}{{\mathbb Q}}
\newcommand{\R}{{\mathbb R}}
\newcommand{\T}{{\mathbb T}}

\newcommand{\Z}{{\mathbb Z}}

\def\B0{{\bold{0}}}

\newcommand{\la}{\langle}
\newcommand{\ra}{\rangle}

\catcode`\@=12

\def\Empty{}
\newcommand\oplabel[1]{
  \def\OpArg{#1} \ifx \OpArg\Empty {} \else
    \label{#1}
  \fi}

\newcommand{\comm}[1]{}
\newcommand{\comment}[1]{}

\begin{document}

\title[Mobius disjointness]{Quantitative almost reducibility and  M\"obius disjointness for analytic quasiperiodic Schrodinger cocycles}

\author{Wen Huang, Jing Wang, Zhiren Wang and Qi Zhou}
\address[W.~Huang]{CAS Wu Wen-Tsun Key Laboratory of Mathematics, School of Mathematical Sciences, University of Science and Technology of China,
	Hefei, Anhui, 230026, P.R. China}
 \email{wenh@mail.ustc.edu.cn}
\address[J.~Wang]{Department of Mathematics, School of Science, Nanjing University of Science and Technology, Nanjing, 210094, P.R. China}
\email{jing.wang@njust.edu.cn}
\address[Z.~Wang]{Department of Mathematics, Pennsylvania State University, University Park,
PA 16802, USA}
\email{zhirenw@psu.edu}
\address[Q.~Zhou]{Chern Institute of Mathematics and LPMC, Nankai University, Tianjin 300071, China}
\email{qizhou@nankai.edu.cn}

\date{\today}

\pagestyle{headings}

\maketitle

\begin{abstract}
Sarnak's M\"obius disjointness conjecture states that M\"obius function is disjoint to any zero entropy dynamics. 
We prove that M\"obius disjointness conjecture holds for one-frequency analytic quasi-periodic cocycles which are  almost reducible, which extend \cite{LS15,W17} to the noncommutative case. The proof relies on quantitative version of almost reducibility. 
\end{abstract}

\section{Introduction}\label{sec-intro}

A quasi-periodic $SL(2,\R)$-cocycle is a linear skew product of the form
\[
(x,\varpi)\mapsto (x+\alpha, A(x)\varpi)
\]
where $x\in\T^1:=\R^1/\Z^1$, $\varpi\in\R^2$, $\alpha\in\R\backslash\Q$ 
 and $A\in C^0(\T^1, SL(2,\R))$, denoted by $(\alpha, A)$.
We consider the projective action of the $SL(2,\R)$-cocycles $(\alpha, A)$ on $\R\mathbb P^1$. Those are quasiperiodically forced (qpf) circle homeomorphisms of the form
\[T_{(\alpha, A)}\ : \ \T^1\times\R\mathbb P^1\rightarrow \T^1\times\R\mathbb P^1,  \ \ (x, \varphi)\mapsto (x+\alpha,\frac{A(x)\cdot\varphi}{\|A(x)\cdot\varphi)\|}),\]
where $\cdot$ denotes the M\"obius transformation.

The topological dynamics of the projective action of a quasi-periodic $SL(2,\R)$ cocycle is a very interesting subject in itself \cite{BJ09, Her81, Her83, HYi09, J09, WJ, WZJ, Z}, which has zero 
 topological entropy \cite{HYi09}. In this paper, we focus on the question that whether this action satisfies Sarnak's M\"obius disjointness conjecture.

The  M\"obius function $\mu:\N\rightarrow \{-1,0,1\}$ is defined by $\mu(1)=1$ and
\[
\mu(n)=\left\{\begin{array}{ll}(-1)^k & \ \textrm{if}\ n\ \textrm{is a product of} \ k \ \textrm{distinct primes};\\
0 & \ \textrm{otherwise}.
\end{array}\right.
\]
The orthogonality of $\mu$ with other sequences is a very important issue 
in number theory. For instance, the disjointness of M\"obius function with constant sequence, that is $\sum_{n\leq N}\mu(n)=o(N)$, is equivalent 
to the prime number theorem. The M\"obius randomness law \cite{IK04} suggests that $\mu(n)$ has significant cancellations with reasonable sequences $\xi(n)$, which means
\begin{equation}\label{EqMD}\sum_{n\leq N}\mu(n)\xi(n)=o(N).
\end{equation}

Sarnak's M\"obius disjointness conjecture \cite{S09} expects the disjointness holds for deterministic sequences, namely:

\begin{Conjecture}[M\"obius disjointness conjecture] For a topological dynamical system (t.d.s.) $(X, T)$ of zero topological entropy, all $f\in C(X)$ and all $x\in X$, the sequence $\xi(n)=f(T^nx)$ satisfies \eqref{EqMD}.\end{Conjecture}

There have been many partial results on M\"obius disjointness conjecture. 
For the progress in this conjecture, we will simply refer to the recent comprehensive survey \cite{FKPL18} for brevity, and we only discuss the historical developments that are more relevant to this paper. We remark that in the special case when $A(x)\in SO(2):$ $$A(x)=R_{g(x)} :=\left(\begin{array}{cc}\cos(2\pi g(x)) & -\sin(2\pi g(x))\\ \sin(2\pi g(x)) & \cos(2\pi g(x))\end{array}\right)$$ for some $g:\T^1\to \R$, i..e., 
 $T_{(\alpha,A)}$ has the form  $(x, \varphi)\mapsto (x+\alpha, \varphi+g(x))$, then the action $T_{(\alpha,A)}$ is known to satisfy the M\"obius disjointness conjecture as long as $g$ is $C^{1+\epsilon}$ by the recent work of de Faveri \cite{dF20} (see also the earlier works \cite{HWY17,KLR19,KL13,LS15,W17}, which assumed higher regularity and/or other conditions).

Therefore the natural question is whether M\"obius  function is disjoint for any cocycle $(\alpha, A)$, or what happens if the non-commutative part appears? It is easy to see, if the cocycles has positive  Lyapunov exponent, then M\"obius disjointness conjecture holds (\S \ref{SecHyperbolic}). Thus the real difficulty lies in if the cocycle has zero Lyapunov exponent. From this respect, we mention the following important concept raised by Eliasson \cite{E92,E01}.  Recall that $(\alpha, A)$  is 
$C^\omega$-conjugated to $(\alpha,\tilde{A})$, if  there exists $B\in C^\omega(\T, PSL(2,\R))$ such that  $$B(\cdot+\alpha)^{-1}A(\cdot)B(\cdot)=\tilde{A}(\cdot).$$ Then $(\alpha, A)$ is said to be  $C^\omega$-almost reducible if its closure of analytic conjugate class contains a constant cocycle. 

%
%


\begin{Theorem}\label{mainthm-1}
Let $(\alpha, A)\in \R\backslash\Q\times C^\omega(\T^1, SL(2,\R))$. If $(\alpha, A)$ is almost reducible, then the M\"obius disjointness conjecture holds for $(\T^1\times\R\mathbb P^1, T_{(\alpha, A)})$.
\end{Theorem}

One may ask, compared to the whole space of analytic cocycle, how large is the space of almost reducible cocycle is. Indeed, Avila's global theory 
of analytic $SL(2,\R)$ cocycles \cite{A3,A2,A15} stated that  typical cocycle is almost reducible or has positive Lyapunov exponent. To precise this, let's take quasi-periodic Schr\"odinger cocycles
\[A(x)=S_{v,E}(x)=\left(\begin{array}{ccc}
E-v(x)& -1\cr
1 & 0
\end{array}\right),\]
as a typical example.  Quasi-periodic Schr\"odinger cocycles are the main 
source of examples for quasi-periodic $SL(2,\R)$ cocycles, because of their relation to
one-dimensional discrete quasi-periodic Schr\"{o}dinger operators
\begin{equation}\label{equ-sch-op}
 (H_{v,\alpha,x}u)_n=u_{n+1}+u_{n-1}+v(x+n\alpha)u_n,   \forall n\in\Z,
\end{equation}
where $x\in\T^1$ is the phase, $v\in C^0(\T^1,\R)$ is the potential, and $\alpha\in\T^1$ is the frequency. Then as a direct  corollary  of Avila's 
global theory and \thmref{mainthm-1}, we have the following:

\begin{Corollary}\label{mainthm-schrodinger}
Let $\alpha\in\R\backslash\Q$.
Then for  a (measure-theoretically) typical $v\in C^\omega(\T^1, \R)$, the M\"obius disjointness conjecture holds for the dynamical system $(\T^1\times\R\mathbb P^1, T_{(\alpha, S_{v,E})})$ induced by the projective action of $(\alpha, S_{v,E})$ for any $E\in\R$.
\end{Corollary}

Let us explain  the meaning of ``typical'' in the theorem as in \cite{A15}. Since in the infinite-dimensional settings one lacks a translation-invariant measure (Haar measure), it is common to replace the notion of  {\it almost every} by {\it prevalence}: one fixes some probability measure $\mu$ of compact support (a set of admissible perturbations $w$), and declares a property to be typical if it is satisfied for almost every perturbation $v+w$ of every starting condition $w$. In finite-dimensional vector 
spaces, prevalence implies full Lebesgue measure.  For instance,  one can 
set $\Delta=\D^{\N}$ endowed with the probability measure $\mu$ given by the product of normalized Lebesgue measure. Given an arbitrary function 
$\epsilon:\N \rightarrow \R_{+}$ which decays exponentially fast, we associate a probability measure $\mu_{\epsilon}$ with compact support on $C^\omega(\T^1, \R)$ by push forward of $\mu$ under the map
$$ \{t_m\}_{m\in\Z} \rightarrow \sum_{m\geq 1} \epsilon(m)2\Re[t_m e^{2\pi i mx}].$$
Thus Corollary  \ref{mainthm-schrodinger} says that, for any $\alpha\in \R\backslash\Q$, and for every $v\in C^\omega(\T^1, \R)$, the M\"obius disjointness conjecture holds for the dynamics $T_{(\alpha, S_{v+w,E})}$ for 
$\mu_{\epsilon}$-almost every $w$ and every $E\in\R$.

\begin{Remark}
In essence, Corollary \ref{mainthm-schrodinger} relies on Avila's \textit{Almost Reducibility Conjecture}(ARC), which asserts that any  \textit{subcritical cocycle}\footnote{Consult  section \ref{global-th} for precise definition.}  is almost reducible. A proof of ARC is announced in \cite{A15}, to appear in \cite{A3,A2}. ARC has many important dynamical and spectral consequences \cite{A3,A15,AJM,AKL,AYZ, GY,LYZZ,MJ}, indeed, it was already stated as  \textit{Almost Reducibility Theorem}(ART) in \cite{AJM}.

\end{Remark}

If we take  $v(x)=2\lambda \cos(2\pi x)$, which corresponds to almost Mathieu operators, then we have the following:
\begin{Corollary}\label{mainthm-schrodinger-1}
Let $\alpha\in\R\backslash\Q$,  the M\"obius disjointness conjecture holds for the dynamics $(\T^1\times\R\mathbb P^1, T_{(\alpha, S_{2\lambda \cos(2\pi\cdot),E})})$  for any $\lambda \neq \pm 1$, and for any $E\in\R$.
\end{Corollary}

We take Schr\"odinger cocycles (and almost Mathieu operators) as typical examples, since it is interesting from both physical and mathematical points of view. Physically, quasi-periodic Schr\"{o}dinger operator is a mathematical model for quasicrystals, and it plays a central role in explaining the quantum Hall effect (which was Thouless's Nobel Prize work)  \cite{L1,TKNN}. From the mathematical side,  Schr\"odinger cocycles (even almost Mathieu cocycles)  could display all the dynamical behavior of general $SL(2,\R)$ cocycles \cite{A15}.  For example, Avila-Fayad-Krikorian \cite{AFK11} proved that for any $\alpha\in\R\backslash\Q$, $v\in C^\omega(\T^1, \R)$ and for almost every energy $E\in\R$ the analytic cocycle $(\alpha, S_{v,E})$ is either $C^\omega$-rotations reducible or has positive Lyapunov exponent. Recall
$(\alpha, A)$ is called \textit{$C^\omega$-rotations reducible}, if there 
exists $B\in C^\omega(\T, PSL(2,\R))$ such that $$B(\cdot+\alpha)^{-1}A(\cdot)B(\cdot)\in SO(2,\R).$$ Thus
Theorem \ref{mainthm-schrodinger} can be viewed as a direct non-commutative extension of the results from \cite{LS15,W17}.  Secondly, to better understand zero  topological entropy,  complexity is a very useful concept, 
while complexity is closely related to the growth of cocycles, and Schr\"odinger cocycles (even almost Mathieu cocycles) are believed to  display high order complexity \cite{ALSZ}.

Finally, we outline the novelty and difficulty of the proof.  The concept 
of almost reducibility generalizes the scope of applicability of the local theories. While quantitative estimates  were not involved in the definition, we stress that the main novelty  of the paper  is that  quantitative version of  almost reducibility implies the disjointness. 
Theorem \ref{mainthm-1} covers all frequencies, and the difficulty in its 
proof lies in two possible issues: the frequency being Liouvillean and the cocycle being almost reducible to parabolic cocycles. 

If the frequency is Diophantine, then in the commutative case  \cite{LS15,W17},  the M\"obius disjointness is a straightforward corollary of Davenport's theorem \cite{D37}; while in our case,  quantitative almost reducibility was obtained  based on almost localization of the dual model, as a 
consequence, we prove that the cocycle has sub-polynomial complexity, which leads to M\"obius disjointness by  \cite{HWY17}. When the frequency is 
Liouvillean, almost reducibility was  previous obtained by different methods: periodic approximation as in \cite{A3}, and non-standard KAM as in \cite{HY12}.  However,  these estimates are far from our needs. Indeed, uniform estimates in the parameters are usually not useful, since counterexamples \cite{ALSZ} and even phase transitions \cite{AYZ} could happen.  
On the other hand,  to get precise estimates for any frequency,  one has to  interpolate  non-standard KAM  schemes \cite{HY12} with classical KAM 
schemes \cite{E92}.

Another difficulty is that $(\alpha,A)$ may be almost reducible to parabolic cocycles, i.e.  there exists conjugacies $B_n$  which conjugate $(\alpha,A)$ to  perturbations of $(\alpha,A_n)$, where 
$A_n=\left(\begin{array}{cc} 1 & c_n\cr 0 & 1 \end{array}\right)$ is parabolic. This is a new difficulty in the non-commutative case. In this case, one needs delicate control on the sizes of $c_n$ and the conjugacy $B_n$ (Proposition \ref{prop-quantitative-almost-reducibility-2}). We point 
out the same difficulty also occurs in Avila-You-Zhou's work \cite{AYZ16} 
on the dry Martini problem. However, Aubry duality, a key tool in \cite{AYZ16},  cannot be used here. To deal with this obstruction, we distinguish two cases: if $c_n$ is relatively small,  then one can consider it as a 
perturbation of identity, and still apply the measurable complexity criterion of Huang-Wang-Ye \cite{HWY17};  if $c_n$ is relatively large, we will approximate the trajectories of the projective cocycles with either rotational or periodic trajectories on sufficiently long intervals,  and use 
the Matom\"aki-Radziwi\l\l-Tao estimate; however, if $c_n$ is too large (even of constant size), then one will lose control of the conjugacy $B_n$, which destroy the whole KAM scheme. Our quantitative argument will provide bounds on $c_n$ that prevents this from happening.

\smallskip

The remainder of the paper is organized as follows. In section \ref{sec-2}, we give some notations and collect some necessary background which will be used in what follows. In section \ref{sec-3}, we give the proof of the main theorems, mainly Theorem \ref{mainthm-1}. If the cocycle is almost reducible and not uniformly hyperbolic, we get the quantitative almost reducibility result for $\beta(\alpha)=0$ and $\beta(\alpha)>0$ respectively, and in the end deduce M\"obius disjointness according to the type of the constant matrices in the conjugated systems (elliptic or parabolic). In section \ref{sec-almost-reducibility-1}, we give the quantitative almost reducibility for $\beta(\alpha)=0$, using Aubry duality and  quantitative almost localization. In section \ref{sec-almost-reducibility-2}, 
we give the quantitative almost reducibility for $\beta(\alpha)>0$ by KAM 
scheme. In section \ref{sec-measure-complexity}, we get the measure complexity of the cocycle for the case that the constant matrices in the conjugated cocycles are elliptic. In section \ref{sec-parabolic}, we prove the 
linear disjointness of the dynamics of the system with M\"obius function for the case  that the constant matrices in the conjugated cocycles are parabolic.

\section{Notations and preliminaries}\label{sec-2}

\subsection{Notations on functions}
Denote by $C^\omega_{h}(\T^d,*)$ the
set of all  $*$-valued functions ($*$ will usually denote $\R$, $sl(2,\R)$
$SL(2,\R)$ and etc.) which are analytic and bounded in  $ \{ x \in \C^d/\Z^d\  : \  | \Im x |< h \}$, and for 
any $F\in C^\omega_{h}(\T^d,*)$, we define the norm
\[\|F\|_h:=  \sup_{ | \Im x |< h } \| F(x)\|,\]
(where $\|\cdot\|$ denotes absolute value or the usual matrix norm). 

Moreover, an integrable  real-valued function $f$ on $\T^{d}$ has the Fourier expansion $f(\phi)=\sum_{ k\in\Z^{d}}\hat{f}( k)e^{2\pi i
\langle  k, \phi\rangle }$ with $\hat{f}(k)=\int_{\T^d}f(\phi)e^{- 2\pi i
\langle  k, \phi\rangle }d\phi$. For any $K>0$, $\mathcal{T}_K$ and $\mathcal{R}_K$  are used
to denote the truncation operators:
\[
\mathcal{T}_N (f)=\sum_{| k|<K} \hat f(k)e^{2\pi i\la  k,\phi\ra},
\  \mathcal{R}_N (f) =\sum_{ |k|\geq K} \hat f(k) e^{2\pi i
\la k,\phi\ra}.
\]

\subsection{Continued fraction expansion and CD bridge}\label{subsec-continued-fraction-expansion}
Let $\alpha\in\R^1$ be irrational. Define $a_0=[\alpha],
\alpha_{0}=\alpha-a_0,$ and inductively for $k\geq 1$,
$$a_k=[\alpha_{k-1}^{-1}],\qquad \alpha_k=\{ \alpha_{k-1}^{-1}\}=\alpha_{k-1}^{-1}-a_k.$$
We define $p_0=0,  p_1=1, q_0=1, q_1=a_1$ and inductively,
\begin{eqnarray*}
 p_k=a_kp_{k-1}+p_{k-2}, \quad  q_k=a_kq_{k-1}+q_{k-2}.
\end{eqnarray*}
Then  the sequence $\{q_n\}_{n\in\N}$ satisfies
\begin{eqnarray*}
&& \forall 1 \leq k < q_n,\quad \|k\alpha\|_{\T} \geq
\|q_{n-1}\alpha\|_{\T},\\
&& \frac{1}{q_n+q_{n+1}}\leq \inf_{p\in\Z}|q_n \alpha-p|\leq {1 \over
q_{n+1}}.
\end{eqnarray*}
Let $\beta(\alpha):=\limsup_{n\rightarrow \infty}\frac{\ln q_{n+1}}{q_n}$. Equivalently, we have 
\begin{equation*}
\beta(\alpha)=\limsup_{k\rightarrow\infty}\frac{1}{|k|}\ln\frac{1}{\|k\alpha\|_{\T}},
\end{equation*}
where $\|k\alpha\|_{\T}:=\inf_{j\in\Z}|k\alpha-j|$.

Let us now introduce the CD bridge.

\begin{Definition}[{\cite{AFK11}}]\label{CDbridge}
Let $0<\mathcal A \leq \mathcal B\leq \mathcal C$. We say that the pair of denominators
$(q_l,q_n)$ forms a $\CD(\cal A,\cal B,\cal C)$ bridge if
\begin{itemize}
\item $q_{i+1}\leq q_i^{ \cal A}, \quad \forall i=l,\ldots,n-1$
\item $q_l^{\cal C}\geq q_n\geq q_l^{\cal B}$.
\end{itemize}
\end{Definition}

 In the following, for simplicity, we will fix a subsequence $\{q_{n_k}\}_k$ of
$\{q_n\}_n$, denoted by $\{Q_k\}_k$, and the subsequence $\{q_{n_k+1}\}_k$,
denoted by $\{\bQ_k\}_k$.

The following lemma was proved by Avila-Fayad-Krikorian in \cite{AFK11}, and we include the proof for completeness.
\begin{Lemma}[{\cite{AFK11}}]\label{CDbridge}
For any ${\cal A}>0$, there exists  a subsequence $\{Q_k\}_{k\in\N_0}$ such that $Q_0=1$ and for each $k\geq 1$,
$Q_{k}\leq \bQ_{k-1}^{{\cal A}^4}$. Furthermore, either $\bQ_k\geq
Q_k^{\cal A}$, or the pairs $(\bQ_{k-1},Q_{k})$ and $(Q_k,Q_{k+1})$
are both $\CD({\cal A},{\cal A},{\cal A}^3)$ bridges. 
\end{Lemma}
\begin{pf}
Suppose for $l\leq k$, we have already selected $Q_l$ satisfying the conditions. If  $q_{n+1}\leq q_n^{\cal A}$ for any $q_n>Q_k$, then we can find $q_{n_0}=\bQ_k, q_{n_1}, q_{n_2},\ldots$ that $(q_{n_j}, q_{n_{j+1}})$ is $\CD(\cal A, \cal A, \cal A^3)$ bridge for $j=0, 1, \ldots.$ In this case, we let $Q_{k+j}=q_{n_j}$ for $j\geq 1$. Otherwise, we let $q_n>Q_k$ be the smallest denominator that satisfies $q_{n+1}>q_n^{\cal A}$. In this situation, if $q_n\leq\bQ_k^{\cal A^4}$, then let $Q_{k+1}=q_{n}$. Otherwise, we can find $q_{n_0}=\bQ_k< q_{n_1}< \ldots$ that $(q_{n_j}, q_{n_{j+1}})$ is $\CD(\cal A, \cal A, \cal A^2)$ bridge.  Let $q_{n_l}<q_{n}$ be the smallest denominator that $q_{n_l}^{\cal A^2}\geq q_n$, by which we know that $l\geq 1$. If $q_{n }\geq q_{n_{l}}^{\cal A}$, let  $Q_{k+j}=q_{n_j}$ for $j=1, \ldots, l$ and $Q_{k+l+1}=q_n$. Then $(\bQ_k, Q_{k+1})$, $(Q_{k+j}, Q_{k+j+1})$ for $j=1, \ldots, l$ are $\CD(\cal A, \cal A, \cal A^2)$ bridges and thus $\CD(\cal A, \cal A, \cal A^3)$ bridges. If $q_{n }<q_{n_l}^{\cal A}$ (in this case, $l\geq 2$), we let $Q_{k+j}=q_{n_j}$ for $j=1, \ldots, l-1$ and $Q_{k+l}=q_n$. Then $(\bQ_k, Q_{k+1})$, $(Q_{k+j}, Q_{k+j+1})$ for $j=1, \ldots, l-2$ are $\CD(\cal A, \cal A, \cal A^2)$ bridges and $(Q_{k+l-1}, Q_{k+l})$ is $\CD(\cal A, \cal A, \cal A^3)$ bridge. This completes the proof.
\end{pf}

\begin{Remark}\label{remark-selection-qk}
By the selection of $\{Q_k\}_{k\in\N_0}$, those $q_n$ satisfying $q_{n+1}>q_n^{\cal A}$ belong to $\{Q_k\}_{k\in\N_0}$. 
\end{Remark}

\begin{Corollary}[{\cite{KWYZ18}}]\label{cd-coroll}
For the subsequence $\{Q_k\}_{k\in\N_0}$ selected in Lemma \ref{CDbridge} we have $Q_k\geq
Q_{k-1}^\cal A$ for every $k\geq 1$. 
\end{Corollary}

\subsection{Lyapunov exponent and uniformly hyperbolic}\label{SecLya}
Given $A\in C^\omega(\T, SL(2,\R))$ and $\alpha\in\R$ irrational, the iterates of the quasi-periodic cocycle $(\alpha, A)$ are of the form $(\alpha, A)^n=(n\alpha, A_n)$, where 
\[A_n(x):=\left\{\begin{array}{ll} 
A(x+(n-1)\alpha)A(x+(n-2)\alpha)\cdots A(x),\ \ & n\geq 0\\
A(x+n\alpha)^{-1}A(x+(n+1)\alpha)^{-1}\cdots A(x-\alpha)^{-1}, \ \ & n<0
 \end{array}.\right.\]
 The \textit{Lyapunov exponent} is defined as $L(\alpha, A):=\lim_{n\rightarrow\infty}\frac{1}{n}\int_{\T^1}\ln \|A_n(x)\|dx\geq 0$.
 
 We say the cocycle $(\alpha, A)$ is \textit{uniformly hyperbolic} if there exists a continuous splitting $\R^2=E^s(x)\oplus E^u(x)$, and $C>0, 0<\lambda<1$ such that for every $n\geq 1$ we have 
 \[\|A_n(x)\cdot w\|\leq C \lambda^n\|w\|, \ \ \ \ w\in E^s(x),\]
 \[\|A_{-n}(x)\cdot w\|\leq C \lambda^n\|w\|, \ \ \ w\in E^u(x).\]
Such a splitting is automatically unique and thus invariant, i.e., $A(x)E^s(x)=E^s(x+\alpha)$ and $A(x)E^u(x)=E^u(x+\alpha)$.

\subsection{Global theory of one-frequency $SL(2,\R)$-cocycles.}\label{global-th}

Avila's global theory classified the analytic $SL(2,\R)$ cocycles that are not uniformly hyperbolic to three cases: \textit{supercritical}, \textit{subcritical} and \textit{critical}. A cocycle $(\alpha, A)$ which is not uniformly hyperbolic is said to be \textit{supercritical} if $L(\alpha, A)>0$; it is \textit{subcritical}, if there exists $\delta>0$ such that $L(\alpha, A(z))=0$ for $|\Im z|\leq \delta$; and it is \textit{critical} otherwise.
In particular, we say the Schr\"odinger operator $H_{v,\alpha}$ is \textit{acritical} if for any $E$ in the spectrum of the cocycle $(\alpha, S_{v,E})$ is not critical. To give an example of subcritical cocycles, we consider a cocycle $(\alpha, R_g)$ which is homotopic to constant. It is subcritical and the projective action is of the form
\[ T: \T^1\times \R\mathbb P^1\circlearrowleft  \ \textrm{via}\ (x,\varphi)\mapsto (x+\alpha, \varphi+g(x)).\]
Moreover, if $(\alpha, A)$ is supercritical, then it is \textit{non-uniformly hyperbolic}, that is, there is a splitting as in the uniformly hyperbolic case, except that this splitting is not continuous.
It is equivalent to the existence of a \textit{strange non-chaotic attractor} $\phi^-$ ( An SNA in a qpf system $T$ is a $T$-invariant graph which has negative Lyapunov exponent and is not continuous.) and a \textit{strange non-chaotic repeller} $\phi^+$ ( An SNR is a non-continuous $T$-invariant graph with positive Lyapunov exponent.) for the qpf circle homeomorphism $T_{(\alpha, A)}$ induced by the projective action of $(\alpha, A)$ (For a detailed discussion of this relation, see \cite[Section 1.3.2]{J09} ).  

The following theorem of Avila shows that the critical case is very rare.
\begin{Theorem}[\cite{A15}]\label{thm-avila-global}
Let $\alpha$ be irrational. Then for a (measure-theoretically) typical $v\in C^\omega(\T^1, \R)$, the operator $H_{v,\alpha}$ is acritical. 
\end{Theorem}


%
%

\subsection{Measure complexity of a t.d.s.}

Let $(X,T)$ be a topological dynamical system(t.d.s.) with a metric $d$ and let $\mathcal M(X,T)$ be the set of all $T$-invariant Borel probability measures on $X$. For $\rho\in\mathcal M(X,T)$ and any $n\in \N$, we consider the metric
\[\bar d_n(x,y)=\frac{1}{n}\sum_{i=0}^{n-1}d(T^i x, T^i y)\]
for any $x,y\in X$. For $\epsilon>0$, let
\[S_n(d,\rho,\epsilon)=\min\{m\in\N\ :\ \exists x_1, x_2,\cdots, x_m\ s.t.\ \rho(\bigcup_{i=1}^m B_{\bar d_n}(x_i, \epsilon))>1-\epsilon\},\]
where $B_{\bar d_n}(x,\epsilon):=\{y\in X\ :\ \bar d_n(x,y)<\epsilon \}$ for any $x\in X$.

\begin{Definition}[\cite{HWY17}]
Let $U\ :\ \N\rightarrow [1,+\infty)$ be an increasing sequence with $\lim_{n\rightarrow +\infty}U(n)=+\infty$. We say the measure complexity of $(X,d,T,\rho)$ is weaker than $U(n)$ if
\[\liminf_{n\rightarrow+\infty}\frac{S_n(d,\rho,\epsilon)}{U(n)}=0,\ \ \forall \epsilon>0.\]
\end{Definition}

\begin{Remark}
By the following proposition, the measure complexity of $(X,d,T,\rho)$ is weaker than $U(n)$ if and only if the measure complexity of $(X,d',T,\rho)$ is weaker than $U(n)$ for any compatible metric $d'$ on $X$. Thus we can simply say the measure complexity of $(X,T,\rho)$ is weaker than $U(n)$.
\end{Remark}

\begin{Proposition}[\cite{HWY17}]\label{prop-hwy}
Assume $(X,\mathcal{B}_X, T, \rho)$ is measurably isomorphic\footnote{We say $(X,\mathcal B_X, T,\rho)$ is \textit{measurably isomorphic} to $(Y,\mathcal B_Y, S, \nu)$, if there are $X'\in\mathcal B_X, Y'\in\mathcal B_Y$ with $\rho(X')=1, \nu(Y')=1$, $TX'\subseteq X'$, $SY'\subseteq Y'$, and an invertible measure-preserving map $H: X'\rightarrow Y'$ with $H\circ T(x)=S\circ H(x)$ for all $x\in X'$. } to $(Y, \mathcal{B}_Y, S, \mu)$, and $U(n): \N\rightarrow [1,+\infty)$ is an increasing sequence with $\lim_{n\rightarrow +\infty}U(n)=+\infty$. Then the measure complexity of $(X,\mathcal{B}_X, T, \rho)$ is weaker than $U(n)$ if and only if the measure complexity of $(Y, \mathcal{B}_Y, S, \mu)$ is weaker than $U(n)$.
\end{Proposition}


\begin{Definition}[\cite{HWY17}]\label{sub-polynomial}
Let $(X,T)$ be a t.d.s. and $\rho\in\mathcal M(X,T)$. The measure complexity of $(X,T,\rho)$ is sub-polynomial if it is weaker than $U_\tau(n)=n^\tau$ for any $\tau>0$.
\end{Definition}

Using the measure complexity, Huang-Wang-Ye\cite{HWY17} provided a criterion for a t.d.s. satisfying the required disjointness, which says:
\begin{Theorem}[\cite{HWY17}]\label{thm-hwy-1}
Let $(X,T)$ be a t.d.s. such that the measure complexity of $(X,T,\rho)$ is  sub-polynomial for any $\rho\in\mathcal M(X,T)$. Then the M\"{o}bius disjointness conjecture holds.
\end{Theorem}

Using Theorem \ref{thm-hwy-1}, \cite{HWY17} got M\"obius disjointness for systems with discrete spectrum. Recall that a t.d.s. $(X,T,\rho)$ has \textit{discrete spectrum} if $L^2(X,\mathcal B_X, \rho)$ is spanned by the set of eigenfunctions for $T$, where $\mathcal B_X$ stands for the Borel $\sigma$-algebra of $X$ and $\rho\in\mathcal M(X,T)$. Let $(X,T)$ and $(Y,S)$ be two t.d.s., and let $\rho\in\mathcal M(X,T)$ and $\nu\in\mathcal M(Y,S)$. If $(X,\mathcal B_X, T,\rho)$ is measurably isomorphic to $(Y, \mathcal B_Y, S, \mu)$, then $(X, T, \rho)$ has discrete spectrum iff $(Y, S, \nu)$ has discrete spectrum. Moreover, the following theorem implies that a t.d.s which has discrete spectrum for any invariant measure is linear disjoint with the M\"obius function.

\begin{Theorem}[\cite{HWY17}]\label{thm-hwy-3}
Let $(X,T)$ be a t.d.s. and $\rho\in\mathcal M(X,T)$. If $\rho$ has discrete spectrum, then the measure complexity of $(X,T,\rho)$ is sub-polynomial.
\end{Theorem}

The following fact is evident and we omit the proof.
\begin{Proposition}\label{pro-sarnak-conjugate}
Let $(X,T), (Y, S)$ be two t.d.s.. If they are conjugate, i.e.  there exists a homeomorphism $h: X\rightarrow Y$ such that $h\circ T=S\circ h$,  then the M\"obius disjointness conjecture holds for $(X,T)$ if and only if it is true for $(Y,S)$.
\end{Proposition}



\section{Proof of the main theorems}\label{sec-3}
In this section, we give the proof of Theorem  \ref{mainthm-1} and Corollary \ref{mainthm-schrodinger},  Corollary \ref{mainthm-schrodinger-1}. 
In order to prove Theorem \ref{mainthm-1}, we first consider the case that $(\alpha, A)$ is not uniformly hyperbolic, since the case that $(\alpha, A)$ is uniformly hyperbolic is relatively easy (Lemma \ref{lem-le-positive}).

\begin{Theorem}\label{mainthm}
Let $(\alpha, A)\in \R\backslash\Q\times C^\omega(\T^1, SL(2,\R))$. If $(\alpha, A)$ is almost reducible and not uniformly hyperbolic, then Sarnak's M\"obius conjecture holds for $T_{(\alpha, A)}$.
\end{Theorem}

\textit{Strategy of the proof.}  
As we mentioned before, quantitative estimates are not involved in the definition of almost reducibility, and the key point in the proof is to get the quantitative almost reducibility of the cocycle with delicate control on the conjugations and perturbations. To get the appropriate quantitative almost reducibility, if $\beta(\alpha)=0$, the proof will be based on   a quantitative version of Aubry duality and almost localization, and if $\beta(\alpha)>0$, we will perform a KAM scheme.  In the end, we deduce the desired result from the quantitative almost reducibility. 
%
%

\subsection{Quantitative almost reducibility} 
First we state the quantitative almost reducibility results that we need, and leave the full proof to  section \ref{sec-almost-reducibility-1} and  section \ref{sec-almost-reducibility-2}. Both the quantitative almost reducibility result and the proof
will depend on the arithmetic property of the frequency $\alpha$. If  $\beta(\alpha)=0$, then we have the following:

\begin{Proposition}\label{prop-quantitative-almost-reducibility-1}
Let $\alpha\in\R\backslash\Q$ with $\beta(\alpha)=0$. If $(\alpha, A)$ is almost reducible and not uniformly hyperbolic, then either it is reducible, or there exist $r_*>0, n_*\in\N,  \varrho\in\R$, and an infinite sequence $\{n_j\}_{j\in\N}\subseteq \N$ such that for $n_j\geq n_*$, there exists $W_j: \T^1\rightarrow  SL(2,\R)$ analytic with $\|W_j\|_{r_*}\leq Ce^{o(n_j)}$ that $$W_j(\cdot+\alpha)^{-1}A(\cdot)W_j(\cdot)=R_{\pm\varrho}(I+G_j(\cdot)),$$  with $\|G_j\|_{r_*}\leq Ce^{-cn_j}$, where  $c, C$ are constants not depending on $j$.
\end{Proposition}

If  $\beta(\alpha)>0$,    let  $\mathcal A>2$ and  $\{Q_k\}_{k\in\N_0}$ be the subsequence selected in Lemma \ref{CDbridge}. Then we will have the following: 

\begin{Proposition}\label{prop-quantitative-almost-reducibility-2}
Let $\alpha\in\R\backslash\Q$ with $\beta(\alpha)>0$. If $(\alpha, A)$ is almost reducible and not uniformly hyperbolic, then
\begin{enumerate}
\item either there exist sequences  $\{W_{j}\}_{j\in\N}\subseteq C^\omega(\T^1, PSL(2,\R))$, $\{\varrho_j\}_{j\in\N}\subseteq\R$ and $\{n_j\}_{j\in\N} \subseteq \N$ 
such that $$W_{j}(\cdot + \alpha)^{-1}A(\cdot)W_{j}(\cdot)=R_{\varrho_j}(I+G_{j}(\cdot))$$ with $\|W_{j}\|_{C^1}\leq e^{n_j\eta_j}$, $\|G_{j}\|_{C^0}\leq e^{-n_j\tau_j}$, where $\eta_j=o(\tau_j)<1$, and $\tau_jn_j^{1/2}>1$;
\item or there is $k_*\in \N$ that for $k\geq k_*$, there exists $W_k\in C^\omega(\T^1, PSL(2,\R))$ that $$W_{k}(\cdot+\alpha)^{-1}A(\cdot)W_k(\cdot)=\left(\begin{array}{cc}1 & c_k \\ 0 & 1\end{array}\right)(I+G_k(\cdot))$$ with estimates $\|G_k\|_{C^0}\leq Ce^{-Q_k\tau_k}$, $\|W_k\|_{C^1}\leq Ce^{CQ_k\eta_k}$,  where $e^{-\frac{Q_k\tau_k}{10}}\leq |c_k|\leq e^{-Q_k\eta_k}$, $\eta_k=o(\tau_k), \tau_k=o(1)$, $\tau_kQ_k^{\frac{1}{2 \mathcal A}}>1$ and $C$ is a global constant not depending on $k$. 
\end{enumerate}
\end{Proposition}

\subsection{Proof of Theorem \ref{mainthm}} Let $(\alpha, A)$ be almost reducible and not uniformly hyperbolic. By Proposition \ref{prop-quantitative-almost-reducibility-1} and Proposition \ref{prop-quantitative-almost-reducibility-2}, the cocycle $(\alpha, A)$ is either reducible, or almost reducible with corresponding quantitative estimations. 

For reducible cocycles, we have the following measure complexity estimate:
\begin{Lemma}\label{lem-reducible-disjointness}
If  $(\alpha, A)$ is reducible and not uniformly hyperbolic, then the M\"obius disjointness holds for the t.d.s. $(\T^1\times\R\mathbb P^1, T_{(\alpha, A)})$ induced by the projective action of $(\alpha, A)$.\end{Lemma}
\begin{pf}
Assume that $(\alpha, A)$ is conjugate to $(\alpha, D)$ with $D\in SL(2,\R)$. By Proposition \ref{pro-sarnak-conjugate}, it suffices to show the M\"obius disjointness conjecture holds for $(\alpha, D)$. As $(\alpha, A)$ is not uniformly hyperbolic, then $D$ is elliptic or parabolic.
Without loss of generality, we assume $D$ is in the real normal form 
 $R_\varrho$ 
 or $\left(\begin{array}{cc} 1 & c\cr 0 & 1 \end{array}\right)$, where $\varrho\in\T^1$, and $0\neq c\in\R$.

Case (1): $D=R_\varrho$.
In this case, 
\[T_{(\alpha, D)}^n(x,\varphi)=(x+n\alpha, \varphi+n\varrho).\]
As $T_{(\alpha, D)}$ is a translation on $\T^1\times\R\mathbb P^1$, it is a straightforward corollary of Davenport's theorem \cite{D37} that it satisfies the M\"obius disjointness conjecture (see e.g \cite{S09}).

Case (2): $D=\left(\begin{array}{cc} 1 & c\cr 0 & 1 \end{array}\right)$. In this case, we have 
\[T_{(\alpha, D)}^n(x,\varphi)=(x+n\alpha, \Pi\circ\left( \frac{1}{2\pi} \arctan \frac{1}{\cot\hat\varphi+nc} \right) ), \]
where $\Pi: (-\frac{1}{4}, \frac{1}{4}]\rightarrow \R\mathbb P^1$ is the canonical projection, and $\hat\varphi:=2\pi \gamma(\varphi)$ with $\gamma: \R\mathbb P^1\rightarrow (-\frac{1}{4}, \frac{1}{4}]$ being the lift of the identity map on $\R\mathbb P^1$.
Then $\Pi\circ\left( \frac{1}{2\pi} \arctan \frac{1}{\cot\hat\varphi+nc} \right)\to 0$ as $n\to \infty$. Given $f\in C(\T^1\times \R\mathbb P^1)$ and $(x,\varphi)\in \T^1\times \R\mathbb P^1$, as $f$ is uniformly continuous, $f(T_{(\alpha, D)}^n(x,\varphi))-f(x+n\alpha,0)\to 0$. It follows that $\lim_{n\to\infty}\frac{1}{N}\sum_{n=1}^N\mu(n)f(T_{(\alpha, D)}^n(x,\varphi))$ coincides with $\lim_{n\to\infty}\frac{1}{N}\sum_{n=1}^N\mu(n)f(x+n\alpha,0)$. Again, one can easily deduce from Davenports theorem as in Case (1) that this later limit vanishes. Hence the M\"obius disjointness holds in this case.
\end{pf}

Lemma \ref{lem-reducible-disjointness} deals with the case that  $(\alpha, A)$ is reducible. If $(\alpha, A)$ is almost reducible with suitable estimates, 
that is 
$$W_j(\cdot +\alpha)^{-1}A(\cdot)W_j(\cdot)=A_j(I+G_j(\cdot)),$$ 
then the proof of M\"obius disjointness conjecture not only depends on the quantitative estimates of the conjugation $W_j$ and the perturbation $G_j$, but most importantly, also depends on the structure of  $A_j$. If $A_j$ is elliptic, then one should estimate  the measure complexity of the dynamics as follows:

\begin{Proposition}\label{prop-measure-complexity-almost}
Let $\alpha\in\R\backslash\Q$ and $A\in C^\omega(\T^1, SL(2,\R))$. If there exists $j_*\in\N$ and a sequence  $\{\varrho_j\}_{j\in\N}\subseteq\R$ such that for $j\geq j_*$, there is $W_j\in  C^\omega(\T^1, PSL(2,\R))$  that 
\[W_j(\cdot+\alpha)^{-1}A(\cdot)W_j(\cdot)=R_{\varrho_j}(I+G_j(\cdot)),\] with $\|G_j\|_{C^0}\rightarrow 0$ and $\|W_j\|_{C^1}\|G_j\|_{C^0}^\eta\rightarrow 0$ as $j\rightarrow \infty$ for any $\eta>0$, then the measure complexity of $(\T^1\times\R\mathbb P^1,T_{(\alpha, A)},\rho)$ is sub-polynomial for any $\rho\in \mathcal M(\T^1\times\R\mathbb P^1, T_{(\alpha, A)})$.
\end{Proposition}

If $A_j$ is parabolic and $\beta(\alpha)>0$, then one should use periodic approximation, and decompose the periodic sequence into short average of Dirichlet characters, reducing the problem to control the average of multiplicative function on a typical interval to finish the proof.

For $\alpha\in\R\backslash\Q$ with $\beta:=\beta(\alpha)>0$, we assume $\{\frac{p_k}{q_k}\}_{k\in\N_0}$ is the convergents of $\alpha$, and denote 
\[\mathfrak Q:=\{q_k \ | \ q_{k+1}\geq e^{\frac{\beta}{2}q_k}\}.\]

\begin{Proposition}\label{prop-parabolic}
For the analytic cocycle $(\alpha, A)$ with $\beta:=\beta(\alpha)>0$, if there exists a sequence $\{q_{k_j}\}_j\subseteq \mathfrak Q$ and $\{W_j\}_j\subseteq C^\omega(\T^1, PSL(2,\R))$ that 
$$W_j(\cdot +\alpha)^{-1}A(\cdot)W_j(\cdot)=A_j(I+G_j(\cdot))$$ with $A_j=\left(\begin{array}{cc} 1 & c_j\\ 0 &1 \end{array}\right)$, $\|G_{j}\|_{C^0}\leq \tilde Ce^{-q_{k_j}\tau_j}$,  $\|W_j\|_{C^1}\leq \tilde C e^{\tilde Cq_{k_j}\eta_j}$, $\tilde ce^{-(\frac{1}{7}-\xi)q_{k_j}\tau_j}\leq |c_j|\leq 1$ for some $0<\xi<\frac{1}{7}$, $\eta_j=o(\tau_j)$, $\tau_j=o(1)$, and $\tau_jq_{k_j}>q_{k_j}^{\frac{1}{2}+\epsilon}$ for some $\epsilon>0$, where $\tilde c, \tilde C>0$ are absolute constants,  then the M\"obius disjointness holds for $(\T^1\times\R\mathbb P^1, T_{(\alpha, A)})$.
\end{Proposition}

We leave the proof of Proposition \ref{prop-measure-complexity-almost} and Proposition  \ref{prop-parabolic} in section \ref{sec-measure-complexity} and  section \ref{sec-parabolic} respectively.

\medskip

\noindent\textit{Proof of Theorem \ref{mainthm}.} For $\beta(\alpha)=0$, if $(\alpha, A)$ is reducible, then by Lemma \ref{lem-reducible-disjointness}, the M\"obius disjointness conjecture holds for $(\T^1\times\R\mathbb P^1, T_{(\alpha, A)})$; otherwise, by Proposition \ref{prop-quantitative-almost-reducibility-1}, there exist $r_*>0, n_*\in\N, \varrho\in\R$, and a sequence $\{n_j\}_{j\in\N}$ such that for $n_j\geq n_*$, there exists $W_j\in C^\omega(\T^1, SL(2,\R))$ such that $W_j(\cdot+\alpha)^{-1}A(\cdot)W_j(\cdot)=R_{\pm\varrho}(I+G_j(\cdot))$ with $\|W_j\|_{C^1}\leq \frac{C}{r_*}e^{o(n_j)}$ and $\|G_j\|_{C^0}\leq Ce^{-cn_j}$, which means that $\|W_j\|_{C^1}\|G_j\|_{C^0}^\eta\leq \frac{C^2}{r_*}e^{-c\eta n_j}e^{o(n_j)}\rightarrow 0$ as $j\rightarrow \infty$ for any $\eta>0$. Then by Proposition \ref{prop-measure-complexity-almost}, the measure complexity of $(\T^1\times\R\mathbb P^1, T_{(\alpha, A)}, \rho)$ is sub-polynomial for any $\rho\in\mathcal M(\T^1\times\R\mathbb P^1, T_{(\alpha,A)})$ and hence the M\"obius disjointness holds by Theorem \ref{thm-hwy-1}.

For $\beta(\alpha)>0$, if there exist sequences $\{W_j\}_{j\in\N}\subseteq C^\omega(\T^1, PSL(2,\R))$, $\{\varrho_j\}_{j\in\N}\subseteq \R$, and $\{n_{j}\}_{j\in\N}\subseteq \N$ such that $W_{j}(\cdot+\alpha)^{-1}A(\cdot)W_j(\cdot)=R_{\varrho_j}(I+G_j(\cdot))$ with $\|W_j\|_{C^1}\leq e^{n_j\eta_j}$, $\|G_j\|_{C^0}\leq e^{-n_j\tau_j}$, $\eta_j=o(\tau_j)$ and $\tau_jn_j>n_j^{1/2}$, then by Proposition \ref{prop-measure-complexity-almost} and Theorem \ref{thm-hwy-1}, the M\"obius disjointness conjecture holds. Otherwise, by Proposition \ref{prop-quantitative-almost-reducibility-2}, there is $k_*\in \N$ that for $k\geq k_*$, there exists $W_k\in C^\omega(\T^1, PSL(2,\R))$ that $W_{k}(\cdot+\alpha)^{-1}A(\cdot)W_k(\cdot)=A_k(I+G_k(\cdot))$ with $\|G_k\|_{C^0}\leq Ce^{-Q_k\tau_k}$, $\|W_k\|_{C^1}\leq Ce^{CQ_k\eta_k}$, $A_k=\left(\begin{array}{cc}1 & c_k \\ 0 & 1\end{array}\right)$, where $e^{-\frac{Q_k\tau_k}{10}}\leq|c_k|\leq e^{-Q_k\eta_k}$, $\eta_k=o(\tau_k)<1$, $\tau_k=o(1)$, $\tau_kQ_k^{\frac{1}{2 \mathcal A}}>1$ and $C$ is a global constant not depending on $k$. Moreover, by Remark \ref{remark-selection-qk}, we have that for large enough $q_{k_j}\in \mathfrak Q$, it also belongs to $\{Q_k\}_{k\in\N_0}$. Hence, 
by Proposition \ref{prop-parabolic}, the M\"obius disjointness holds for the corresponding t.d.s. $(\T^1\times\R\mathbb P^1, T_{(\alpha,A)})$.

\qed

\subsection{Proof of Theorem \ref{mainthm-1}} \label{SecHyperbolic}

If $(\alpha, A)$ is almost reducible and not uniformly hyperbolic, then by Theorem \ref{mainthm}, the M\"obius disjointness conjecture holds.  
Now we consider the cases $(\alpha, A)$ is uniformly hyperbolic. Actually, for any cocycle $(\alpha, A)$ that has positive Lyapunov exponent, the M\"obius disjointness holds for its projective action, which is the content of the following lemma. Thus, we finish the proof.

\begin{Lemma}\label{lem-le-positive}
Let $(\alpha, A)$ be a continuous $SL(2,\R)$-cocycle. If $L(\alpha,A)>0$, then the M\"obius disjointness holds for $(\T^1\times\R\mathbb P^1, T_{(\alpha,A)})$.
\end{Lemma}
\begin{pf}

%

It follows from the Oseledets Theorem that there exist two $T_{(\alpha,A)}$-invariant graphs $\phi^{\pm}$\cite[Section 1.3.2]{J09}.
Then we can associate a measure $\mu_{\phi^\pm}$ respectively  given by 
\[\mu_{\phi^\pm}(\Omega)=\textrm{Leb}_{\T^1}(\pi_1(\Omega\cap \Phi^\pm)) \]
for every Lebesgue-measurable set $\Omega\subseteq \T^1\times\R\mathbb P^1$, where $\pi_j$ is the projection to the $j$-th variable, and $\Phi^\pm:=\{ (x, \phi^\pm(x))  : x\in\T^1\}\subseteq\T^1\times\R\mathbb P^1$. It is easy to see that $\mu_{\phi^\pm}$ are $T_{(\alpha, A)}$-invariant and ergodic. Moreover, $\mu_{\phi^\pm}$ are the exactly two $T_{(\alpha, A)}$-invariant ergodic measures, which means all the $T_{(\alpha, A)}$-invariant measures are of the form $\rho=\lambda \mu_{\phi^+}+(1-\lambda)\mu_{\phi^-}$ with $\lambda\in [0,1]$. 
Let 
\[H^\pm : \Phi^{\pm}\rightarrow \T^1, \ via\ \  (x, \phi^\pm(x)) \mapsto x.\] Then  $(\T^1\times\R\mathbb P^1, \mathcal B_{\T^1\times\R\mathbb P^1}, T_{(\alpha, A)}, \mu_{\phi^\pm})$ are measurably-isomorphic to  $(\T^1, \mathcal B_{\T^1},  R_\alpha,  \textrm{Leb}_{\T^1})$ by $H^\pm$, where $R_\alpha  :  \T^1\rightarrow\T^1,$ via $x\mapsto x+\alpha$. Since $(\T^1, R_\alpha, \textrm{Leb}_{\T^1})$ has discrete spectrum, then owing to 
Theorem \ref{thm-hwy-3}, the measure complexity of $(\T^1, R_\alpha, \textrm{Leb}_{\T^1})$ is sub-polynomial, and it is the same for 
 $(\T^1\times\R\mathbb P^1, T_{(\alpha, A)}, \mu_{\phi^\pm})$ by Proposition \ref{prop-hwy}.
Thus, the measure complexity of $(\T^1\times\R\mathbb P^1, T_{(\alpha, A)}, \rho)$ is sub-polynomial for any $\rho\in\mathcal M(\T^1\times\R\mathbb P^1, T_{(\alpha, A)})$, which implies that the M\"obius disjointness conjecture holds by Theorem \ref{thm-hwy-1}.
\end{pf}

\subsection{Proof of Corollary \ref{mainthm-schrodinger}}


By Theorem \ref{thm-avila-global}, for (measure-theoretically) typical $v\in C^\omega(\T^1,\R)$, for any $E\in \Sigma_{v,\alpha}$, where $\Sigma_{v,\alpha}$ is the spectrum of $H_{v,\alpha}$, the cocycle $(\alpha, S_{v,E})$ is either supercritical or subcritical.
Moreover, it is well known (see \cite{JM82}) that
\begin{equation}\label{equ-not-spectrum}
\Sigma_{v,\alpha}=\{E\in\R\ :\ (\alpha, S_{v,E}) \ \textrm{is not uniformly hyperbolic}\}.
\end{equation}
Therefore, if $E\notin \Sigma_{v,\alpha}$, then $(\alpha, S_{v,E})$ is uniformly hyperbolic. Thus, for (measure-theoretically) typical $v\in C^\omega(\T^1,\R)$ and for any $E\in\R$, the cocycle $(\alpha, S_{v,E})$ is {either uniformly hyperbolic, supercritical or subcritical. For the case $(\alpha, S_{v,E})$ is supercritical or uniformly hyperbolic, the M\"obius disjointness conjecture holds for $T_{(\alpha, S_{v,E})}$ by Lemma \ref{lem-le-positive}, since $L(\alpha, S_{v,E})>0$. If $(\alpha, S_{v,E})$ is subcritical, then by Avila's ART \cite{A3,A2}, it is almost reducible. Then Theorem \ref{mainthm-1} yields Corollary \ref{mainthm-schrodinger}.

\subsection{Proof of Corollary \ref{mainthm-schrodinger-1}} Denote $v(x)=\cos(2\pi x)$.  We only need to consider the case $\lambda\neq 0$. Otherwise if $\lambda=0$, then the cocycle is a constant one,  and of course it is almost reducible. 


 If $E\notin\Sigma_{\lambda v, \alpha}$, then by (\ref{equ-not-spectrum}), the cocycle $(\alpha, S_{\lambda v, E})$ is uniformly hyperbolic. And if $|\lambda|>1$, then by Corollary 2 in \cite{BJ}, 
\[L(\alpha, S_{\lambda v, E})=\max\{0, \ln |\lambda|\}>0.\]
Then these two cases follow from Lemma \ref{lem-le-positive}.

If $0<|\lambda|<1$, then by Theorem 19 in \cite{A15}, the cocycle $(\alpha, S_{\lambda v, E})$ is subcritical. We 
 divide the proof into two cases: $\beta(\alpha)=0$ and $\beta(\alpha)>0$. For $\beta(\alpha)>0$,  the cocycle $(\alpha, S_{\lambda v, E})$ is  almost reducible by Theorem 1.1 in \cite{A3}. For $\beta(\alpha)=0$,  the cocycle $(\alpha, S_{\lambda v, E})$ is  almost reducible by Theorem 3.8 in  \cite{A1}.  Then the result follows from Theorem \ref{mainthm-1}.

\section{Quantitative almost reducibility for $\beta(\alpha)=0$: proof of Proposition \ref{prop-quantitative-almost-reducibility-1}}\label{sec-almost-reducibility-1}

In the case  $\beta(\alpha)=0$, we will  use  the quantitative almost localization of the dual model to get the desired result of the original system. 

We now give some concepts and known results related to our proof.

\subsection{Aubry duality and almost localization}
Let $v: \T^1 \rightarrow \R$ be analytic and $H=H_{v,\alpha,x}: \ell^2(\Z)\rightarrow \ell^2(\Z)$ be the quasiperiodic Schr\"{o}dinger operator with $v(x)=\sum_{k\in\Z} \hat v_k e^{2\pi i kx}$. Then the dual operator $\hat H_{v,\alpha,\theta}$  given by
\[(\hat H_{v,\alpha,\theta} \hat u)_n= \sum_{k\in\Z} \hat v_k \hat u_{n-k}+ 2\cos(2\pi(\theta+n\alpha)) \hat u_n,\] 
is defined on $\ell^2(\Z)$. It has the property that if $u :\T^1\rightarrow \C$ is an $L^2$ function such that its Fourier coefficients satisfy $\hat H_{v, \alpha, \theta}\hat u=E \hat u$, then
\[S_{v,E}(x)\cdot U(x) =e^{2\pi i \theta}U(x+\alpha),\]
where $U(x)=\left(\begin{array}{cc} e^{2\pi i\theta}u(x)\cr
u(x-\alpha)
\end{array}\right)$.

Now, we  give some necessary and useful concepts. 

\begin{definition}[Resonances]\label{def-reson-eps}
Fix $\epsilon_0>0$ and $\theta\in\R$. An integer $k\in \Z$ is called an $\epsilon_0$-resonance of $\theta$ if $\|2\theta-k\alpha\|_{\T}\leq e^{-\epsilon_0 |k|}$ and $\|2\theta-k\alpha\|_{\T}=\min_{|j|\leq |k|} \|2\theta-j\alpha\|_{\T}$. 
\end{definition}
We denote by $\{n_j\}_{j\in\Z}$ the set of $\epsilon_0$-resonances of $\theta$, ordered in such a way that $|n_1|\leq |n_2|\leq \cdots$. We say that $\theta$ is $\epsilon_0$-resonant if the set $\{n_j\}_j$ is infinite.
In particular,  by direct computation, one can see that if $\beta(\alpha)=0$, then $\|2\theta-k\alpha\|_{\T}\leq e^{-\epsilon_0|k|}$ implies $\|2\theta-k\alpha\|_{\T}=\min_{|j|\leq |k|}\|2\theta-j\alpha\|_{\T}$ for $k$ large. Moreover, the following inequality holds.

\begin{Lemma}[\cite{A1}]\label{lem-est-L}
If $\beta(\alpha)=0$, then we have  
\[e^{\epsilon_0|n_j|}\leq L_j:=\|2\theta-n_j\alpha\|_{\T}^{-1}\leq e^{o(|n_{j+1}|)}.\]
\end{Lemma}

\begin{definition}[Almost localization]\label{def-almost-loca}
The family $\{\hat H_{v,\alpha,\theta}\}_{\theta\in\R}$ is almost localized if there exist constants $C_0, C_1, \epsilon_0, \epsilon_1>0$ such that for every solution $\hat u$ of $\hat H \hat u=E\hat u$ satisfying $\hat u_0=1$ and $|\hat u_k|\leq 1+|k|$,  we have
\begin{equation}\label{est-almost-loca}
|\hat u_k|\leq C_1e^{-\epsilon_1|k|},\ \ \forall \ C_0(1+|n_j|)\leq |k|\leq C_0^{-1}|n_{j+1}|,
\end{equation}
where $\{n_j\}_j$ is the set of $\epsilon_0$-resonances of $\theta$.
\end{definition}

 In order to get the almost reducibility of the Schr\"{o}dinger cocycles, we will actually need the almost localization of the dual operators, which is the following result of Avila and Jitomirskaya \cite{AJ10}:

 \begin{Theorem}[\cite{AJ10}]\label{thm-almost-loca}
 Let $\alpha\in\R\backslash \Q$ satisfy $\beta(\alpha)=0$. There exists an absolute constant $c_0>0$ such that for any given $0<r_0<1, C_0>1$, there exist $\epsilon_0=\epsilon_0(r_0)>0, \epsilon_1=\epsilon_1(r_0, C_0)\in (0,r_0)$ and $C_1=C_1(\alpha, r_0, C_0)>0$ such that the following holds: given any $v\in C^\omega(\T^1, \R)$ satisfying $\|v\|_{r_0}\leq c_0 r_0^3$, the family $\{\hat H_{v,\alpha,\theta}\}_{\theta\in\R}$ is almost localized with parameters $C_0, C_1, \epsilon_0, \epsilon_1$ as in (\ref{est-almost-loca}). For $v(x)=2\lambda cos(2\pi x)$, the conclusion holds for $0<|\lambda|<1$.
 \end{Theorem}

While almost reducibility allows one to conjugate the dynamics of cocycles close to constant ones, it is rather convenient to have the conjugated cocycles in Schr\"{o}dinger form, since many results (particularly the ones depending on Aubry duality, c.f.\cite{AJ10}) are obtained only in this setting. The following lemma obtained by Avila and Jitomirskaya \cite{AJ11} takes care of this:

\begin{Lemma}[\cite{AJ11}]\label{lem-almost-redu-schro-avila-jito}
Let $(\alpha, A)\in \R\backslash\Q\times C^\omega(\T^1, SL(2,\R))$ be almost reducible. Then there exists $h_*>0$ such that for every $\kappa>0$, there is $v\in C_{h_*}^\omega(\T^1,\R)$ with $\|v\|_{h_*}<\kappa$, $E\in\R$ and $B\in C_{h_*}^\omega(\T^1, PSL(2,\R))$ such that $$B(\cdot+\alpha)A(\cdot)B(\cdot)^{-1}=S_{v,E}(\cdot).$$
\end{Lemma}

\subsection{Proof of  Proposition \ref{prop-quantitative-almost-reducibility-1}}\label{subsec-proof-mainthm}

By Lemma \ref{lem-almost-redu-schro-avila-jito}, there exists $h_*>0$ such that for any $\kappa>0$, there is $v\in C_{h_*}^\omega(\T^1, \R)$ with $\|v\|_{h_*}<\kappa$, $E\in\R$ and $B\in C_{h_*}^\omega(\T^1, PSL(2,\R))$  such that \begin{equation}\label{schro}B(x+\alpha)A(x)B(x)^{-1}=S_{v,E}(x).\end{equation} Now, we let $r_0=\min\{ h_*, \frac{1}{2}\}$, $\kappa=c_0r_0^3$ and $c_0, \epsilon_0, \epsilon_1,C_0, C_1$ be the parameter defined as in Theorem \ref{thm-almost-loca}. Then by Theorem \ref{thm-almost-loca}, the family $\{\hat H_{v,\alpha,\theta}\}_{\theta\in\R}$ is almost localized with parameters $\epsilon_0, \epsilon_1, C_0, C_1$. Since $(\alpha, A)$ is not uniformly hyperbolic, it is the same for $(\alpha, S_{v,E})$, and then $E\in \Sigma_{v,\alpha}$. By Theorem 3.3 in \cite{AJ10}, there exist some $\theta=\theta(E)\in\R$ and $\hat u=(\hat u_k)_{k\in\Z}$ such that $\hat H_{v,\alpha,\theta}\hat u=E\hat u$ with $\hat u_0=1$ and $|\hat u_k|\leq 1, \ (\forall k\in\Z)$. Moreover,  Theorem  \ref{thm-almost-loca} implies that  $|\hat u_k|\leq C_1e^{-\epsilon_1|k|}, \ \forall \ C_0(1+|n_j|)\leq |k|\leq C_0^{-1}|n_{j+1}|$, where $\{n_j\}_j$ is the set of $\epsilon_0$-resonances of $\theta$. We fix $\theta=\theta(E)$ in the sequel. We divide the proof into two cases: non-$\epsilon_0$-resonant case and $\epsilon_0$-resonant case. 

\medskip

\textit{Case (i): $\theta$ is not $\epsilon_0$-resonant.}  

If $\theta$ is not $\epsilon_0$-resonant, then by almost localization estimate, we obtain that $\hat u$ is localized, i.e. $|\hat u_k|\leq C_1e^{-\epsilon_1|k|}$ for large enough $k$,  that is, it is the Fourier coefficients of an analytic function. Classical Aubry duality yields a connection between localization and reducibility: Indeed by  Theorem 2.5 in \cite{AJ10}, one has the following:

If $2\theta \notin \alpha\Z+\Z$, then there exists  $B: \T^1\rightarrow SL(2,\R)$ analytic such that $$B(x+\alpha)A(x)B(x)^{-1}=R_{\pm \theta}.$$
If $2\theta \in \alpha\Z+\Z$, then there exist analytic $B: \T^1\rightarrow PSL(2,\R)$ and analytic $\kappa:\T^1\rightarrow \R $ such that $$B(x+\alpha)A(x)B(x)^{-1}=\left(\begin{array}{cc}\pm1 & \kappa(x) \\ 0 & \pm1\end{array}\right).$$ Now since $\beta(\alpha)=0$, we can further conjugate the cocycle to a constant 
matrix by solving 
$\phi(x+\alpha)-\phi(x)=\kappa(x)-\int_{0}^1 \kappa(x)dx$ with $\int_{0}^1 \phi(x)dx=0$. Letting $B'(x)=\left(\begin{array}{cc}\pm1 & -\phi(x) \\ 0 & \pm1\end{array}\right)B(x)$, we have 
$$B'(x+\alpha)A(x)B'(x)^{-1}=\left(\begin{array}{cc}\pm 1 & \int_{0}^1 \kappa(x)dx \\ 0 & \pm 1\end{array}\right).$$ 
Therefore, in any case,  the cocycle $(\alpha, A)$ is reducible.
  
  \medskip
 
\textit{Case (ii): $\theta$ is $\epsilon_0$-resonant.}
 
We denote by $\{n_j\}_j$ the infinite set of $\epsilon_0$-resonances of $\theta$, and we can actually get the following local almost reducibility lemma. This lemma was presented in \cite{A1}, and we give the proof here for completeness.
 
 \begin{Lemma}\label{thm-real-conjugacy}
Given $r_0\in(0,1), C_0>1 $, let $v\in C_{r_0}^\omega(\T^1,\R)$
with $\|v\|_{r_0}<c_0r_0^3$ and take $\epsilon_0=\epsilon_0(r_0)>0$, $\epsilon_1=\epsilon_1(r_0,C_0)\in(0,r_0)$, $C_1=C_1(\alpha, r_0, C_0)>1$ as in Theorem \ref{thm-almost-loca}. 
Fix some $n=|n_j|<\infty$ and let $N=|n_{j+1}|<\infty$. Then for any $r_1\in (0, \frac{\epsilon_1}{30\pi})$, there exists $n_*=n_*(\alpha, r_0, r_1, \epsilon_0, \epsilon_1, C_0, C_1)$, 
such that if $n\geq n_*$ then there is $W: \T^1\rightarrow SL(2,\R)$ analytic with $\|W\|_{r_1/2}\leq Ce^{o(N)}$ and $$\|W(x+\alpha)^{-1}S_{v,E}(x)W(x)-R_{\sigma \theta}\|_{r_1/2}\leq Ce^{-cN},$$ where $\sigma\in \{\pm 1\}$, $C$ is a large constant and $c$ is a small constant that both depend on $\alpha, r_0, r_1, \epsilon_0, \epsilon_1, C_0, C_1$, but not on $E$ or $\theta$.
\end{Lemma}
\begin{Remark}
Same as in Theorem \ref{thm-almost-loca}, if $v=2\lambda \cos(2\pi x)$, then the conclusion holds for $0<|\lambda|<1$.
\end{Remark} 

\medskip

Then in the $\epsilon_0$-resonant case, by Lemma \ref{thm-real-conjugacy}, 
 for any $r_1\in (0, \frac{\epsilon_1}{30\pi})$ there exists $n_*$ only depending on $\alpha, r_0, r_1, \epsilon_0, \epsilon_1, C_0, C_1$, that for any $|n_j|\geq n_*$ there exists $\tilde W_j\in C_{r_1/2}^\omega(\T^1, SL(2,\R))$ with $\|\tilde W_j\|_{r_1/2}\leq Ce^{o(|n_{j+1}|)}$ and $\|\tilde W_j(x+\alpha)^{-1}S_{v,E}(x)\tilde W_j(x)- R_{\sigma \theta}\|_{r_1/2}  \leq Ce^{-c |n_{j+1}|}$. Let $W_j(x)=B(x)^{-1}\tilde W_j(x)$, where $B(x)$ is defined in \eqref{schro}. Then we have $\|W_j\|_{r_1/2}\leq \|B\|_{r_1/2}\cdot\| \tilde W_j\|_{r_1/2}<Ce^{o(|n_{j+1}|)}$ for $|n_{j+1}|$ large enough. This completes the proof of Proposition \ref{prop-quantitative-almost-reducibility-1} as long as the proof of Lemma \ref{thm-real-conjugacy} is finished. 

\subsection{Proof of Lemma \ref{thm-real-conjugacy}}\label{sec-quantitative-almost-reducibility}

In this subsection, all the constants may depend on $\alpha, r_0, r_1, \epsilon_0,$  $ \epsilon_1, C_0, C_1$, but not on $E$, $\theta$, $n$ or $N$. In the following $C, c$ represent big and small constant respectively, and $n_*\in \N$ is also a constant.

%
%
%
%
%
%
%
%
%

First by almost localization result of the dual operator, we have the following:

\begin{Lemma}\label{lem-almost-loca}
For any $m\in [4C_0 (1+n), C_0^{-1}N]$, let $I=\left[\ -\left[\frac{m}{2}\right],  m-\left[\frac{m}{2}\right]\ \right]$. Then 
\begin{equation}\label{equ-U-I}
S_{v,E}(x)U^{I}(x)=e^{2\pi i \theta} U^{I}(x+\alpha)+ e^{2\pi i \theta}\left(\begin{array}{cc} h(x)\cr 0 \end{array}\right),
\end{equation}
with $\|h\|_{r_1}\leq Ce^{-\frac{\epsilon_1}{3}m}$, where $u^{I}(x)=\sum_{k\in I}\hat u_ke^{2\pi i kx}$, $U^I(x)=\left(\begin{array}{cc} e^{2\pi i \theta}u^{I}(x)\cr  u^{I}(x-\alpha)
\end{array}\right)$. 
Moreover, we have \begin{equation}\label{equ-est-upper-U}
\|U^I\|_{r'}\leq C(\epsilon_1, C_1,r')e^{2\pi C_0 r' n},
\end{equation}
for any $0<r'\leq r_1$.
\end{Lemma}
\begin{pf}
Since $\hat H_{v,\alpha, \theta} \hat u=E \hat u $, a direct computation shows that (\ref{equ-U-I}) holds, where the Fourier coefficients $(\hat h_k)_{k\in\Z}$ of $h$ satisfy
\begin{equation}\label{equ-h-fourier-1}
\hat h_k=\chi_{I}(k)( \ E- 2\cos( 2\pi(\theta+k\alpha) ) \ ) \hat u_k-\sum_{l\in\Z}\chi_{I}(k-l)\hat u_{k-l}\hat v_l,
\end{equation}
and $\chi_{I}$ is the characteristic function of $I$. Since $\hat H_{v,\alpha,\theta}\hat u=E\hat u$, we also have 
\begin{equation}\label{equ-h-fourier-2}
\hat h_k=-\chi_{\Z\backslash I}(k) (\ E- 2\cos( 2\pi (\theta+k\alpha)  )\ )\hat u_k + \sum_{l\in\Z} \chi_{\Z\backslash I}(k-l)\hat u_{k-l}\hat v_l.
\end{equation}
Now we claim that  for all $k\in\Z$, 
\[|\hat h_k|\leq C(C_1, r_0, \epsilon_1)e^{-\frac{\epsilon_1}{3}|k|}e^{-\frac{\epsilon_1}{3}m}.\]
Recall that by Theorem \ref{thm-almost-loca}, for $\frac{m}{4}\leq |k|\leq m$ we have $|\hat u_k|\leq C_1 e^{-\epsilon_1|k|}$, and $|\hat u_k|\leq 1$ for all $k\in\Z$. Moreover, the Fourier coefficients of $v$ satisfy $|\hat v_k|\leq \|v\|_{r_0}e^{-2\pi r_0|k|}\leq c_0r_0^3 e^{-2\pi r_0 |k|}$. 

For $k\in I$, we have $|k|\leq \frac{m}{2}+1$. Then by (\ref{equ-h-fourier-2}),
\begin{eqnarray*}
|\hat h_k|&=&|\sum_{l\in\Z} \chi_{\Z\backslash I}(k-l)  \hat u_{k-l}\hat v_l|
\leq  \sum_{|k-l|\geq \frac{m}{2}} |\hat u_{k-l}\hat v_l|\\
&\leq & \sum_{\frac{m}{2}\leq |k-l|\leq m} C_1 c_0r_0^3 e^{-\epsilon_1|k-l|}e^{-2\pi r_0|l|}+ \sum_{|l|>\frac{m}{2}-1}c_0r_0^3e^{-2\pi r_0|l|}\\
&\leq & C_1 c_0r_0^3 e^{-\frac{m}{2}\epsilon_1} \sum_{l\in\Z}e^{-2\pi r_0 |l|}+c_0r_0^3 \sum_{|l|>\frac{m}{2}-1}e^{-2\pi r_0 |l|}
\\
&\leq& C(C_1, r_0) e^{-\frac{\epsilon_1}{2}m}
\\
&\leq &C(C_1,r_0,\epsilon_1) e^{-\frac{\epsilon_1}{3}m} e^{-\frac{\epsilon_1}{3}|k|}.
\end{eqnarray*}

If $k\notin I$, then $|k|\geq \frac{m}{2}$. By (\ref{equ-h-fourier-1}),
\begin{eqnarray*}
|\hat h_k|&=&|\sum_{l\in\Z} \chi_{I}(k-l) \hat u_{k-l} \hat v_l |
\leq \sum_{|k-l|\leq \frac{m}{2}+1} |\hat u_{k-l}\hat v_l |\\
&=&\left(\sum_{|k-l|<\frac{m}{4}}+ \sum_{\frac{m}{4}\leq |k-l|\leq \frac{m}{2}+1}\right) |\hat u_{k-l}\hat v_l |
\\
&\leq & \sum_{|l|>\frac{|k|}{2}}c_0r_0^3 e^{-2\pi r_0 |l|} +\sum_{l\in\Z} C_1 c_0r_0^3 e^{-\epsilon_1|k-l|}e^{-2\pi r_0 |l|}\\
&\leq & C(r_0,  C_1, \epsilon_1)e^{-\epsilon_1|k| }\\
&\leq& C(r_0,  C_1, \epsilon_1)e^{-\frac{\epsilon_1}{3}m }e^{-\frac{\epsilon_1}{3}|k|}.
\end{eqnarray*}
Therefore, 
\begin{eqnarray*}
\|h\|_{r_1}&\leq& \sum_{k\in\Z}|\hat h_k| e^{2\pi r_1|k|}
\leq C e^{-\frac{\epsilon_1}{3}m}\sum_{k\in\Z}e^{-(\frac{\epsilon_1}{3}-2\pi r_1)|k|}\leq Ce^{-\frac{\epsilon_1}{3}m}.
\end{eqnarray*} 
Now we can give the upper bound of $\|U^I(x)\|$ for $|\Im x|\leq r'$: $\forall 0<r'\leq r_1$
\begin{eqnarray*}
 \|u^I\|_{r'}&\leq& \sum_{k\in I} |\hat u_k|e^{2\pi r' |k|}\leq \sum_{|k|<C_0 (1+n)}e^{2\pi r' |k|} +\sum_{C_0(1+ n)\leq |k|\leq \frac{m}{2}+1} C_1 e^{-\epsilon_1 |k|}e^{2\pi r'|k|}\\
&\leq&\frac{2e^{2\pi r'}}{e^{2\pi r'}-1} e^{2\pi r' C_0(1+ n)}+\frac{2C_1}{1-e^{-(\epsilon_1-2\pi r')}}e^{-(\epsilon_1-2\pi r')C_0 (1+n)}\\
&\leq & C(\epsilon_1, C_1,r')e^{2\pi r' C_0 n}.
\end{eqnarray*}
%

\end{pf}

Once we have this, we can  fix $n=|n_j|, N=|n_{j+1}|$, $I_0=[\ -[\frac{C^{-1}_0}{2} N], C_0^{-1}N-[\frac{C_0^{-1}}{2}N]\ ]$. Let $U^{I_0}(x)=\left(\begin{array}{cc} e^{2\pi i \theta}u^{I_0}(x)\cr u^{I_0}(x-\alpha)\end{array}\right)$,  where $u^{I_0}(x)=\sum_{k\in I_0}\hat u_k e^{2\pi i kx}$. Let $B(x)=(  U^{I_0}(x),\ \overline{  U^{I_0}(x)})$. Recall that $L^{-1}:=L_j^{-1}=\|2\theta-n_j\alpha\|_{\T}$. Since $\theta$ is $\epsilon_0$-resonant, then 
\[e^{\epsilon_0 n}\leq L\leq e^{o(N)}\]
for $n$ large enough by Lemma \ref{lem-est-L}. 

We need the following auxiliary lemma to prove Lemma \ref{thm-real-conjugacy}, which is proved in \cite{AJ10}.

\begin{Lemma}[\cite{AJ10}]\label{lem-inf-det-B}
There exists $n_*\in\N$ 
such that if $n\geq n_*$, then 
\begin{equation*}\label{equ-est-inf-det-B}
\inf_{x\in\T^1}|\det B(x)|\geq cL^{-C_*}.
\end{equation*}
where $C_*$ is a constant that depends on $\alpha, \epsilon_0, \epsilon_1, r_0, r_1, C_0, C_1$.
\end{Lemma}

 The following  lemmas were also mentioned in \cite{AJ10}. 
 \begin{Lemma}[\cite{AJ10}]\label{lem-est-det-B-difference}
 Let $x_0\in \T^1$. Then there exists $n_*\in\N$ 
 that if $n\geq n_*$, then
 \[ \sup_{|\Im z|\leq r_1/2}|\det B(z)-\det B(x_0)|\leq e^{-cN}. \]
 \end{Lemma}

   \begin{Lemma}[\cite{AJ10}]\label{lem-est-inf-B-strip}
 There exists $n_*\in\N$ 
 that for $n\geq n_*$ we have 
 \[\inf_{|\Im z|\leq r_1/2}|\Im \det B(z)|\geq cL^{-C_*},\]
 where $C_*$ is the same as in Lemma \ref{lem-inf-det-B}.
 \end{Lemma}
 \begin{pf}
 Note that $\Re \det B(x)=0$ for all $x\in \T^1$.
 The result follows directly from Lemma \ref{lem-inf-det-B}, Lemma \ref{lem-est-det-B-difference} and Lemma \ref{lem-est-L}.
\end{pf} 
 
\noindent\textit{Proof of Lemma \ref{thm-real-conjugacy}.}  Take $S=\Re U^{I_0}(x), T=\Im U^{I_0}(x)$ on $\T^1$. 
Then we have $B(x)=(S(x), \pm T(x))\left(\begin{array}{cc}1 & 1\cr \pm i & \mp i\end{array} \right)$, and thus $\textrm{det}B(x)=\mp 2i \textrm{det}(S(x), \pm T(x))$, which implies that $\textrm{det} (S(x), \pm T(x))$ does not change the sign as $x\in\T^1$ changes by Lemma \ref{lem-inf-det-B}. Let $\sigma\in\{\pm 1\}$ be chosen such that $\textrm{det} (S, \sigma T)>0$ and we denote $\tilde W(x):=(S(x), \sigma T(x))$. Since 
\[\left(\begin{array}{cc} 1 & 1\cr \sigma i & -\sigma i \end{array}\right)\left(\begin{array}{cc} e^{2\pi i \theta} & 0 \cr 0 & e^{-2\pi i \theta} \end{array}\right) \left( \begin{array}{cc} 1 & 1\cr \sigma i & -\sigma i \end{array}\right)^{-1}=R_{-\sigma \theta}, \]
by Lemma \ref{lem-almost-loca} we have 
\[\|S_{v,E}(x)\tilde W(x)- \tilde W(x+\alpha)R_{-\sigma \theta}\|_{r_1}\leq Ce^{-\frac{\epsilon_1}{3C_0}N},\] 
and the complex extension considered here is the holomorphic one. We define $W(x)=\frac{\tilde W(x)}{|\textrm{det}B(x)/2|^{1/2}}$ on $\T^1$, so that $\textrm{det}W(x)=\frac{\textrm{det}\tilde W(x)}{|\textrm{det}B(x)/2|}=1$ (by Lemma \ref{lem-est-inf-B-strip}, there is no problem with branching when extending $|\textrm{det}B(x)|^{-1/2}$ to $|\Im x|\leq r_1/2$, since $\Re \textrm{det} B(z)=0$ for all $z\in\C$). Then  by  Lemma \ref{lem-est-inf-B-strip}, and \eqref{equ-est-upper-U} of Lemma \ref{lem-almost-loca}, we have 
\begin{eqnarray*}
\|W\|_{r_1/2} &=&\left\| \frac{\tilde W(x)}{|\textrm{det}B(x)/2|^{1/2}} \right\|_{r_1/2} \leq \frac{2\|\tilde W\|_{r_1/2}}{\inf_{|\Im x|\leq r_1/2}|\det B(x)|^{1/2}}\\
 &\leq& 
 CL^{C_*}e^{\pi C_0r_1n}\leq C L^{C_*+\frac{\pi C_0r_1}{\epsilon_0}}\leq Ce^{o(N)},
\end{eqnarray*}
and furthermore, we have
\begin{eqnarray*}
\lefteqn{\left\|S_{v,E}(x)W(x)-\left|\frac{\textrm{det}B(x+\alpha)}{\textrm{det}B(x)}\right|^{1/2}W(x+\alpha)R_{-\sigma \theta}\right\|_{\frac{r_1}{2}}}\\
&=&\left\|  \frac{S_{v,E}(x)\tilde W(x)-\tilde W(x+\alpha)R_{-\sigma\theta}}{|\textrm{det}B(x)|^{1/2}}  \right\|_{\frac{r_1}{2}}
\leq Ce^{-\frac{\epsilon_1}{3C_0}N}L^{\frac{C_*}{2}}
\end{eqnarray*}
for $n\geq n_*=n_{*}(\alpha, \epsilon_0,\epsilon_1, r_0, r_1, C_0, C_1)$. Moreover, by Lemma \ref{lem-est-det-B-difference} and \ref{lem-est-inf-B-strip}, we have 
\begin{eqnarray*}
\left\|\left| \frac{\textrm{det}B(x+\alpha)}{\textrm{det}B(x)} \right|^{1/2}-1\right\|_{\frac{r_1}{2}}
\leq\frac{\|\textrm{det}B(x+\alpha)-\textrm{det}B(x)\|_{r_1/2}}{\inf_{|\Im x|\leq r_1/2}|\textrm{det}B(x)|}
\leq Ce^{-cN}L^{C_*}
\end{eqnarray*}
for $n\geq n_{*}$. In conclusion, we obtain that 
\begin{eqnarray*}
\|S_{v,E}(x)W(x)-W(x+\alpha)R_{-\sigma \theta}\|_{r_1/2}
&\leq&Ce^{-\frac{\epsilon_1}{3C_0}N}L^{\frac{C_*}{2}}+Ce^{-cN}L^{C_*}e^{o(N)}\\
&\leq& Ce^{-cN}.
\end{eqnarray*}
That is, 
\[W(x+\alpha)^{-1}S_{v,E}(x)W(x)=R_{-\sigma\theta}(I+G(x)), \]
with \[\|G(x)\|_{r_1/2}\leq Ce^{o(N)}e^{-cN}\leq Ce^{-cN}\] 
for 
$n\geq n_{*}(\alpha, \epsilon_0, \epsilon_1, C_0, C_1, r_0, r_1)$.

\qed

\section{Quantitative almost reducibility for $\beta(\alpha)>0$: proof of Proposition \ref{prop-quantitative-almost-reducibility-2}}\label{sec-almost-reducibility-2}

For the case $\beta(\alpha)>0$, we first 
embed the local cocycle to a continuous $sl(2,\R)$ linear system \cite{YZ13} for technical convenience. Since the linear system is close to a constant system, then it allows us to perform the KAM scheme. To get the quantitative almost reducibility for continuous system, we divide the proof to several different cases in the KAM iteration, which we will explain explicitly later. Since the cocycle is almost reducible if and only if the embedded linear system is almost reducible, in the end, we obtain the quantitative almost reducibility of the cocycle.

As we will show below, any analytic almost reducible $SL(2,\R)$-cocycle that is not uniformly hyperbolic can be reduced to a cocycle in the form $(\alpha, R_\varrho e^{\tilde G(\cdot)})$ with $\|\tilde G\|$ small enough (consult section \ref{sec5.2} for details).  
 Therefore, 
we consider the reduced cocycle as a perturbation of a rotation, and deal with the
corresponding continuous quasi-periodic linear $sl(2,\R)$  system by the following  local embedding theorem.


\begin{Theorem}[\cite{YZ13}]\label{thm-embedding}
Suppose that $\alpha\in\T^d$, $h>0$, $G\in C_h^\omega(\T^d, sl(2,\R))$, $A=2\pi\varrho J$ where $\varrho\in\T$ and $J=\left(\begin{array}{cc} 0 & 1\\ -1 & 0\end{array}\right)$. Then there exist $\epsilon_*=\epsilon_*(A, h, |\alpha|)>0$, $F\in C^\omega_{\frac{h}{1+|\alpha|}}(\T^{d+1}, sl(2,\R))$ such that the cocycle $(\alpha, e^Ae^{G(\cdot)})$ is the Poincar\'e map of 
\[
\left\{ \begin{array}{ll} \dot x=(A+F( \theta))x\\
\dot{ \theta}=\omega=(1,\alpha)
\end{array} \right.
\]
with $\|F\|_{\frac{h}{1+|\alpha|}}\leq \tilde C\|G\|_h$, provided that $\|G\|_h<\epsilon_*$, where $\tilde C=\pi e^{\frac{2\pi h(2|\varrho|+1)}{1+|\alpha|}}$.  
\end{Theorem}
\begin{Remark}
By the proof in \cite{YZ13}, we know that in Theorem \ref{thm-embedding}, for $A=2\pi\varrho J$  we only need $\epsilon_*<\frac{1}{\tilde C^2}$.
\end{Remark}

For the local continuous  quasi-periodic linear $sl(2,\R)$  system, we have the  following
iterative lemma.

\subsection{Iterative Lemma}

Let $\mathcal A>2$ and $\{Q_\iota\}_{\iota\in\N_0}$ be the selected sequence in Lemma \ref{CDbridge}. 
Let  $c_0=\frac{1}{2\cdot 48^2}$, $0<h<\frac{1}{120}, h_+=\frac{h}{64}$. If  
 \begin{equation}\label{inequ-parameters}
 c_0Q_\iota^{\frac{1}{2\mathcal A}}h^2\geq \cal A^4,\ \ \  24 \ln (2Q_\iota)<Q_\iota^{\frac{1}{2\mathcal A}},
 \end{equation}
then we can get the following iterative lemma  with $C> 1$ being a large absolute constant:
\begin{Proposition}[Iterative Lemma]\label{lem-c-iterative-lemma}
We consider the continuous $sl(2,\R)$ linear differential equation
\begin{equation}\label{sys-continuous-original}
\left\{\begin{array}{ll} \dot x=(A+F(\theta))x\\
\dot\theta=\omega=(1,\alpha)
\end{array}
\right.
\end{equation}
where $\|F\|_h\leq \varepsilon\leq 2e^{-c_0 Q_\iota^{\frac{1}{\mathcal A}}h^2}$, $A=2\pi \varrho J$ or $\left(\begin{array}{cc}0 & c^*\\ 0 & 0 \end{array}\right)$ with $e^{-\frac{1}{6}Q_\iota h^2}\leq|c^*|\leq e^{-\frac{9}{4}Q_\iota h^4}$. Then there exists a conjugation map $B\in C_{h_+}^\omega(\T^2, PSL(2,\R))$ that conjugates (\ref{sys-continuous-original}) to 
\begin{equation}\label{sys-c-after-iteration}
\left\{\begin{array}{ll}
\dot x=(A_++F_+(\theta))x\\
\dot \theta=\omega=(1,\alpha)
\end{array}
\right.
\end{equation}
where 
\begin{enumerate}
\item either $A_+=2\pi\varrho_+ J$ with one of the following estimates:
\begin{enumerate}[(i)]
\item \[\ \ \ \left\{\begin{array}{ll} \|B\|_{C^1}\leq e^{C\tilde Q_{\iota+1}h_+^2},\\
\|F_+\|_{h_+}\leq e^{-  \tilde Q_{\iota+1} h_+}\leq e^{- Q_{\iota+1}^{\frac{1}{\mathcal A}}h_+}\ \ \ \end{array}\right.\] 
\item  \[\ \ \ \ \ \left\{\begin{array}{ll} \|B\|_{C^1}\leq e^{C\tilde Q_{\iota+1}h_+^4},\\
\|F_+\|_{h_+}\leq 2e^{-c_0\tilde Q_{\iota+1} h_+^2}\leq 2e^{-c_0 Q_{\iota+1}^{\frac{1}{\mathcal A}}h_+^2}\end{array}\right.\] 
\end{enumerate}
where $\tilde Q_{\iota+1}$ equals either $\bQ_\iota$ or $Q_{\iota+1}$ according to different situations, and $\tilde Q_{\iota+1}\geq Q_\iota^{\cal A}$ in both cases. 
\item or $A_+=\left(\begin{array}{cc} 0 & c_{+}^*\\
0 & 0
\end{array}\right)$ with 
\[ \left\{\begin{array}{lll}
\|B\|_{C^1}\leq e^{CQ_{\iota+1}h_+^4}\\
\|F_+\|_{h_+}\leq 9e^{-\frac{9}{4}Q_{\iota+1}h_+^2}\\
e^{-\frac{1}{6}Q_{\iota+1}h_+^2}\leq|c^*_{+}|\leq e^{-\frac{9}{4}Q_{\iota+1}h_+^4}
\end{array}.
\right.
\]
\end{enumerate}
\end{Proposition}

\begin{Remark}\label{uniformnorm}
As we can see, we have uniform control of the size of the conjugation as  $\|B\|_{C^1}\leq e^{CQ_{\iota+1}h_+^2}$.
\end{Remark}

\medskip

The iterative lemma will be proved in different situations.
We summarize different cases and state them in several propositions, leaving the proof in section \ref{sec-propositions}. We first introduce some notions. 
For any $h>0$, denote by $\mathcal B_h$ the set of $G\in C_h^\omega(\T^2, sl(2,\R))$ satisfying 
\[|G|_h:=\sum_{k\in\Z^2} \|\hat G(k)\|e^{2\pi |k|h}<\infty.\]
It is direct to check,  for any $G\in C^\omega_h(\T^2, sl(2,\R))$ and $0<h'<h$, we have 
\[|G|_{h'}\leq \frac{36}{\min\{1, (h-h')^2\}}\|G\|_h.\]

In the following propositions,  let $0<h'<\frac{1}{160}$,   $\tilde h=\frac{h'}{6}$,   and we suppose that $q_n$  is large enough such that 
\begin{equation}\label{equ-parameter}
 2c_0q_n^{\frac{1}{2\mathcal A}}h'^2\geq \cal A^4,  \ \ \ 24 \ln (2q_n)<q_n^{\frac{1}{2\mathcal A}}.
 \end{equation}
 Now according to the types of the constant matrix $A$(elliptic or parabolic), and the arithmetic property of the frequency  $\alpha$ (Diophantine or Liouvillean), we divide the full proof into four different cases: \\
 
 If $A$ is elliptic, and the frequency $\alpha$ is relatively Diophantine (i.e. $q_{n+1}\leq q_n^{\mathcal A}$), then we have the following:
 
 \begin{Proposition}\label{propo-c-ellip-cd}
 Consider the system (\ref{sys-continuous-original}) with $A=2\pi\varrho J$ and $|F|_{h'}\leq \epsilon\leq 2e^{- c_0q_n^{\frac{1}{\mathcal A}}h'^2}$. Suppose that  $q_{n+l}\leq q_n^{\mathcal A^4}$ for some $l\in\N$. Let $\Lambda_2:=\{k\in\Z^2 :  |2\varrho \pm\la k,\omega\ra|\geq \epsilon^{\frac{1}{4}}\}$. Then there exists $B\in C^\omega_{\tilde h}(\T^2, PSL(2,\R))$ that conjugates $(\ref{sys-continuous-original})$ to 
\begin{equation}\label{sys-c-ellip-cd-final}
\left\{\begin{array}{ll} \dot x=(\tilde A_1+\tilde F_1(\theta))x\\
\dot \theta=\omega=(1,\alpha)
\end{array},\right.
\end{equation}
such that
\begin{enumerate}[(a)]
\item if $\Lambda_2^c\cap\{k\in\Z^2 : |k|<\frac{q_{n+l}}{2}\}=\emptyset$, then $\tilde A_1=2\pi\bar\varrho J$ for some $\bar\varrho\in\R$, $\|B\|_{C^1}<2$, and $\|\tilde F_1\|_{\tilde h}\leq \epsilon e^{-q_{n+l}h'}$;
\item if $\Lambda_2^c\cap\{k\in\Z^2 : |k|<\frac{q_{n+l}}{2}\}\neq\emptyset$, then  we have the following:
\begin{enumerate}[(i)]
\item either $\tilde A_1=2\pi\bar\varrho  J$ for some $\bar\varrho\in\R$, with 
\begin{itemize}
\item either $\|B\|_{C^1}\leq e^{2q_{n+l}h'^2}$, $\|\tilde F_1\|_{\tilde h}\leq 4e^{-\frac{1}{2}q_{n+l}h'}$,
\item or $\|B\|_{C^1}\leq 8\pi q_{n+l}$, $\|\tilde F_1\|_{\tilde h}\leq 2e^{-\frac{c_0}{3}q_{n+l}h'^2}$;
\end{itemize}
\item or  $\tilde A_1=\left(\begin{array}{cc}0 & \bar c\\ 0 & 0\end{array}\right)$ with $\bar c\in \R$, $\|B\|_{C^1}\leq e^{q_{n+l}h'^4}$, $\|F_1\|_{\tilde h}\leq 9e^{-\frac{3}{4}q_{n+l}h'^2}$, $e^{-\frac{c_0}{3}q_{n+l}h'^2}\leq |\bar c|\leq 4e^{-\frac{3}{4}q_{n+l}h'^4}$.
\end{enumerate}
\end{enumerate}
\end{Proposition}

\medskip
 
 If $A$ is elliptic, and the frequency $\alpha$ is relatively Liouvillean (i.e. $q_{n+1}>q_n^{\mathcal A}$), then we have the following:

\begin{Proposition}\label{propo-c-ellip-liouv}
 Consider the system (\ref{sys-continuous-original}) with $A=2\pi\varrho J$ and $|F|_{h'}\leq \epsilon\leq 2e^{- c_0q_n^{\frac{1}{\mathcal A}}h'^2}$. Suppose that  $q_{n+1}>q_n^{\cal A}$. Let $\Lambda_2:=\{k\in\Z^2 :  |2\varrho \pm\la k,\omega\ra|\geq \epsilon^{\frac{1}{4}}\}$. Then there exists $B\in C^\omega_{\tilde h}(\T^2, PSL(2,\R))$ that conjugates $(\ref{sys-continuous-original})$ to 
\begin{equation}\label{sys-c-ellip-liouv-final}
\left\{\begin{array}{ll} \dot x=(\tilde A_2+\tilde F_2(\theta))x\\
\dot \theta=\omega=(1,\alpha)
\end{array},\right.
\end{equation}
such that
\begin{enumerate}[(a)]
\item if $\Lambda_2^c\cap\{k\in\Z^2 : |k|<\frac{q_{n+1}}{6}\}=\emptyset$, then $\tilde A_2=2\pi\bar\varrho J$ for some $\bar\varrho\in\R$, and $\|B\|_{C^1}< e^{q_{n+1}h'^4}$, $\|\tilde F_2\|_{\frac{h'}{2}}\leq \epsilon e^{-\frac{q_{n+1}h'}{2}}$;
\item if $\Lambda_2^c\cap\{k\in\Z^2 : |k|<\frac{q_{n+1}}{6}\}\neq\emptyset$,   then  we have the following:
\begin{enumerate}[(i)]
\item either $\tilde A_2=2\pi\bar\varrho  J$ for some $\bar\varrho\in\R$, with 
\begin{itemize}
\item either $\|B\|_{C^1}\leq e^{2q_{n+1}\tilde h^2}$, $\|\tilde F_2\|_{\tilde h}\leq 4e^{-\frac{\tilde h}{2}q_{n+1}}$,
\item or $\|B\|_{C^1}\leq e^{q_{n+1}\tilde h^4}$, $\|\tilde F_2\|_{\tilde h}\leq 2e^{-\frac{c_0}{3}q_{n+1}\tilde h^2}$;
\end{itemize}
\item or  $\tilde A_2=\left(\begin{array}{cc}0 & \bar c\\ 0 & 0\end{array}\right)$ with $\bar c\in\R$, $\|B\|_{C^1}\leq e^{q_{n+1}\tilde h^4}$, $\|\tilde F_2\|_{\tilde h}\leq 9e^{-\frac{3}{4}q_{n+1}\tilde h^2}$, $e^{-\frac{c_0}{3}q_{n+1}\tilde h^2}\leq |\bar c|\leq 4e^{-\frac{3 }{4}q_{n+1}\tilde h^4}$.
\end{enumerate}
\end{enumerate}
\end{Proposition}

\medskip
 If $A$ is parabolic, and the frequency $\alpha$ is relatively Diophantine (i.e. $q_{n+1}\leq q_n^{\mathcal A}$), then we have the following:

\begin{Proposition}\label{propo-c-parab-cd}
Consider the system (\ref{sys-continuous-original}) with $A=\left(\begin{array}{cc}0 &  c^*\\ 0 & 0\end{array}\right)$, $|c^*|\leq 1$,  and $|F|_{h'}\leq \epsilon\leq 2e^{- c_0q_n^{\frac{1}{\mathcal A}}h'^2}$. Suppose that  $q_{n+l}\leq q_n^{\mathcal A^4}$ for some $l\in\N$. Then there exists $B\in C^\omega_{\tilde h}(\T^2, PSL(2,\R))$ that conjugates $(\ref{sys-continuous-original})$ to 
\begin{equation}\label{sys-c-parab-cd-final}
\left\{\begin{array}{ll} \dot x=(\tilde A_3+\tilde F_3(\theta))x\\
\dot \theta=\omega=(1,\alpha)
\end{array},\right.
\end{equation}
such that
\begin{enumerate}[(i)]
\item either $\tilde A_3=2\pi\bar\varrho  J$ for some $\bar\varrho\in\R$, with 
\begin{itemize}
\item either $\|B\|_{C^1}\leq    4e^{ q_{n+l}  h'^2}$, $\|\tilde F_3\|_{\tilde h}\leq 4e^{-\frac{1 }{2}q_{n+l} h'}$,
\item or $\|B\|_{C^1}\leq 4$, $\|\tilde F_3\|_{\tilde h}\leq 2e^{-\frac{ c_0}{3}q_{n+l} h'^2}$;
\end{itemize}
\item or  $A_3=\left(\begin{array}{cc}0 & \bar c\\ 0 & 0\end{array}\right)$ with $\bar c\in\R$, $\|B\|_{C^1}\leq 4e^{\frac{1 }{2} q_{n+l}h'^4}$, $\|\tilde F_3\|_{\tilde h}\leq 9e^{-\frac{3 }{4}q_{n+l}h'^2}$, 
$e^{-\frac{c_0}{3}q_{n+l}h'^2}\leq |\bar c|\leq 4e^{-\frac{3}{4}q_{n+l} h'^4}$.
\end{enumerate}
\end{Proposition}

\medskip

 If $A$ is parabolic, and the frequency $\alpha$ is relatively Liouvillean (i.e. $q_{n+1}> q_n^{\mathcal A}$), then we have the following:

\begin{Proposition}\label{propo-c-parab-liouv}
Consider the system (\ref{sys-continuous-original}) with $A=\left(\begin{array}{cc}0 &  c^*\\ 0 & 0\end{array}\right)$, $|c^*|\leq e^{-\frac{9}{4}q_{n}h'^4}$,  and $|F|_{h'}\leq \epsilon\leq 2e^{- c_0q_n^{\frac{1}{\mathcal A}}h'^2}$. If $q_{n+1}>q_n^{\cal A}$, then there exists $B\in C^\omega_{\tilde h}(\T^2, PSL(2,\R))$  that conjugates $(\ref{sys-continuous-original})$ to 
\begin{equation}\label{sys-c-parab-liouv-final}
\left\{\begin{array}{ll} \dot x=(\tilde A_4+\tilde F_4(\theta))x\\
\dot \theta=\omega=(1,\alpha)
\end{array},\right.
\end{equation}
such that
\begin{enumerate}[(i)]
\item either $\tilde A_4=2\pi\bar\varrho  J$ for some $\bar\varrho\in\R$, with 
\begin{itemize}
\item either $\|B\|_{C^1}\leq  e^{2 q_{n+1}\tilde h^2 }$, $\|\tilde F_4\|_{\tilde h}\leq  4e^{-\frac{1 }{2}q_{n+1}\tilde h}$,
\item or $\|B\|_{C^1}\leq e^{ q_{n+1}\tilde h^4}$, $\|\tilde F_4\|_{\tilde h}\leq  2e^{-\frac{ c_0 }{3}q_{n+1}\tilde h^2}$;
\end{itemize}
\item or  $\tilde A_4=\left(\begin{array}{cc}0 & \bar c\\ 0 & 0\end{array}\right)$ with $\bar c\in\R$, $\|B\|_{C^1}\leq   e^{ q_{n+1}\tilde h^4   }$, $\|\tilde F_4\|_{\tilde h}\leq  9e^{-\frac{3 }{ 4}q_{n+1}\tilde h^2}$, 
$ e^{-\frac{c_0}{3}q_{n+1}\tilde h^2}\leq |\bar c|\leq  4e^{-\frac{3}{4}q_{n+1} \tilde h^4}$.
\end{enumerate}
\end{Proposition}

\medskip

\noindent\textit{Proof of Propositon \ref{lem-c-iterative-lemma}.} 
Once we have these preparing propostions, we can now finish the proof of Propositon \ref{lem-c-iterative-lemma}. 
Let $h'=\frac{3}{4}h, \tilde h=\frac{h'}{6}, h_+=\frac{\tilde h}{8}$. By the definition of $\|\cdot\|_h$ and $|\cdot|_h$, we have 
\[|F|_{ h'}\leq \frac{36}{(h/4)^2}\|F\|_h\overset{(\ref{inequ-parameters})}{\leq} \varepsilon^{\frac{4}{5}}<2e^{- c_0Q_\iota^{\frac{1}{\mathcal A}}h'^2 }.\]
By (\ref{inequ-parameters}), we have $Q_\iota$ and $h'$ satisfy (\ref{equ-parameter}) with $q_n$ replaced by $Q_\iota$. 

Now we divide the proof to several cases. \\



\noindent\textit{Case 1: $(Q_\iota, Q_{\iota+1})$ is a $CD(\mathcal A, \mathcal A, \mathcal{A}^3)$ bridge}. Then $  Q_{\iota+1}\leq Q_\iota^{\mathcal{A}^3}<Q_\iota^{\mathcal{A}^4}$. 
According to the types of $A$ (elliptic or parabolic) , we apply Proposition \ref{propo-c-ellip-cd} or Proposition \ref{propo-c-parab-cd} respectively. Then there exists $B\in C^\omega_{\tilde h}(\T^2, PSL(2,\R))$ that conjugates (\ref{sys-continuous-original}) to (\ref{sys-c-after-iteration}) with desired estimates.

\medskip

\noindent\textit{Case 2: $\bQ_\iota>Q_{\iota}^{\cal A}$}. In this case, Proposition \ref{propo-c-ellip-liouv} or Proposition \ref{propo-c-parab-liouv} is applied respectively according to whether $A$ is elliptic or parabolic. Then there exists $B_1\in C_{\tilde h}^\omega(\T^2, PSL(2,\R))$ that conjugates (\ref{sys-continuous-original}) to 
\begin{equation}\label{sys-c-after-cd}
\left\{
\begin{array}{lll}
\dot x=(\tilde A+\tilde F(\theta))x\\
\dot \theta=\omega=(1,\alpha)
\end{array},
\right.
\end{equation}
where 
\begin{itemize}
\item either $ \tilde A=2\pi\bar \varrho J$ with
\begin{itemize}
\item either $\left\{\begin{array}{ll}
\|B_1\|_{C^1}\leq  e^{2\bQ_\iota\tilde h^2}\\
\|\tilde F\|_{\tilde h}\leq 4e^{-\frac{1 }{2}\bQ_\iota\tilde h}
\end{array}\right.$,
\item or $\left\{\begin{array}{ll}
\|B_1\|_{C^1}\leq e^{\bQ_{\iota}\tilde h^4}\\
\|\tilde F\|_{\tilde h}\leq 2e^{-\frac{ c_0 }{3}\bQ_{\iota}\tilde h^2}
\end{array}\right.$;
\end{itemize}
\item or $\tilde A=\left(\begin{array}{cc} 0 & \bar c \\ 0 & 0\end{array}\right)$ with 
\[\left\{\begin{array}{lll}
\|B_1\|_{C^1}\leq  e^{\bQ_\iota \tilde h^4}\\
\|\tilde F\|_{ \tilde h}\leq 9e^{-\frac{3 }{4}\bQ_\iota\tilde h^2}\\
e^{-\frac{c_0}{3}\bQ_\iota\tilde h^2}\leq |\bar c|\leq 4e^{-\frac{3 }{4}\bQ_\iota\tilde h^4}
\end{array}.\right.\]
\end{itemize}
 
 In this case, we can check that $\bQ_\iota$ and $\frac{3}{4}\tilde h$ satisfy (\ref{equ-parameter}) with $\bQ_\iota$ and $\frac{3}{4}\tilde h$ in place of $q_n$ and $h'$ respectively. Moreover, we have 
 \[|\tilde F|_{\frac{3}{4}\tilde h}\leq \frac{36}{(\tilde h/4)^2}\|\tilde F\|_{\tilde h}\leq 2e^{-\frac{4c_0}{15}\bQ_\iota\tilde h^2}<2e^{-c_0\bQ_\iota^{\frac{1}{2\cal A}}\tilde h^2}.\]
 
 
If $\tilde A$ is elliptic and $Q_{\iota+1}\leq \bQ_\iota^{\mathcal A}$, then we let $A_+:=\tilde A$, $F_+:=\tilde F$, $B:=  B_1$, and finish one step of iteration.

If $\tilde A$ is elliptic and $Q_{\iota+1}> \bQ_\iota^{\mathcal A}$, then since $Q_{\iota+1}\leq \bQ_\iota^{\cal A^4}$ by Lemma \ref{CDbridge}, we can apply Proposition \ref{propo-c-ellip-cd}. There exists $B_2\in C_{h_+}^\omega(\T^2,PSL(2,\R))$ that conjugates 
(\ref{sys-c-after-cd}) to (\ref{sys-c-after-iteration}). Let $B:=B_1B_2$. Then we get the result.

If $\tilde A$ is parabolic, and $Q_{\iota+1}=\bQ_\iota$, then we let $B:=B_1, A_+:=\tilde A, F_+:=\tilde F$, and get the result.

Otherwise for $\tilde A$ being parabolic, since $Q_{\iota+1}\leq \bQ_\iota^{\cal A^4}$, then by Proposition \ref{propo-c-parab-cd} there exists $B_3\in C_{h_+}^\omega(\T^2, PSL(2,\R))$ that conjugates (\ref{sys-c-after-cd}) to (\ref{sys-c-after-iteration}). Let $B:=B_1B_3$. Then we finish the proof.
\qed

\subsection{Proof of Proposition \ref{prop-quantitative-almost-reducibility-2}}\label{sec5.2}

Suppose $A\in C_h^\omega(\T^1, SL(2,\R))$ for some $h>0$. 
Since $(\alpha, A)$ is $C^\omega$-almost reducible, then in the closure of its analytic conjugate class, there exists a constant cocycle $(\alpha, D)$ and the convergence is uniform in some fixed strip $\{|\Im z|<h'\}$, where $0<h'<h$.

Let $\mathcal A>2$ and $\{Q_\iota\}_{\iota\in\N}$ be the selected sequence in Lemma \ref{CDbridge}. 
Without loss of generality, we assume that $0<h'<\frac{1}{120}$. Let $h_0=\frac{h'}{1+|\alpha|}$, $c_0=\frac{1}{2\cdot 48^2}$. Let $n_*\in\N$ be large enough that $c_0q_{n_*}^{\frac{1}{2\mathcal A}}h_0^2\geq\cal A^4$,  $q_{n_*}\in\mathfrak Q \cap \{Q_\iota\}_\iota$ (this is possible by Remark \ref{remark-selection-qk}), and for any $n\geq n_*$, $q_{n}$ satisfies 
\[  24\ln (2q_n)<q_n^{\frac{1}{2\mathcal A}}.\]
Assume $q_{n_*}=Q_{\iota_*}$. Let $h_{\iota+1}=\frac{h_\iota}{64}$. Then by the selection of $\{Q_\iota\}_\iota$, we have for any $\iota\in\N_0$, 
\begin{equation*}\label{equ-parameter-selection}
c_0Q_{\iota_*+\iota}^{\frac{1}{2\mathcal A}}h_\iota^2\geq\cal A^4, \ \    24 \ln (2Q_{\iota_*+\iota}) < Q_{\iota_*+\iota}^{\frac{1}{2\mathcal A}}.
\end{equation*}

 As for $D\in SL(2,\R)$, there exists $P\in SL(2,\R)$ that $\tilde L_*:=P^{-1}DP$ is  in the normal form $R_\varrho$, $\left(\begin{array}{ll}\lambda &0\\0 &\frac{1}{\lambda} \end{array}\right)$, or  $\left(\begin{array}{ll}1 &\tilde c\\0 & 1 \end{array}\right)$, where $\lambda, \varrho, \tilde c\in\R$ depend on $D$. Then for any $0<\varepsilon<\|P\|^{-8}$, there exists $\tilde B_1\in C_{h'}^\omega(\T^1, PSL(2,\R))$ such that 
\[\tilde B_1(\cdot+\alpha)^{-1}A(\cdot)\tilde B_1(\cdot)=D+\tilde G_1(\cdot),\]
with $\|\tilde G_1\|_{h'}<\varepsilon$. Let $\tilde B_2=\tilde B_1P$. Then
\[\tilde B_2(\cdot+\alpha)^{-1}A(\cdot)\tilde B_2(\cdot)=\tilde L_*+\tilde G_2(\cdot),\]
with $\|\tilde G_2\|_{h'}=\|P^{-1}\tilde G_1P\|_{h'}<\varepsilon^{\frac{3}{4}}$. If $\tilde L_*=\left(\begin{array}{ll}1 &\tilde c\\0 & 1 \end{array}\right)$, then we let $K=\left(\begin{array}{ll}\varepsilon^{-\frac{3}{16}} &0\\0 & \varepsilon^{\frac{3}{16}} \end{array}\right)$, and get $K^{-1}\tilde L_*K=\left(\begin{array}{ll}1 &\tilde c\varepsilon^{\frac{3}{8}}\\0 & 1 \end{array}\right)$. 
Now if $\tilde L_*$ is parabolic, we let $\tilde B=\tilde B_2K$ and $\bar G=K^{-1}\tilde G_2K+\left(\begin{array}{ll}0 &\tilde c\varepsilon^{\frac{3}{8}}\\0 & 0 \end{array}\right)$; otherwise, let $\tilde B=
\tilde B_2 $ and $\bar  G=
\tilde G_2$. Then we have  
\[\tilde B(\cdot+\alpha)^{-1}A(\cdot)\tilde B(\cdot)=L_*+\bar G(\cdot),\]
with $\|\bar G\|_{h'}<C(D)\varepsilon^{\frac{3}{8}}$, and $L_*$ being a rotation or diagonal matrix.


If $L_*$ is diagonal, i.e., $L_*=\left(\begin{array}{cc} \lambda &0 \\ 0 & \lambda^{-1}\end{array}\right)$, without loss of generality, we assume that $|\lambda|>1$.   If $  |\lambda|-1 \geq \|\bar G\|_{h'}^{1/2}$, then $(\alpha, A)$ is uniformly hyperbolic by the usual cone-field criterion \cite{Yoc}, which contradicts with our hypothesis. Therefore in this case $  |\lambda|-1 \leq \|\bar G\|_{h'}^{1/2}$, and if $\lambda>1$,   we have $$L_*+\bar G=R_0(I+L_*+\bar G-I)$$ with $\|L_*+\bar G- I\|_{h'}\leq 2\|\bar G\|_{h'}^{1/2}$, and if $\lambda<-1$, we have $$L_*+\bar G=R_{\frac{1}{2}}(I-(L_*+\bar G+I) )$$ with $\|L_*+\bar G+ I\|_{h'}\leq 2\|\bar G\|_{h'}^{1/2}$.  If $L_*$ is a rotation, denoting $L_*=R_{-\varrho_0}$, then $L_*+\bar G=R_{-\varrho_0}(I+R_{\varrho_0}\bar G)$. Thus in both cases for $L_*$, by implicit function theorem, for $\|\bar G\|_{h'}$ small enough, i.e., $\varepsilon$ small enough, there exists $\tilde G\in C_{h'}^\omega(\T^1, sl(2,\R))$ such that $L_*+\bar G=R_{-\varrho_0}e^{\tilde G}$ with $\varrho_0\in\T^1$, $\|\tilde G\|_{h'}\leq 2\|\bar G\|_{h'}^{\frac{1}{2}}< \varepsilon^{\frac{1}{8}}<e^{-Q_{\iota_*}h_0^2 }$. 
%
%
%

Then by Theorem \ref{thm-embedding}, cocycle $(\alpha, R_{-\varrho_0}e^{\tilde G})$ can be embedded into the continuous linear differential equation 
\begin{equation}\label{sys-c-final-proof-original}
\left\{ \begin{array}{ll} \dot x=(\tilde A+\tilde F(\theta))x\\
\dot{\theta}=\omega=(1,\alpha)
\end{array} \right.
\end{equation}
with $\tilde A=2\pi \varrho_0 J$, $\|\tilde F\|_{h_0}< e^{-c_0Q_{\iota_*}h_0^2}$, and $\Phi^1(0,\cdot)=R_{-\varrho_0}e^{\tilde G(\cdot)}$, where $\Phi^t(\theta)$ is the fundamental solution matrix of (\ref{sys-c-final-proof-original}). Let $A_0:=\tilde A, F_0:=\tilde F$. By Proposition \ref{lem-c-iterative-lemma}, for $\iota\geq 1$ there exist  $B_\iota\in C^\omega(\T^2, PSL(2,\R))$, $ A_\iota\in sl(2,\R)$, and $ F_\iota\in C^\omega(\T^2, sl(2,\R))$, such that $B_\iota$ conjugates 
\[\left\{\begin{array}{ll}
\dot x=(A_{\iota-1}+F_{\iota-1}(\theta))x\\
\dot \theta=\omega=(1,\alpha)
\end{array}\right.\]   
to 
\begin{equation}\label{sys-c-final-proof-1}
\left\{\begin{array}{ll}
\dot x=(A_{\iota}+F_{\iota}(\theta))x\\
\dot \theta=\omega=(1,\alpha)
\end{array}\right.
\end{equation}
for either $A_\iota$ elliptic or parabolic, with corresponding estimates as in Proposition \ref{lem-c-iterative-lemma}. 
Let $\tilde W_\iota=B_1B_2\cdots B_{\iota}$. By Remark \ref{uniformnorm}, we have for any $\iota\in\N$
\begin{equation}\label{size}
\|B_{\iota}\|_{C^1}\leq e^{C Q_{\iota_*+\iota}h^2_{\iota}}.
\end{equation}

If there exists an infinite sequence $\{\iota_j\}_{j\in\N}$ such that $A_{\iota_j}$ is elliptic, then  
$\tilde W_{\iota_j}$ conjugates (\ref{sys-c-final-proof-original}) to (\ref{sys-c-final-proof-1}) with following two cases:

\smallskip
{\it Case 1:} $\|F_{\iota_j}\|_{C^0}\leq e^{-\tilde Q_{\iota_*+\iota_j}h_{\iota_j} }, \|B_{\iota_j}\|_{C^1}\leq e^{C\tilde Q_{\iota_*+\iota_j}h^2_{\iota_j}}$.
In this case, by (\ref{equ-parameter-selection}), \eqref{size} and the selection of $\{Q_\iota\}_{\iota}$, 
we have 
\[\|\tilde W_{\iota_j}\|_{C^1}\leq \iota_je^{C\tilde Q_{\iota_*+\iota_j}h^2_{\iota_j}}e^{\sum_{\iota=1}^{\iota_j-1}C Q_{\iota_*+\iota}h^2_{\iota}}<e^{2C\tilde Q_{\iota_*+\iota_j}h^2_{\iota_j}}\]
since $\tilde Q_{\iota_*+\iota_j}\geq Q_{\iota_*+\iota_j-1}^{\cal A}$, and $Q_{\iota+1}\geq Q_{\iota}^{\mathcal A}$ for any $\iota\in\N_0$.

\smallskip
{\it Case 2:} $\|F_{\iota_j}\|_{C^0}\leq 2e^{-c_0\tilde Q_{\iota_*+\iota_j}h^2_{\iota_j}}, \|B_{\iota_j}\|_{C^1}\leq e^{C\tilde Q_{\iota_*+\iota_j}h^4_{\iota_j}}$. 
Similarly as Case 1, one can estimate
\[\|\tilde W_{\iota_j}\|_{C^1}\leq \iota_je^{C\tilde Q_{\iota_*+\iota_j}h^4_{\iota_j}}e^{\sum_{\iota=1}^{\iota_j-1}C Q_{\iota_*+\iota}h^2_{\iota}}<e^{2C\tilde Q_{\iota_*+\iota_j}h^4_{\iota_j}}.\]

Otherwise, there exists $N_*\in\N$ such that for $\iota\geq N_*$,   we have that $A_\iota$ is parabolic with $\|F_\iota\|_{C^0}\leq 9 e^{-\frac{9}{4}Q_{\iota_*+\iota}h^2_{\iota}}$, $\|B_\iota\|_{C^1}\leq e^{CQ_{\iota_*+\iota}h^4_\iota}$, and $e^{-\frac{1}{6}Q_{\iota_*+\iota}h_\iota^2}\leq |c_\iota^*|\leq  e^{-\frac{9}{4}Q_{\iota_*+\iota}h_\iota^4}$. For $\iota\geq N_*$,   combined with estimate \eqref{size}, we obtain that 
\[\|\tilde W_\iota\|_{C^1}\leq \iota e^{C Q_{\iota_*+\iota}h^4_{\iota}}
e^{\sum_{k=1}^{\iota-1}C Q_{\iota_*+k}h^2_{k}}<e^{2C  Q_{\iota_*+\iota}h^4_{\iota}}.\]

Now we assume that $ \Phi_\iota^t(\theta)$ is the fundamental solution matrix of  (\ref{sys-c-final-proof-1}). Then 
\begin{equation}\label{equ-intergral}
\Phi_\iota^t(\theta)=e^{A_\iota t}\left(I+\int_{0}^t e^{-A_\iota s}F_{\iota}(\theta+s\omega)\Phi_\iota^s(\theta)ds\right).
\end{equation}
Let $\tilde G_\iota^t(\theta)=e^{-A_\iota t}\Phi_\iota^t(\theta)$ and $g_\iota(t)=\|\tilde G_\iota^t\|_{C^0}$. Then,
\[\tilde G_\iota^t(\theta)=I+\int_{0}^te^{-A_\iota s}F_\iota(\theta+\omega s)e^{A_\iota s}\tilde G_\iota^s(\theta)ds,\]
and thus $g_\iota(t)\leq 1+\int_{0}^t(1+s)^2\|F_\iota\|_{C^0}g_\iota(s)ds$ for $A_\iota$ both elliptic and parabolic with $|c_\iota^*|<1$. By Gronwall's inequality, we have 
\[g_\iota(t)\leq e^{\|F_\iota\|_{C^0}\int_{0}^t (1+s)^2ds}\leq e^{t(1+t+\frac{t^2}{3})\|F_\iota\|_{C^0}}.\]
Let $t=1$ in (\ref{equ-intergral}).  Then 
$\Phi_\iota^1(0,\tilde\theta)=e^{A_\iota}(I+G_\iota(\tilde\theta))$, with 
\[G_\iota(\tilde\theta)=\int_{0}^1e^{-A_\iota s}F_{\iota}(s, \tilde\theta+s\alpha)e^{A_\iota s}\tilde G_\iota^s(0,\tilde\theta)ds.\]
We get 
\[\|G_\iota\|_{C^0}\leq  \int_{0}^1(1+s)^2\|F_\iota\|_{C^0}g_\iota(s)ds<8\|F_\iota\|_{C^0}.\]
Since $\tilde W_\iota$ conjugates (\ref{sys-c-final-proof-original}) to (\ref{sys-c-final-proof-1}), then
\[\Phi^t(0,\tilde\theta)=\tilde W_\iota(t,\tilde\theta+t\alpha)\Phi_\iota^t(0,\tilde\theta)\tilde W_\iota^{-1}(0,\tilde\theta).\]
Let $\bar W_\iota(\tilde\theta)=\tilde W_\iota(0,\tilde\theta)$ and $t=1$. Then $\bar W_\iota\in C^\omega(\T^1, PSL(2,\R))$ and we have 
\[\bar W_\iota(\cdot+\alpha)^{-1}R_{-\varrho_0}e^{\tilde G(\cdot)}\bar W_\iota(\cdot)=e^{A_\iota}(I+G_\iota(\cdot)).\]
Let $W_\iota=\tilde B \bar W_\iota$. Then we get the desired result.

\subsection{Modified KAM scheme}\label{kam}

The proof of the   Proposition \ref{propo-c-ellip-cd} - Proposition \ref{propo-c-parab-liouv} is based on  classical KAM scheme \cite{E92} and   non-standard KAM scheme developed in \cite{HY12}. 
There are four steps altogether in the  proofs of the KAM iteration. In this subsection, we will first revisit the KAM scheme, and provide necessary ingredients for the proof, and leave the full proof of Proposition \ref{propo-c-ellip-cd} - Proposition \ref{propo-c-parab-liouv} to the next subsection. 

We will start from the system 
\begin{equation}\label{initial}
\left\{\begin{array}{ll} \dot x=(A+F(\theta))x\\
\dot \theta=\omega=(1,\alpha)
\end{array}\right.
\end{equation}
where $|F|_{h'}\leq \epsilon$ for some $0<h'<\frac{1}{160}$ with $0<\epsilon<\min\{10^{-8}, 2e^{-c_0q_n^{\frac{1}{\cal A}}h'^2}\}$, and $q_n, h'$ satisfy (\ref{equ-parameter}).

\medskip
\noindent\textbf{Step 1: Elimination of the non-resonant terms}.
\smallskip

In this step, we will try to eliminate the non-resonant terms of the perturbation $F(\theta)$ in (\ref{initial}).  First we give the basic settings. 
For any $\omega\in\R^2, A\in sl(2,\R), \eta>0$, we decompose $\mathcal B_{h'}=\mathcal B_{h'}^{(nre)}(\eta)\oplus\mathcal B_{h'}^{(re)}(\eta)$ (depending on $A, \omega, \eta$) in such a way that for any $Y\in\mathcal B_{h'}^{(nre)}(\eta)$, 
\begin{equation*}\label{equ-condition-non-resonant}
\partial_\omega Y, \ [A, Y] \in \mathcal B_{h'}^{(nre)}(\eta),\ \ \ |\partial_\omega Y-[A, Y]|_{h'}\geq \eta |Y|_{h'}.
\end{equation*}

\begin{Lemma}[\cite{HY12}, Lemma 3.1]\label{lem-hou-you}
Let $\epsilon\in (0, (\frac{1}{10})^{8})$ and $\epsilon^{1/4}\leq \eta<1 $. Then for any $F\in\mathcal B_{h'}$ satisfying $|F|_{h'}\leq \epsilon$, there exist $Y\in\mathcal B_{h'}$ and $F^{(re)}\in\mathcal B_{h'}^{(re)}(\eta)$ such that the system 
\[\left\{\begin{array}{ll} \dot x=(A+F(\theta))x\\
\dot \theta=\omega=(1,\alpha)
\end{array}\right.\]
can be conjugate to the system 
\begin{equation}\label{sys-c-after-eliminate-non-resonant}
\left\{\begin{array}{ll}
 \dot x=(A+F^{(re)}(\theta))x\\
\dot\theta=\omega=(1,\alpha)
\end{array}
\right.
\end{equation}
by the conjugation map $e^Y$, with $|Y|_{h'}\leq \epsilon^{1/2}, |F^{(re)}|_{h'}\leq 2\epsilon$.
\end{Lemma}

For $\eta>0, \tilde\varrho\in\R$, let $\Lambda_1(\eta), \Lambda_2(\tilde\varrho,\eta)$ be the subsets of $\Z^2$ with $\Lambda_1(\eta)=-\Lambda_1(\eta)$, $\Lambda_2(\tilde\varrho,\eta)=-\Lambda_2(\tilde\varrho,\eta)$ , such that 
\begin{equation*}\label{equ-lambda}
k\in\Lambda_1(\eta)\Rightarrow |\la k,\omega\ra|\geq \eta,\ \ \ k\in\Lambda_2(\tilde\varrho,\eta)\Rightarrow |2\tilde\varrho\pm \la k,\omega\ra|\geq \eta.
\end{equation*}
Furthermore, we let 
\[\Lambda_{2,1}(\tilde\varrho,\eta):=\{k\in\Z^2 \ :\ |2\tilde\varrho-\la k,\omega\ra|\geq \eta\}, \ \ \Lambda_{2,2}(\tilde\varrho,\eta):=\{k\in\Z^2\ : \ |2\tilde\varrho+\la k, \omega\ra|\geq \eta\}.\]
Then, $\Lambda_{2,2}(\tilde\varrho,\eta)=-\Lambda_{2,1}(\tilde\varrho,\eta)$. These sets are very useful when one tries to analyze the structure of $F^{(re)}$. In the following, we will analyze this.

Obviously, the structure of  $F^{(re)}$ depends on the constant matrix $A$. If $A=2\pi\tilde\varrho J$, then in order to explore the structure clearly,  it is better to state the result in $su(1,1)$. 
Recall that  $su(1,1)$ is  the space consisting of matrices of the
form $\left(
\begin{array}{l}
\mathrm{i} t\ \ \ \ \ v
\\
\overline{v}\ \ -\mathrm{i}t
\end{array}
\right)$ $(t\in\mathbb{R},\ v\in\mathbb{C}),$ 
and  $sl(2,\mathbb{R})$ is isomorphic to $su(1,1)$ by the rule $A\mapsto MAM^{-1},$ where
$M=\frac{1}{1+\mathrm{i}}\left(
\begin{array}{l}
1\ \ -\mathrm{i}
\\
1 \quad\ \ \mathrm{i}
\end{array}
\right),$
since a simple calculation yields
\begin{equation*}
\begin{split}
M\left(
\begin{array}{l}
x\ \qquad y+z
\\
y-z \quad\ -x
\end{array}
\right)M^{-1}=\left(
\begin{array}{l}
\mathrm{i}z\ \qquad x-\mathrm{i}y
\\
x+\mathrm{i} y \quad\ -\mathrm{i} z
\end{array}
\right), x,y,z\in\mathbb{R}.
\end{split}
\end{equation*}
Within this context, we have the following:

\begin{Corollary}[\cite{HY12}, Corollary 3.2]\label{cor-hou-you}
Let 
$A=2\pi\tilde\varrho J$. Then the conclusion of Lemma \ref{lem-hou-you} holds with 
$F^{(re)}$ in the form
\begin{eqnarray*}\label{equ-form-elliptic-non-resonant}
\nonumber M F^{(re)} M^{-1}&=&\sum_{k\in\Lambda_1^c} \left(\begin{array}{cc} i \hat F_{-}^{(re)}(k) & 0  \\ 0 & - i \hat F_{-}^{(re)}(k) \end{array}\right)e^{2\pi i\la k,\theta\ra} \\& +& \sum_{k\in\Lambda_{2,1}^c}  \hat F_{+,1}^{(re)}(k) \left(\begin{array}{cc} 0 &  1  \\ 0  & 0 \end{array}\right)  e^{2\pi i\la k,\theta\ra}+  \sum_{k\in\Lambda_{2,2}^c}  \hat F_{+,2}^{(re)}(k)\left(\begin{array}{cc} 0 &  0  \\ 1  & 0 \end{array}\right) e^{2\pi i\la k,\theta\ra},
\end{eqnarray*}
where $\hat F^{(re)}_{\pm}(k)=\frac{\hat F^{(re)}_{12}(k)\pm \hat F^{(re)}_{21}(k)}{2}$, $\hat F_{+,1}^{(re)}(k)=\hat F_{11}^{(re)}(k)-i\hat F^{(re)}_+(k)$, $\hat F_{+,2}^{(re)}(k)=\hat F_{11}^{(re)}(k)+i\hat F^{(re)}_+(k)$, and $\Lambda_1:=\Lambda_1(\eta), \Lambda_{2,j}:=\Lambda_{2,j}(\tilde\varrho, \eta)$ with $j=1, 2$.
\end{Corollary}

\smallskip


If $A$ is parabolic, then we have the following:

\begin{Corollary}\label{cor-parabolic-non-resonant}
Let $A=\left(\begin{array}{cc}0 & c^* \\ 0&0\end{array}\right)$ with $|c^*|\leq 1$. Then the conclusion of Lemma \ref{lem-hou-you} holds with $F^{(re)}$ in the form
\begin{equation*}\label{equ-form-parabolic-non-resonant}
F^{(re)}=\sum_{k\in\Lambda_1^c}\hat F^{(re)}(k)e^{2\pi i\la k,\theta\ra},
\end{equation*}
where $\Lambda_1:=\Lambda_1(\eta^{1/3})$.
\end{Corollary}
\begin{pf}

Denote $\tilde{\mathcal B}_{h'}^{(nre)}=\{Y\in \mathcal B_{h'}\ | \ Y=\sum_{k\in\Lambda_1}\hat Y(k)e^{2\pi i\la k,\theta\ra} \}$. Then for any $Y\in \tilde{\mathcal B}_{h'}^{(nre)}$, we have $\partial_\omega Y=\sum_{k\in\Lambda_1}2\pi i\la k,\omega\ra\hat Y(k)e^{2\pi i \la k,\theta\ra}\in \tilde{ \mathcal B}_{h'}^{(nre)}$,
\[[A, Y]=\left(\begin{array}{cc} c^*Y_{21} & -2c^*Y_{11}\\
0 & -c^*Y_{21}
\end{array}\right)\in\tilde{ \mathcal B}_{h'}^{(nre)},\]
and then 
\[\partial_\omega Y-[A, Y]=\left(\begin{array}{cc}\partial_\omega Y_{11} - c^*Y_{21}& \partial_\omega Y_{12}+2c^* Y_{11}\\
\partial_\omega Y_{21} & -\partial_\omega Y_{11}+c^*Y_{21} \end{array}\right).\]
Then one can check that $|\partial_\omega Y-[A, Y]|_{h'}\geq \eta|Y|_{h'}$, which means $\tilde{\cal B}_{h'}^{(nre)}\subseteq \cal B_{h'}^{(nre)}(\eta)$, and the result follows.
\end{pf}

Therefore, by Corollary \ref{cor-hou-you} and  Corollary \ref{cor-parabolic-non-resonant},  in order to analyze the structure of $F^{(re)}$, one only need to analyze the structure of 
$\Lambda_1^c(\eta)$ and $\Lambda_{2,j}^c(\tilde\varrho, \eta)$ $(j=1,2)$.

\begin{Lemma}\label{lem-lambda}
Let $\tilde\varrho\in\R, 0<\eta\leq\frac{1}{4q_n^{ \mathcal A^4}}$. If $q_{n+l}\leq q_n^{\mathcal A^4}$ for some $l\in\N$, then
\begin{equation}\label{equ-lambda-1}
\Lambda_1^c(\eta)\cap \{k\in\Z^2 : |k|< q_{n+l} \}=\{0\},
\end{equation}
\begin{equation}\label{equ-lambda-2}
\#\left(\Lambda_{2,j}^c(\tilde\varrho,\eta)\cap \left\{k\in\Z^2 : |k|<\frac{q_{n+l}}{2}\right\}\right)\leq 1,  \ \  (j=1,2).
\end{equation}
\end{Lemma}
\begin{pf}
For any $k\in\Z^2$ with $0<|k|< q_{n+l} $, we have 
\[|\la k,\omega\ra|\geq \frac{1}{2q_{n+l}}\geq \frac{1}{2q_n^{\mathcal A^4}}> \eta,\]
which implies (\ref{equ-lambda-1}).

Moreover,  if there exist distinct $k, k'\in \Lambda_{2,1}^c(\tilde\varrho, \eta)$ with $|k|, |k'|<\frac{q_{n+l}}{2}$, then we have 
$|\la k-k',\omega\ra|<2\eta$, which contradicts with the fact 
\[|\la k-k',\omega\ra|\geq \frac{1}{2q_{n+l}}\geq \frac{1}{2q_n^{\mathcal{A}^4}}\geq 2\eta.\]
This implies (\ref{equ-lambda-2}) with $j=1$. Since $\Lambda_{2,2}(\tilde\varrho,\eta)=-\Lambda_{2,1}(\tilde\varrho,\eta)$, we also obtain (\ref{equ-lambda-2}) with $j=2$.
\end{pf}

In order to analyze the structure of $\Lambda_1^c(\eta)$ and $\Lambda_{2,j}^c(\tilde\varrho, \eta)$ $(j=1,2)$ for the case $q_{n+1}>q_n^{\mathcal A}$, we will need the following lemma.

\begin{Lemma}[\cite{HY12}, Lemma 4.1]\label{lem-hou-you-2}
For any $k=(k_1, k_2)\in \Z^2$ satisfying 
\begin{itemize}
\item $|k|:=|k_1|+|k_2|\leq \frac{q_{n+1}}{6}$, and
\item $k\neq l(p_n, -q_n)$, $l\in\Z$,
\end{itemize}
we have $|\la k,\omega\ra|\geq \frac{1}{7q_n}$.
\end{Lemma}

\begin{Remark}
In fact, in the case $q_{n+1}>q_n^{\mathcal A}$ with $n$ sufficiently large, we have $|\la k, \omega\ra|\geq\frac{1}{7q_{n}}$ for any $|k|\leq \frac{q_{n+1}}{3}$ and $k\neq l (p_n, -q_n), l\in\Z$.
\end{Remark}

\begin{Remark}\label{remark-liouvillean-lambda-structure}
Let $\eta=\epsilon^{1/4}$. If $\eta^{1/3}\leq \frac{1}{7q_n}$, by the above lemma, it is obvious that
\begin{equation*}
\Lambda_1^c(\eta^{1/3})\cap \left\{k\in\Z^2: |k|< \frac{q_{n+1}}{6}\right\}\subseteq \left\{k=l(p_n, -q_n): l\in\Z, |k|< \frac{q_{n+1}}{6}\right\}.
\end{equation*}
Moreover, if $|\tilde\varrho|<\frac{1}{2}\epsilon^{1/4}$, then we have 
\begin{equation*}
\Lambda_2^c(\tilde\varrho,\eta)\cap \{k\in\Z^2: |k|<\frac{q_{n+1}}{6}\}\subseteq \left\{k=l(p_n, -q_n): l\in\Z, |k|< \frac{q_{n+1}}{6}\right\}.
\end{equation*}
\end{Remark}

\medskip

\noindent
\textbf{Step 2: Rotation.}
\smallskip

For the Diophantine case ($q_{n+1}\leq q_n^{\mathcal A}$), we follow the ideas in \cite{E92}, taking a rotation to eliminate the resonant term. 
Otherwise, we follow the ideas in \cite{HY12}:  
In general, if we don't have the assumption $|\tilde\varrho|<\frac{1}{2}\epsilon^{1/4}$, 
then we have 
\begin{equation*}
\Lambda_{2,1}^c(\tilde\varrho,\epsilon^{1/4})\cap \{k\in\Z^2: |k|<\frac{q_{n+1}}{6}\}\subseteq \left\{k=k_*+l(p_n, -q_n): l\in\Z, |k|< \frac{q_{n+1}}{6}\right\},
\end{equation*}
for some $k_*\in\Lambda_{2,1}^c(\tilde\varrho,\epsilon^{1/4})\cap \{k\in\Z^2: |k|<\frac{q_{n+1}}{6}\}$.
In this case, we can take a rotation to make the truncated system a periodic system:

\begin{Lemma}[\cite{HY12}, Lemma 5.2]\label{lem-hou-you-3}
Let $A=2\pi \tilde\varrho J$ and $\Lambda_2:=\Lambda_2(\tilde\varrho, \epsilon^{1/4})$, $\Lambda_{2,j}:=\Lambda_{2,j}(\tilde\varrho, \epsilon^{1/4})$ (j=1,2). If $\Lambda_2^c\cap \{k\in\Z^2 : |k|<\frac{q_{n+1}}{6}\}\neq \emptyset$, then there exists $\tilde{\cal Q}\in C_{\frac{h'}{3}}^\omega(\T^2, PSL(2,\R))$ that conjugates (\ref{sys-c-after-eliminate-non-resonant}) to
\begin{equation}\label{sys-after-rotation}
\left\{\begin{array}{ll}
\dot x =( A' + F' (\theta))x\\
\dot\theta=\omega
\end{array}\right.
\end{equation}
where $$ A' =2\pi (\tilde\varrho-\frac{\la k_*,\omega\ra}{2} )J, \qquad \mathcal{ T}_{\frac{q_{n+1}}{6}} F' =\sum_{k=l(p_n, -q_n) \atop{|k|< \frac{q_{n+1}}{6}}}\hat { F'} (k)e^{2\pi i\la k,\theta\ra},$$
with estimate $| F' |_{h'/3}\leq 2\epsilon^{3/4}$. Moreover,  $\tilde{\cal Q}(\theta)=\cal Q(\theta)e^{Y(\theta)}$ takes the form 
 $$\cal Q(\theta)=\left(\begin{array}{cc} \cos(\pi\la k_*,\theta\ra) & \sin(\pi\la k_*,\theta\ra)\\ -\sin(\pi\la k_*,\theta\ra) & \cos(\pi \la k_*,\theta\ra)
\end{array}\right),$$ 
where $k_*\in \Lambda_{2,1}^c\cap \{k\in\Z^2 : |k|<\frac{q_{n+1}}{6}\}$. Furthermore, we have the following estimates
 $$|Y|_{\frac{h'}{3}}\leq \epsilon^{\frac{3}{8}}, \qquad \|\cal Q\|_{C^1}\leq 2\pi |k_*|<\frac{\pi}{3}q_{n+1}.$$
\end{Lemma}

\medskip

\noindent\textbf{Step 3: Floquet's theorem.}

By Lemma \ref{lem-hou-you-3}, one can observe that due to the special form of the perturbation,  the truncated system 
\begin{equation}\label{trunkam}
\left\{\begin{array}{ll}
\dot x =( A' +\mathcal{ T}_{\frac{q_{n+1}}{6}} F'(\theta))x = ( A' + \sum_{k=l(p_n, -q_n) \atop{|k|< \frac{q_{n+1}}{6}}}\hat {F'} (k)e^{2\pi i\la k,\theta\ra})x\\
\dot\theta=\omega=(1,\alpha)
\end{array}\right.
\end{equation}
is in fact a periodic system, which can be conjugate to the constant by the famous Floquet's theorem. The following lemma will provide estimates that are sufficient for us.

\begin{Lemma}[\cite{HY12}, Lemma 7.1]\label{lem-hou-you-floquet}
The system 
\[\left\{
\begin{array}{ll}
\dot x=G(\theta)x\\
\dot \theta=\omega=(1,\alpha)
\end{array}
\right.\]
where $G\in\mathcal B_{h''}$ for some $h''>0$, $|G|_{h''}<\varepsilon$ and is of the form
\[G(\theta_1,\theta_2)=\sum_{l\in\Z}\hat G(lp, -lq)e^{2\pi i l(p\theta_1-q\theta_2)},\] 
with $(p, -q)\in\Z^2$ fixed, can be conjugate to some constant system  
\[\left\{\begin{array}{ll}
\dot x=D\\
\dot \theta=\omega
\end{array}
\right.\]
by a conjugation map $B\in C^\omega(\T^2, PSL(2,\R))$, with $D\in sl(2,\R)$, and the estimations 
\[\|B\|_{h''}\leq e^{2\frac{|G|_{h''}}{|\tau|}( 2+(|p|+|q|)h'' )},\ \
\textrm{and}\ \ \|D\|\leq |\tau|e^{\frac{2|G|_{h''}}{|\tau|}},\]
where $\tau=p-q\alpha$.
\end{Lemma}

\medskip


\noindent\textbf{Step 4: Normalization of the constant matrix.}

By Floquet's theorem (Lermma \ref{lem-hou-you-floquet}), the truncated system \eqref{trunkam} can be conjugate to system with constant matrix $\tilde A$, and consequently, the non-truncated system \eqref{sys-after-rotation} is conjugate to 
\begin{equation}\label{sys-c-before-normalization}
\left\{\begin{array}{ll}
\dot x=(\tilde A+\tilde F (\theta))x\\
\dot\theta=\omega
\end{array}.\right.
\end{equation}
However, the constant matrix $\tilde A$ may be out of control. Before finishing one step of KAM iteration, we will try to normalize the constant matrix. If the constant matrix $\tilde A$ is elliptic or hyperbolic, then we need the following:

\begin{Lemma}[\cite{HY12}, Lemma 8.1]\label{lem-hou-you-normalization}
Let $A\in sl(2,\R)$ satisfy $spec(A)=\{i\varrho, -i\varrho\}$ with $0\neq \varrho\in\R$. There exists $P\in SL(2,\R)$ such that $\|P\|\leq 2(\frac{\|A\|}{|\varrho|})^{1/2}$ and $P^{-1}AP=\left(\begin{array}{cc} 0 &\varrho \\ -\varrho & 0\end{array}\right)$.
\end{Lemma}

\begin{Lemma}[\cite{Puig06}, Proposition 18]\label{lem-puig-normalization}
Let $A\in sl(2,\R)$ satisfy $spec(A)=\{\lambda, -\lambda\}$ with $0\neq \lambda\in\R$. There exists $P\in SL(2,\R)$ such that $\|P\|\leq (\frac{\|A\|}{|\lambda|})^{1/2}$ and $P^{-1}AP=\left(\begin{array}{cc} \lambda &0 \\ 0 & -\lambda \end{array}\right)$.
\end{Lemma}

If the initial system    \eqref{initial} is not uniformly hyperbolic, then the conjugated system   \eqref{sys-c-before-normalization}  is also not uniformly hyperbolic, since  uniformly hyperbolic is a conjugate invariant. This will provide us additional information, when we try to normalize the constant matrix.

\begin{Lemma}\label{lem-c-normalization}
Let 
$\mathcal K>0$, $0<h'<\frac{1}{160}$, and $  \cal Kh'^4>4$. 
For the system (\ref{sys-c-before-normalization}), with estimates $\|\tilde A\|\leq e^{\frac{\mathcal Kh'^4}{4}}$ and $\|\tilde F\|_{h''}\leq   e^{-\mathcal Kh'}$ for some $h''>0$. If   (\ref{sys-c-before-normalization}) is not uniformly hyperbolic,  then there exists $P\in SL(2,\R)$ that conjugates (\ref{sys-c-before-normalization}) to 
\begin{equation}\label{sys-c-after-normalization}
\left\{\begin{array}{ll} 
\dot x=(\bar A +\bar F(\theta))x\\
\dot\theta=\omega
\end{array}\right.
\end{equation}
where 
\begin{enumerate}
\item either $\bar A=2\pi\bar\varrho J$ for some $\bar\varrho\in\R$ with
\begin{enumerate}[(i)]
\item either $\left\{\begin{array}{ll}
\|P\|\leq 2e^{\mathcal Kh'^2}\\
\|\bar F\|_{ h''}\leq 4e^{-\frac{\mathcal Kh'}{2}}
\end{array},\right.$
\item or $\left\{\begin{array}{ll}
\|P\|\leq 2\\
\|\bar F\|_{h''}\leq 2e^{-\frac{c_0\mathcal Kh'^2}{3}}
\end{array};\right.$
\end{enumerate}
\item or $\bar A=\left(\begin{array}{cc} 0 & \bar c \\ 0 & 0\end{array}\right)$ where $\bar c\in\R$ with 
\[\left\{\begin{array}{lll}
\|P\|\leq 2e^{\frac{\mathcal Kh'^4}{2}}\\
\|\bar F\|_{h''}\leq 9e^{-\frac{3}{4}\mathcal Kh'^2}\\
e^{-\frac{c_0}{3}\mathcal Kh'^2}\leq |\bar c|\leq 4e^{-\frac{3}{4}\mathcal Kh'^4}
\end{array}.\right.\]
\end{enumerate}
\end{Lemma}

\begin{pf}
We denote  $\tilde A =\left(\begin{array}{cc} a_{11} & a_{12}\\ a_{21} & -a_{11}\end{array}\right)$, and divide the proof into several different cases.

\medskip
\textbf{Case 1:} $\det \tilde A =(2\pi \tilde\varrho)^2>0$. That is, $\tilde A$ is elliptic.
\smallskip

{\it Case 1.1: $|2\pi\tilde \varrho|\geq e^{-\mathcal Kh'^2}$}. In this case, by Lemma \ref{lem-hou-you-normalization}, there exists $P_1\in SL(2,\R)$ with $\|P_1\|\leq 2(\frac{\|\tilde A\|}{2\pi\tilde\varrho})^{1/2}<2e^{\mathcal Kh'^2}$, such that $P_1^{-1}\tilde A P_1=2\pi\tilde\varrho J$. We let $\bar A:=2\pi\tilde\varrho J$, $P:=P_1$ and $\bar F:=P_1^{-1}\tilde FP_1$. Then
\[\|\bar F\|_{h''}\leq \|P_1\|^2\|\tilde F \|_{h''}<4 e^{-\frac{\mathcal Kh'}{2}}.\]

{\it Case 1.2:  $|2\pi\tilde\rho|<e^{-\mathcal Kh'^2}$.}   In this case, $\det \tilde A=-a_{11}^2-a_{12}a_{21}=(2\pi\tilde\varrho)^2>0$, which implies that 
\[|a_{12}a_{21}|=a_{11}^2+(2\pi\tilde\varrho)^2\geq (2\pi\tilde\varrho)^2>0.\]
Hence, we have $\max\{|a_{12}|, |a_{21}|\}\geq |2\pi\tilde\varrho|>0$. Without loss of generality, we assume that $|a_{12}|\geq |2\pi\tilde\varrho|$. Then the system (\ref{sys-c-before-normalization}) can be rewritten as
\begin{equation}
\nonumber \left\{\begin{array}{ll}
\dot x= ( \tilde A ^{(1)}+ \tilde F ^{(1)}(\theta))x\\
\dot \theta=\omega
\end{array}
\right.
\end{equation}
where
\[ \tilde A^{(1)}=\left(\begin{array}{cc} a_{11} & a_{12}\\
a_{21}+\frac{(2\pi\tilde\varrho)^2}{a_{12}} & -a_{11}
\end{array}\right), \ \ \ \    \tilde F^{(1)}(\theta)=\tilde F(\theta)+\left( \begin{array}{cc}
0 & 0\\  
-\frac{(2\pi\tilde\varrho)^2}{a_{12}} & 0
\end{array} \right),\]
with $\det  \tilde A^{(1)}=0$, and 
\[\|  \tilde A^{(1)}\|\leq \|\tilde A\|+2\pi|\tilde\rho|< 2e^{\frac{\mathcal Kh'^4}{4}},\] 
\[\|\tilde F^{(1)}\|_{h''}\leq \|\tilde F\|_{h''}+2\pi|\tilde\varrho|<2e^{-\mathcal Kh'^2}.\]
Then we will reduce it to the following case:

\medskip

\textbf{Case 2:} $\det \tilde A=-\lambda^2\leq 0$.  That is, $\tilde A$ is hyperbolic or parabolic. First we claim that we only need to consider the case $|\lambda|< (\|\tilde F \|_{h''}\|\tilde A\|)^{1/3}$, since we assume the system   \eqref{sys-c-before-normalization}  is  not uniformly hyperbolic. Otherwise, by  Lemma \ref{lem-puig-normalization},  there exists $P_1\in SL(2,\R)$ such that $\|P_1\|\leq  (\frac{\|\tilde A\|}{|\lambda|})^{1/2}$ and $P_1^{-1}\tilde AP_1=\left(\begin{array}{cc} \lambda &0 \\ 0 & -\lambda \end{array}\right)$. Then   \eqref{sys-c-before-normalization}  is conjugate to 
\begin{equation}\label{sys-c-in-normalization-1}
 \left\{\begin{array}{ll}
\dot x= ( \left(\begin{array}{cc} \lambda &0 \\ 0 & -\lambda \end{array}\right)+P_1^{-1} \tilde F(\theta)P_1)x\\
\dot \theta=\omega
\end{array}
\right.
\end{equation}
with $\|P_1^{-1} \tilde F(\theta)P_1\|_{h''} \leq  \frac{\|\tilde A\| \|\tilde F \|_{h''}}{|\lambda|}$. Therefore, if $|\lambda|\geq   ( \frac{\|\tilde A\| \|\tilde F \|_{h''}}{|\lambda|}  )^{1/2}$, i.e., $|\lambda|\geq(\|\tilde F \|_{h''}\|\tilde A\|)^{1/3},$ the system  \eqref{sys-c-in-normalization-1} is uniformly hyperbolic by the usual cone-field criterion  \cite{Yoc}, which contradicts with our assumption.   

Then either in this case, or in Case 1.2, we can write the system  (\ref{sys-c-before-normalization}) as 
\begin{equation}\label{sys-c-in-normalization}
\left\{\begin{array}{ll}
\dot x= (\tilde A^{(2)}+  \tilde F^{(2)}(\theta))x\\
\dot \theta=\omega
\end{array}
\right.
\end{equation}
with  $\det \tilde A^{(2)}=-\lambda^2\leq 0$, $|\lambda|<(\|\tilde A \|\|\tilde F\|_{h''})^{\frac{1}{3}}<e^{-\frac{\mathcal Kh'}{4}}$, $\| \tilde A^{(2)}\|<2e^{\frac{\mathcal Kh'^4}{4}}$, $\| \tilde F^{(2)}\|_{h''}<2e^{-\mathcal Kh'^2}$. We assume that  $ \tilde A^{(2)}=\left(\begin{array}{cc}\tilde a_{11} & \tilde a_{12}\\ \tilde a_{21}  & -\tilde a_{11}\end{array}\right)$. Then there exists $P_2\in SL(2,\R)$ with $\|P_2\|\leq 2$ such that $P_2^{-1} \tilde A^{(2)}P_2=\left(\begin{array}{cc}\lambda & \bar a_{12} \\
 0 & -\lambda \end{array}\right)$ with $\bar a_{12}=\tilde a_{12}-\tilde a_{21}$. Under this conjugation map, (\ref{sys-c-in-normalization}) is conjugate to
 \begin{equation}\label{sys-c-in-normalization-2}
 \left\{\begin{array}{ll}
 \dot x=(\tilde A^{(3)}+ \tilde F^{(3)}(\theta))x\\
 \dot\theta=\omega
 \end{array}
 \right.
 \end{equation}
 where $\tilde A^{(3)}=\left(\begin{array}{cc}0 & \bar a_{12}\\ 0 & 0 \end{array}\right)$, $ \tilde F^{(3)}=P_2^{-1} \tilde F^{(2)}P_2+\left(\begin{array}{cc}\lambda & 0 \\ 0 & -\lambda \end{array}\right)$, with \[|\bar a_{12}|\leq 4e^{\frac{\mathcal Kh'^4}{4}}, \  \  \| \tilde F^{(3)}\|_{h''}\leq 9e^{-\mathcal Kh'^2}.\]

Recall that $c_0=\frac{1}{2\cdot 48^2}$.  Now we consider the following  three sub-cases according to  the value of $|\bar a_{12}|$:\\
  
 {\it Case 2.1: $|\bar a_{12}|\leq e^{-\frac{c_0}{3}\mathcal Kh'^2}$.}  We let $\bar A:=0$, $\bar F:=\tilde F^{(3)}+\tilde A^{(3)}$, $P:=P_2$. Then $\|P\|\leq 2$ and 
 \[\|\bar F\|_{h''}\leq \|\tilde F^{(3)}\|_{h''}+|\bar a_{12}|<2e^{-\frac{c_0}{3}\mathcal Kh'^2}.\]  
 
 \smallskip

 {\it Case 2.2: $e^{-\frac{c_0}{3}\mathcal Kh'^2}<|\bar a_{12}|\leq e^{-\frac{3}{4}\mathcal Kh'^4}$.} Let $\bar A:=\tilde A^{(3)}$ with $\bar c=\bar a_{12}$, $\bar F:= \tilde F^{(3)}$, $P:=P_2$. Then $\|P\|\leq 2$, $e^{-\frac{c_0}{3}\mathcal Kh'^2}<|\bar c|\leq e^{-\frac{3}{4}\mathcal Kh'^4}$, and
 \[\|\bar F\|_{h''}= \| \tilde F^{(3)}\|_{h''}\leq 9e^{-\mathcal Kh'^2}.\]
 
 \smallskip

 {\it Case 2.3: $|\bar a_{12}|>e^{-\frac{3}{4}\mathcal Kh'^4}$.} Let $H=\left(\begin{array}{cc} e^{\frac{\mathcal Kh'^4}{2}} & 0\\
0 & e^{-\frac{\mathcal Kh'^4}{2}}
\end{array}\right)$. Then  $H$ conjugates the  system (\ref{sys-c-in-normalization-2}) to (\ref{sys-c-after-normalization}) where $\bar A=\left(\begin{array}{cc}
0 & \bar c\\
0 & 0
\end{array}\right)$ with $ \bar c=e^{-\mathcal Kh'^4}\bar a_{12}$,  and $\bar F= H^{-1} \tilde F^{(3)} H$. Let $P:=P_2H$. Then $\|P\|\leq 2e^{\frac{\mathcal Kh'^4}{2}}$, and
\[\|\bar F\|_{h''}\leq \|H\|^2\|\tilde F^{(3)}\|_{h''}\leq9e^{\mathcal Kh'^4}e^{-\mathcal Kh'^2}<9e^{- \frac{3}{4}\mathcal Kh'^2 }.\]  Moreover, 
we have $|\bar c|>e^{-\frac{7}{4}\mathcal Kh'^4}> e^{-\frac{c_0}{3}\mathcal Kh'^2}$, and \[|\bar c|\leq 4e^{-\mathcal Kh'^4} e^{\frac{\mathcal Kh'^4}{4}}=4e^{-\frac{3}{4}\mathcal Kh'^4}.\]

\end{pf}

\subsection{Proof of Proposition \ref{propo-c-ellip-cd}-\ref{propo-c-parab-liouv}}\label{sec-propositions}
With the above preparing lemmas, now we can finish the whole proof.

\subsubsection{Proof of Proposition \ref{propo-c-ellip-cd}}
By (\ref{equ-parameter}), we can check that $2e^{-c_0q_n^{\frac{1}{ \cal A}}h'^2}<10^{-8}$. Then by  Lemma \ref{lem-hou-you}, there exist $Y_1\in \mathcal B_{h'}$ and $F^{(re)}\in\mathcal B_{h'}^{(re)}(\epsilon^{\frac{1}{4}})$ such that (\ref{sys-continuous-original}) can be conjugate to (\ref{sys-c-after-eliminate-non-resonant}) by the conjugation map $e^{Y_1}$, and $|Y_1|_{h'}\leq \epsilon^{1/2}$, $|F^{(re)}|_{h'}\leq 2\epsilon$. 

\textit{(a)} If $\Lambda_2^c\cap\{k\in\Z^2 : |k|<\frac{q_{n+l}}{2}\}=\emptyset$,  
then by Corollary \ref{cor-hou-you} and Lemma \ref{lem-lambda}, $F^{(re)}$ has the form
\begin{equation*}
\mathcal T_{\frac{q_{n+l}}{2}}F^{(re)}=\left(\begin{array}{cc}
0 & \hat F_-^{(re)}(0)  \\  -\hat F^{(re)}_-(0) & 0
\end{array}\right),
\end{equation*}
since $\epsilon^{\frac{1}{4}}<\frac{1}{4q_n^{\cal A^4}}$ by (\ref{equ-parameter}).
Let $B=e^{Y_1}$, $\tilde A_1=A+\hat F_-^{(re)}(0)J=:2\pi\bar \varrho  J$, $\tilde F_1=\mathcal R_{\frac{q_{n+l}}{2}}F^{(re)}$. Then we have 
\[\|B\|_{C^1}=\|e^{ Y_1 }\|_{C^1}< 2 ,\]
and
\begin{eqnarray*}
\|\tilde F_1\|_{h'/2}&\leq& \sum_{|k|\geq q_{n+l}/2}|\hat F^{(re)}(k)| e^{ \pi |k|h'}\\
&=&\sum_{|k|\geq q_{n+l}/2}|\hat F^{(re)}(k)| e^{2 \pi |k|h'}e^{-\pi |k|h'}\\
&<&\epsilon  e^{-q_{n+l}h'}.\end{eqnarray*}

\textit{(b)} If $\Lambda_{2}^c\cap\{k\in\Z^2 : |k|<\frac{q_{n+l}}{2}\}\neq\emptyset$, then by Lemma \ref{lem-lambda}, there exists a unique $k_*\in \Lambda_{2,1}^c\cap\{k\in\Z^2 : |k|<\frac{q_{n+l}}{2}\}$. Furthermore, by Corollary \ref{cor-hou-you},  $F^{(re)}$ has the form
\begin{equation*}
M \mathcal T_{\frac{q_{n+l}}{2}}F^{(re)} M^{-1}=\left(\begin{array}{cc}
i \hat F_-^{(re)}(0) &  0 \\ 0 &  -i \hat F^{(re)}_-(0)
\end{array}\right)+   \left(\begin{array}{cc} 0 &   \hat F_{+,1}^{(re)}(k_*)e^{2\pi i\la k_*,\theta\ra}  \\ \overline{\hat F_{+,1}^{(re)}(k_*)}e^{-2\pi i\la k_*,\theta\ra} & 0 \end{array}\right). 
\end{equation*}
%
%
Let $\cal Q(\theta)=\left(\begin{array}{cc} \cos(\pi\la k_*,\theta\ra) & \sin(\pi\la k_*,\theta\ra)\\ -\sin(\pi\la k_*,\theta\ra) & \cos(\pi \la k_*,\theta\ra)
\end{array}\right)$. Then direct calculation shows that  $\cal Q$ conjugates (\ref{sys-c-after-eliminate-non-resonant}) to
\begin{equation}\label{sys-c-ellip-cd-after-rotation}
\left\{\begin{array}{ll}\dot x=(\tilde A+\tilde F(\theta))x\\
\dot \theta=\omega,
\end{array}\right.
\end{equation}
where $\tilde F=\cal Q^{-1}\mathcal R_{\frac{q_{n+l}}{2}}F^{(re)}\cal Q$,  and $\tilde A$ has the form
$$ \tilde A=(2\pi\varrho-\pi\la k_*,\omega\ra +\hat F^{(re)}_-(0) )J + M^{-1}  \left(\begin{array}{cc} 0 &   \hat F_{+,1}^{(re)}(k_*) \\ \overline{\hat F_{+,1}^{(re)}(k_*)} & 0 \end{array}\right) M.$$
Then one has 
\begin{eqnarray*}
\| \cal Q\|_{C^1} &\leq&  2\pi |k_*|< \pi q_{n+l}\\
\|\tilde A \|&\leq& \pi |2\varrho- \la k_*,\omega\ra|+ 4\epsilon < \pi \epsilon^{1/4} +4\epsilon < 8 \epsilon^{1/4},\end{eqnarray*} and furthermore, we have 
\begin{eqnarray*}
\|\tilde F\|_{h'/6}&\leq& \|\mathcal Q\|_{h'/6}^2\sum_{|k|\geq \frac{q_{n+l}}{2}}|\hat F^{(re)}(k)|e^{2\pi |k|\frac{h'}{6}}\\
&\leq& 2\epsilon e^{ \frac{\pi |k_*|h'}{3}} \sum_{|k|\geq \frac{q_{n+l}}{2}} e^{-2\pi |k| h'} e^{2\pi |k|\frac{h'}{6}}<  \epsilon e^{-q_{n+l}h'}.
\end{eqnarray*}
Then by Lemma \ref{lem-c-normalization}, there exists $P\in SL(2,\R)$ that conjugates (\ref{sys-c-ellip-cd-after-rotation}) to (\ref{sys-c-ellip-cd-final}), where 
\begin{itemize}
\item either $ \tilde A_1=2\pi\bar \varrho J$ for some $\bar\varrho\in\R$, with
\begin{itemize}
\item either $\left\{\begin{array}{ll}
\|P\|\leq 2e^{q_{n+l}h'^2}\\
\|\tilde F_1\|_{\tilde h}\leq 4e^{-\frac{1}{2}q_{n+l}h'^2}
\end{array}\right.$,
\item or $\left\{\begin{array}{ll}
\|P\|\leq 2\\
\|\tilde F_1\|_{\tilde h}\leq 2e^{-\frac{ c_0}{3} q_{n+l}  h'^2}
\end{array}\right.$;
\end{itemize}
\item or $\tilde A_1=\left(\begin{array}{cc} 0 & \bar c \\ 0 & 0\end{array}\right)$ where $\bar c\in\R$, with 
\[\left\{\begin{array}{lll}
\|P\|\leq 2e^{\frac{1}{2}q_{n+l}  h'^4}\\
\|\tilde F_1\|_{ \tilde h}\leq 9e^{-\frac{3 }{4}q_{n+l}  h'^2}\\
e^{-\frac{c_0}{3} q_{n+l}  h'^2}\leq |\bar c|\leq 4e^{-\frac{3}{4}q_{n+l}  h'^4}
\end{array}.\right.\]
\end{itemize}
Let $B=e^{Y_1}\mathcal Q P$. Then we get the desired result. \qed

\subsubsection{Proof of Proposition \ref{propo-c-ellip-liouv}}
Since $2e^{-c_0q_n^{\frac{1}{ \cal A}}h'^2}<10^{-8}$, then by  Lemma \ref{lem-hou-you},  there exist $Y_1\in \mathcal B_{h'}$ and $F^{(re)}\in\mathcal B_{h'}^{(re)}(\epsilon^{\frac{1}{4}})$ such that (\ref{sys-continuous-original}) can be conjugate to (\ref{sys-c-after-eliminate-non-resonant}) by the conjugation map $e^{Y_1}$, and $|Y_1|_{h'}\leq \epsilon^{1/2}$, $|F^{(re)}|_{h'}\leq 2\epsilon$. 

\textit{(a)} If $\Lambda_2^c\cap\{k\in\Z^2 : |k|<\frac{q_{n+1}}{6}\}=\emptyset$,  
then by Corollary \ref{cor-hou-you}, we have
\[\cal T_{\frac{q_{n+1}}{6}}F^{(re)}=\sum_{k\in\Lambda_1^c, \atop{|k|<\frac{q_{n+1}}{6}}}\left(\begin{array}{cc} 0 &\hat F_-^{(re)}(k) \\  -\hat F_-^{(re)}(k)   & 0 \end{array}\right)e^{2\pi i\la k,\theta\ra}.\]
Let 
\[E(\theta)=\sum_{ k\in\Lambda_1^c, \atop{0< |k|< \frac{q_{n+1}}{6} }} \left(\begin{array}{cc} 0 &\hat F_-^{(re)}(k) \\  -\hat F_-^{(re)}(k)   & 0 \end{array}\right)\frac{e^{2\pi i\la k,\theta\ra}}{2\pi i\la k,\omega\ra}.\]
Then $e^{E(\theta)}$ conjugates the system (\ref{sys-c-after-eliminate-non-resonant}) to 
\[\left\{ \begin{array}{ll} \dot x=  ( A+\hat F_{-}^{(re)}(0)J+ 
e^{-E(\theta)}\mathcal R_{\frac{q_{n+1}}{6}} F^{(re)}(\theta) e^{E(\theta)} \ ) x\\
\dot \theta=\omega
\end{array}.\right.\]
Since for any $0< |k|< \frac{q_{n+1}}{6} $, we have  $|\la k,\omega\ra|\geq \frac{1}{2q_{n+1}}$, then 
$$|E|_{h'} \leq \frac{2}{\pi}\epsilon q_{n+1}.$$
Let $B=e^{Y_1}e^{E}$, $ \tilde A_2=A+\hat F_-^{(re)}(0)J=:2\pi \bar\varrho  J$, $ \tilde F_2=e^{-E}\mathcal R_{\frac{q_{n+1}}{6}}F^{(re)}e^E$. Then we have 
\[\|B\|_{C^1}\leq\frac{1}{h'}\|B\|_{h'}\leq \frac{1}{h'}e^{|Y_1|_{h'}}e^{|E|_{h'}}\leq \frac{2}{h'} e^{\frac{2}{\pi}\epsilon q_{n+1} }<e^{q_{n+1}h'^4},\]
\[\| \tilde F_2\|_{h'/2}\leq e^{4q_{n+1}\epsilon/\pi}\sum_{|k|\geq q_{n+1}/6}|\hat F^{(re)}(k)| e^{ \pi |k|h'}<\epsilon e^{-\frac{q_{n+1}h'}{2}}.\]

\textit{(b)} If $\Lambda_2^c\cap\{k\in\Z^2 : |k|<\frac{q_{n+1}}{6}\}\neq\emptyset$, then by Lemma \ref{lem-hou-you-3},
 there exists $\tilde{\cal Q}_1\in C_{\frac{h'}{3}}^\omega(\T^2, PSL(2,\R))$ that conjugates (\ref{sys-c-after-eliminate-non-resonant}) to
\begin{equation}\label{sys-c-ellip-liouv-after-rotation}
\left\{ \begin{array}{ll} \dot x=  ( \tilde A+\tilde F(\theta) ) x\\
\dot \theta=\omega
\end{array},\right.
\end{equation}
where $$\tilde A =2\pi (\varrho-\frac{\la k_*,\omega\ra}{2} )J, \qquad \mathcal{ T}_{\frac{q_{n+1}}{6}} \tilde F =\sum_{k=l(p_n, -q_n) \atop{|k|< \frac{q_{n+1}}{6}}}\hat { \tilde F} (k)e^{2\pi i\la k,\theta\ra},$$
with estimate $| \tilde F |_{h'/3}\leq 2\epsilon^{3/4}$. Moreover, $\tilde{\cal Q}_1(\theta)=\cal Q(\theta)e^{Y_2}$ takes the form 
 $$\cal Q(\theta)=\left(\begin{array}{cc} \cos(\pi\la k_*,\theta\ra) & \sin(\pi\la k_*,\theta\ra)\\ -\sin(\pi\la k_*,\theta\ra) & \cos(\pi \la k_*,\theta\ra)
\end{array}\right),$$ 
where $k_*\in \Lambda_{2,1}^c\cap \{k\in\Z^2 : |k|<\frac{q_{n+1}}{6}\}$. Furthermore, we have the following estimates
 $$|Y_2|_{\frac{h'}{3}}\leq \epsilon^{\frac{3}{8}}, \qquad \|\cal Q\|_{C^1}\leq 2\pi |k_*|<\frac{\pi}{3}q_{n+1}.$$

Now we consider the system 
\begin{equation}\label{sys-c-ellip-liouv-before-floquet}
\left\{ \begin{array}{ll} \dot x=  ( \tilde A+\cal T_{\frac{q_{n+1}}{6}}\tilde F ) x\\
\dot \theta=\omega
\end{array}.\right.
\end{equation}
By Lemma \ref{lem-hou-you-floquet}, there exist $D\in sl(2,\R)$ and $L_1\in C_{\frac{h'}{3}}^\omega(\T^2, PSL(2,\R))$ that $L_1$ conjugates (\ref{sys-c-ellip-liouv-before-floquet}) to 
\[\left\{ \begin{array}{ll} \dot x= D x\\
\dot \theta=\omega
\end{array}.\right.
\]
Since we have $$|\tilde A+\mathcal T_{\frac{q_{n+1}}{6}}\tilde F|_{\frac{h'}{3}}\leq 4\epsilon^{\frac{1}{4}},$$
then one can compute $\|L_1\|_{\frac{h'}{3}}<e^{64\epsilon^{1/4} h'q_nq_{n+1}}$, and
\begin{equation}\label{cod}
\|D\|\leq \frac{1}{q_{n+1}} e^{16q_{n+1}\epsilon^{1/4}} \leq e^{ \frac{q_{n+1}}{4}\tilde{h}^4}.
\end{equation}  Then under the conjugation map $L_1$, the system (\ref{sys-c-ellip-liouv-after-rotation}) becomes 
\begin{equation}\label{sys-c-ellip-liouv-after-floquet}
\left\{ \begin{array}{ll} \dot x= (D+L_1^{-1}\mathcal R_{\frac{q_{n+1}}{6}}\tilde FL_1) x\\
\dot \theta=\omega
\end{array}.\right.
\end{equation}
By (\ref{equ-parameter}), we have 
\begin{eqnarray*}
\|L_1^{-1}\mathcal R_{\frac{q_{n+1}}{6}}\tilde FL_1\|_{\frac{h'}{6}}\leq e^{128\epsilon^{1/4} q_nh'q_{n+1}}\sum_{|k|\geq \frac{q_{n+1}}{6}}|\hat{\tilde F}(k)|e^{ 2\pi |k|\frac{h'}{6}}<\epsilon^{\frac{3}{4}}e^{-q_{n+1}\tilde h}.
\end{eqnarray*}
Then by  \eqref{cod}, one can   apply  Lemma \ref{lem-c-normalization}, and  there exists $P\in SL(2,\R)$ that conjugates (\ref{sys-c-ellip-liouv-after-floquet}) to (\ref{sys-c-ellip-liouv-final}), where 
\begin{itemize}
\item either $ \tilde A_2=2\pi\bar \varrho J$ for some $\bar\varrho\in\R$ with
\begin{itemize}
\item either $\left\{\begin{array}{ll}
\|P\|\leq 2e^{q_{n+1}\tilde h^2}\\
\|\tilde F_2\|_{\tilde h}\leq 4e^{-\frac{1 }{2}q_{n+1}\tilde h}
\end{array}\right.$,
\item or $\left\{\begin{array}{ll}
\|P\|\leq 2\\
\|\tilde F_2\|_{\tilde h}\leq 2e^{-\frac{ c_0 }{3}q_{n+1}\tilde h^2}
\end{array}\right.$;
\end{itemize}
\item or $\tilde A_2=\left(\begin{array}{cc} 0 & \bar c \\ 0 & 0\end{array}\right)$, where $\bar c\in\R$, with 
\[\left\{\begin{array}{lll}
\|P\|\leq 2e^{\frac{q_{n+1} }{2} \tilde h^4}\\
\|\tilde F_2\|_{ \tilde h}\leq 9e^{-\frac{3 }{4}q_{n+1}\tilde h^2}\\
e^{-\frac{c_0}{3}q_{n+1}\tilde h^2}\leq |\bar c|\leq 4e^{-\frac{3 }{4}q_{n+1}\tilde h^4}
\end{array}.\right.\]
\end{itemize}
Let $B=e^{Y_1}\tilde{\mathcal Q}_1L_1P$. Then we get the desired result.\qed

\subsubsection{Proof of Proposition \ref{propo-c-parab-cd}}


Let $\eta=\epsilon^{\frac{1}{4}}$. Since $2e^{-c_0q_n^{\frac{1}{ \cal A}}h'^2}<10^{-8}$, by  Lemma \ref{lem-hou-you},
there exist $Y_1\in \mathcal B_{h'}$ and $F^{(re)}\in\mathcal B_{h'}^{(re)}(\eta)$ such that (\ref{sys-continuous-original}) can be conjugate to (\ref{sys-c-after-eliminate-non-resonant}) by the conjugation map $e^{Y_1}$, and $|Y_1|_{h'}\leq \epsilon^{1/2}$, $|F^{(re)}|_{h'}\leq 2\epsilon$. Then by Lemma \ref{lem-lambda}, we have $\Lambda_1^c(\eta^{\frac{1}{3}})\cap \{k\in\Z^2 : |k|< q_{n+l} \}=\{0\}$, since $\epsilon^{\frac{1}{12}}<\frac{1}{4q_{n}^{\mathcal A^4}}$ by (\ref{equ-parameter}).  Furthermore, by 
Corollary \ref{cor-parabolic-non-resonant}, 
one has $F^{(re)}=\hat F^{(re)}(0)+\cal R_{ q_{n+l} }F^{(re)}$. Rewrite the system (\ref{sys-c-after-eliminate-non-resonant}) as 
\begin{equation}\label{sys-c-parab-cd-before-normalization}
\left\{\begin{array}{ll}
\dot x=(\tilde A+\tilde F(\theta))x\\
\dot \theta=\omega
\end{array},\right.
\end{equation}
where $\tilde A=A+\hat F^{(re)}(0)$, $\tilde F=\mathcal R_{ q_{n+l} }F^{(re)}$. Then $\|\tilde A\|\leq 2$, and 
\[\|\tilde F\|_{\frac{h'}{6}}\leq  \sum_{|k|\geq q_{n+l} }|\hat F^{(re)}(k)| e^{2\pi |k|\frac{h'}{6}}\leq |F^{(re)}|_{h'}e^{-\frac{5}{3}\pi q_{n+l}h'}<\epsilon e^{-q_{n+l}h'}.\]
By Lemma \ref{lem-c-normalization}, there exists $P\in SL(2,\R)$ that conjugates (\ref{sys-c-parab-cd-before-normalization}) to (\ref{sys-c-parab-cd-final}),  and let $B=e^{Y_1}P$. Then we get the desired estimates. Since the estimates are similar as  estimates in Proposition \ref{propo-c-ellip-cd}, we omit the details.
\qed


\subsubsection{Proof of Proposition \ref{propo-c-parab-liouv}}
Let $\eta=\epsilon^{1/4}$. Since $2e^{-c_0q_n^{\frac{1}{ \cal A}}h'^2}<10^{-8}$, then  by  Lemma \ref{lem-hou-you}, there exists $Y_1\in \mathcal B_{h'}$ and $F^{(re)}\in\mathcal B_{h'}^{(re)}(\eta)$ such that (\ref{sys-continuous-original}) can be conjugate to (\ref{sys-c-after-eliminate-non-resonant}) by the conjugation map $e^{Y_1}$, and $|Y_1|_{h'}\leq \epsilon^{1/2}$, $|F^{(re)}|_{h'}\leq 2\epsilon$. 

On the other hand, note $\epsilon^{\frac{1}{12}}<\frac{1}{7q_n}$ by (\ref{equ-parameter}). Then by Corollary \ref{cor-parabolic-non-resonant} and Remark \ref{remark-liouvillean-lambda-structure}, $F^{(re)}$ takes the form 
$$\mathcal{ T}_{\frac{q_{n+1}}{6}}  F^{(re)} (\theta)=\sum_{k=l(p_n, -q_n) \atop{|k|< \frac{q_{n+1}}{6}}}\hat { F}^{(re)} (k)e^{2\pi i\la k,\theta\ra}.$$
 Now we consider the system
\begin{equation}\label{sys-c-parab-liouv-before-floquet}
\left\{\begin{array}{ll}
\dot x=(  A+\cal T_{\frac{q_{n+1}}{6}} F^{(re)}(\theta))x\\
\dot \theta=\omega
\end{array},\right.
\end{equation}
with $|A+\cal T_{\frac{q_{n+1}}{6}} F^{(re)}|_{h'}\leq |c^*|+|F^{(re)}|_{h'}\leq 3\tilde \epsilon$ with $\tilde\epsilon=2e^{-c_0q_n^{\frac{1}{\cal A}}h'^2}$. 
Then by Lemma \ref{lem-hou-you-floquet}, there exist $D\in sl(2,\R)$ and  $L_2\in C_{h'}^\omega(\T^2, PSL(2,\R))$ that conjugates (\ref{sys-c-parab-liouv-before-floquet}) to 
\[
\left\{ \begin{array}{ll} \dot x= D x\\
\dot \theta=\omega
\end{array},\right.
\]
with $\|L_2\|_{h'}<e^{48\tilde \epsilon h'q_nq_{n+1}}$, and
\begin{equation}\label{cod2}
\|D\|\leq \frac{1}{q_{n+1} }e^{12 \tilde \epsilon q_{n+1}} \leq e^{ \frac{q_{n+1}}{4}\tilde{h}^4}.
\end{equation}
Then under the conjugation map $L_2$, the system (\ref{sys-c-after-eliminate-non-resonant}) becomes 
\begin{equation}\label{sys-c-parab-liouv-after-floquet}
\left\{ \begin{array}{ll} \dot x= (D+L_2^{-1}\mathcal R_{\frac{q_{n+1}}{6}} F^{(re)}L_2) x\\
\dot \theta=\omega
\end{array},\right.
\end{equation}
with estimate
\begin{eqnarray*}
\|L_2^{-1}\mathcal R_{\frac{q_{n+1}}{6}}F^{(re)}L_2\|_{\frac{h'}{6}}&\leq& e^{96\tilde \epsilon h'q_nq_{n+1}}\sum_{|k|\geq \frac{q_{n+1}}{6}}| \hat F^{(re)}(k)| e^{2\pi |k|\frac{h'}{6}}\\
&\leq&e^{96 \tilde\epsilon h'q_nq_{n+1}}e^{-2\pi\frac{q_{n+1}}{6}\frac{5h'}{6}}|F^{(re)}|_{h'}<\epsilon e^{-q_{n+1}\tilde h}.
\end{eqnarray*}
Then by  \eqref{cod2}, one can   apply  Lemma \ref{lem-c-normalization}, and  there exists $P\in SL(2,\R)$ that conjugates (\ref{sys-c-parab-liouv-after-floquet}) to (\ref{sys-c-parab-liouv-final}).  Let  $B=e^{Y_1}L_2P$. Then we get the desired estimates. Since the estimates are similar as  estimates in Proposition \ref{propo-c-ellip-liouv}, we omit the details. \qed


\section{Measure complexity estimation for elliptic case: proof of Proposition \ref{prop-measure-complexity-almost}}\label{sec-measure-complexity}

In this section, we will prove the measure complexity of a cocycle is sub-polynomial if the constant matrices in the conjugated cocycles are elliptic, and the estimates of the conjugations $W_j$ and perturbations $G_j$ satisfy $\|G_j\|_{C^0}^{\eta}\|W_j\|_{C^1}\rightarrow 0$ as $j\rightarrow \infty$ for any $\eta>0$.

The following lemma gives the difference of two distinct points under the projective action of a cocycle. 
\begin{Lemma}\label{lem-conjugate-estimate}
Let $a\in\T^1$, $B\in C^1(\T^1, PSL(2,\R))$. Then for any $(\theta,\varphi), (\tilde \theta,\tilde\varphi)\in\T^1\times\R\mathbb P^1$, we have 
\[d( T_{(a,B)}(\theta,\varphi), T_{(a,B)}(\tilde \theta, \tilde\varphi) )\leq C^* \|B\|_{C^1}^{4} d( (\theta,\varphi), (\tilde \theta, \tilde\varphi) ),\]
where $C^*$ is an absolute constant, and 
$d( (\theta,\varphi), (\tilde\theta, \tilde\varphi)):= \max\{\|\theta-\tilde\theta\|, \|\varphi-\tilde\varphi\|\}$.
In particular, if $B=R_\varrho$, then 
\[d( T_{(a,B)}(\theta,\varphi), T_{(a,B)}(\tilde \theta, \tilde\varphi) )= d( (\theta,\varphi), (\tilde \theta, \tilde\varphi) ).\]
\end{Lemma}
\begin{pf}
Suppose $B(\theta)=\left(\begin{array}{cc} b_{11}(\theta) & b_{12}(\theta)\\ b_{21}(\theta) & b_{22}(\theta)\end{array}\right)$.
Denote $T_{(a,B)}(\theta,\varphi)=(\theta+a, f(\theta,\varphi))$, where $\tan(\hat f(\theta,\varphi))=\frac{b_{21}(\theta)+b_{22}(\theta)\tan\hat\varphi}{b_{11}(\theta)+b_{12}(\theta)\tan\hat\varphi}=\frac{b_{21}(\theta)\cot\hat\varphi+b_{22}(\theta)}{b_{11}(\theta)\cot\hat\varphi+b_{12}(\theta)}$, with $\hat f:=2\pi\gamma(f)$ and $\hat \varphi:=2\pi\gamma(\varphi)$, where $\gamma: \R\mathbb P^1\rightarrow (-\frac{1}{4}, \frac{1}{4}]$ is the lift of the identity map on $\R\mathbb P^1$. Then, we get
\begin{eqnarray*}
\|f(\theta,\varphi)-f(\theta,\tilde\varphi)\|\leq |\frac{1}{2\pi}(\hat f(\theta,\varphi) -\hat f(\theta,\tilde\varphi))|\leq \frac{1}{2\pi}\left\|\frac{\partial \hat f}{\partial \varphi}\right\|_{C^0}\|\varphi-\tilde\varphi\|.
\end{eqnarray*}
Now we give the estimation of $\|\frac{\partial \hat f}{\partial\varphi}\|_{C^0}$: By a simple computation, 
\begin{eqnarray}
\label{est-partial-varphi-1}\frac{\partial \hat f}{\partial\varphi}(\theta,\varphi)&=&\frac{2\pi (1+(\tan\hat\varphi)^2)}{(b_{11}(\theta)+b_{12}(\theta)\tan\hat\varphi)^2+(b_{21}(\theta)+b_{22}(\theta)\tan\hat\varphi)^2}\\
\label{est-partial-varphi-2}&=&\frac{2\pi (1+(\cot\hat\varphi)^2)}{(b_{11}(\theta)\cot\hat\varphi+b_{12}(\theta))^2+(b_{21}(\theta)\cot\hat\varphi+b_{22}(\theta))^2}.
\end{eqnarray}
Moreover, for any $y\in\R$, we have 
 \begin{eqnarray*}
 (b_{11}(\theta)+b_{12}(\theta)y)^2+(b_{21}(\theta)+b_{22}(\theta)y)^2
 &\geq&\frac{1}{b_{12}^2(\theta)+b_{22}^2(\theta)}
 \geq \frac{1}{2\|B\|^2_{C^1}},\\
 (b_{11}(\theta)y+b_{12}(\theta))^2+(b_{21}(\theta)y+b_{22}(\theta))^2
& \geq&\frac{1}{b_{11}^2(\theta)+b_{21}^2(\theta)}
 \geq \frac{1}{2\|B\|^2_{C^1}}.
  \end{eqnarray*}
  Therefore, if $|\tan\hat\varphi|\leq 1$, then by (\ref{est-partial-varphi-1}), we have 
  $\left|\frac{\partial\hat f}{\partial\varphi}\right|\leq 8\pi\|B\|_{C^1}^2.$
  Otherwise, by (\ref{est-partial-varphi-2}), we also have 
  $\left|\frac{\partial\hat f}{\partial\varphi}\right|\leq 8\pi\|B\|_{C^1}^2.$
  Hence, for any $\varphi,\tilde\varphi\in\R\mathbb P^1$,
  \[ \|f(\theta,\varphi)-f(\theta,\tilde\varphi)\|\leq 4\|B\|_{C^1}^{2}\|\varphi-\tilde\varphi\|.\]

 Similarly, we get
 \[\|f(\theta,\varphi)-f(\tilde \theta,\varphi)\|\leq |\frac{1}{2\pi}(\hat f(\theta,\varphi) -\hat f(\tilde \theta,\varphi))|\leq \frac{1}{2\pi}\left\|\frac{\partial \hat f}{\partial \theta}\right\|_{C^0} \|\theta-\tilde \theta\|,\]
 and 
 \begin{eqnarray*}
\left| \frac{\partial \hat f}{\partial \theta}(\theta,\varphi)\right|
 \leq  8\|B\|_{C^1}^2\cdot \frac{1}{2\pi}\left| \frac{\partial\hat f}{\partial\varphi}(\theta,\varphi) \right|
 \leq 32\|B\|_{C^1}^{4},\ \ \forall (\theta,\varphi)\in \T^1\times\R\mathbb P^1.
 \end{eqnarray*}
Therefore, 
\begin{eqnarray*}
\lefteqn{d(T_{(a,B)}(\theta,\varphi), T_{(a,B)}(\tilde \theta,\tilde\varphi))=\max\{\|f(\theta,\varphi)-f(\tilde \theta,\tilde\varphi)\|, \|\theta-\tilde \theta\|\}}\\
&\leq&\max\{ \|f(\theta,\varphi)-f(\theta,\tilde\varphi)\|+\|f(\theta,\tilde\varphi)-f(\tilde \theta,\tilde\varphi)\|, \|\theta-\tilde \theta\|\}\\
&\leq &C^*\|B\|_{C^1}^{4}d((\theta,\varphi), (\tilde \theta, \tilde\varphi)).
\end{eqnarray*}

If $B=R_\varrho$, then we have $T_{(a,B)}(\theta,\varphi)=(\theta+a, \varphi+\varrho)$ for any $(\theta,\varphi)\in \T^1\times\R\mathbb P^1$, which implies 
\[d( T_{(a,B)}(\theta,\varphi), T_{(a,B)}(\tilde\theta,\tilde\varphi))=d((\theta,\varphi), (\tilde\theta,\tilde\varphi)).\]

\end{pf}

\begin{Lemma}\label{lem-distance-perturb}
Let $a\in \T^1$. For any $A\in C^1(\T^1, PSL(2,\R))$, $F\in C^0(\T^1, SL(2,\R))$, 
if $\|F-I\|_{C^0}= \varepsilon\leq \frac{1}{6}$, then 
\[d(T_{(a, AF)}(\theta, \varphi), T_{(a,A)}(\theta,\varphi) )\leq 2C^*\|A\|_{C^1}^{4}\|F-I\|_{C^0}, \ \ \forall \ (\theta,\varphi)\in\T^1\times\R\mathbb P^1,\]
where $C^*$ is the global constant in Lemma \ref{lem-conjugate-estimate}.
\end{Lemma}
\begin{pf}
Since $(a,AF)=(a,A)\circ (0,F)$, we have $T_{(a,AF)}=T_{(a,A)}\circ T_{(0, F)}$. Then by Lemma \ref{lem-conjugate-estimate}, for any $(\theta,\varphi)\in \T^1\times\R\mathbb P^1$,
\begin{eqnarray*}
d(T_{(a,AF)}(\theta,\varphi),T_{(a,A)}(\theta,\varphi))&=&d(T_{(a,A)}(T_{(0,F)}(\theta,\varphi)),  T_{(a,A)}(\theta,\varphi))\\
&\leq& C^*\|A\|_{C^1}^{4}d( T_{(0,F)}(\theta,\varphi), (\theta,\varphi) ).
 \end{eqnarray*}
 Denote $T_{(0,F)}(\theta,\varphi)=:(\theta,\varphi_1)$ and $F(\theta)=\left(\begin{array}{cc}1+f_1(\theta) & f_2(\theta)\\ f_3(\theta) & 1+f_4(\theta)\end{array}\right)$. Let $\hat \varphi:=2\pi\gamma(\varphi)$ and $\hat \varphi_1:=2\pi\gamma(\varphi_1)$, where $\gamma: \R\mathbb P^1\rightarrow (-\frac{1}{4}, \frac{1}{4}]$ is the lift of the identity map on $\R\mathbb P^1$.  
 If $|\sin\hat\varphi|\leq \frac{\sqrt 2}{2}$, then $|\cos\hat\varphi|\geq \frac{\sqrt 2}{2}$ and 
  \begin{eqnarray*}
   \lefteqn{\|\varphi-\varphi_1\|=\frac{1}{2\pi}|\hat\varphi-\hat\varphi_1|\leq\frac{1}{2\pi}|\tan\hat\varphi-\tan\hat\varphi_1|}\\
 &=&\frac{1}{2\pi}\left|\frac{f_3\cos\hat\varphi+(1+f_4)\sin\hat\varphi}{(1+f_1)\cos\hat\varphi+f_2\sin\hat\varphi}-\frac{\sin\hat\varphi}{\cos\hat\varphi}\right|\\
 &=&\frac{1}{2\pi}\left|\frac{(f_4-f_1)\cos\hat\varphi\sin\hat\varphi-f_2\sin^2\hat\varphi+f_3\cos^2\hat\varphi}{(1+f_1)\cos^2\hat\varphi+f_2\sin\hat\varphi\cos\hat\varphi}\right|\\
 &\leq&\frac{12\varepsilon}{2\pi}< 2\varepsilon,
 \end{eqnarray*}
 which implies that 
 \[d(T_{(0,F)}(\theta,\varphi), (\theta,\varphi))=\|\varphi-\varphi_1\|<2\varepsilon.\]
 If $|\sin\hat\varphi|>\frac{\sqrt 2}{2}$, then similarly,
 \begin{eqnarray*}
 \lefteqn{\|\varphi-\varphi_1\|=\frac{1}{2\pi}|\hat\varphi-\hat\varphi_1|\leq\frac{1}{2\pi} |\cot\hat\varphi-\cot\hat\varphi_1|}\\
 &=&\frac{1}{2\pi}\left|\frac{(f_1-f_4)\cos\hat\varphi\sin\hat\varphi+f_2\sin^2\hat\varphi-f_3\cos^2\hat\varphi}{(1+f_4)\sin^2\hat\varphi+f_3\sin\hat\varphi\cos\hat\varphi}\right|\\
 &\leq& \frac{12\varepsilon}{2\pi}<2\varepsilon,
 \end{eqnarray*}
which also implies that 
\[d(T_{(0,F)}(\theta,\varphi), (\theta,\varphi))=\|\varphi-\varphi_1\|<2\varepsilon.\]
Therefore, 
\[d(T_{(a, AF)}(\theta, \varphi), T_{(a,A)}(\theta,\varphi) )\leq 2C^*\|A\|_{C^1}^{4}\|F-I\|_{C^0}.\]
\end{pf}

\begin{Lemma}[\cite{AFK11}]\label{claim-afk}
For matrices $\{M_l\}_{l}\subseteq SL(2,\R), \{ I+\xi_l\}_l\subseteq C(\T^1, SL(2,\R))$, we have that 
\[M_l(I+\xi_l)\ldots M_1(I+\xi_1)=M^{(l)}(I+\xi^{(l)}),\]
with $M^{(l)}=M_l\ldots M_1$ and $\xi^{(l)}$ satisfying 
\[\|\xi^{(l)}\|\leq e^{\sum_{k=1}^l\|M^{(k)}\|^2\|\xi_k\|}-1.\]
\end{Lemma}

By Lemma \ref{claim-afk} and Lemma \ref{lem-distance-perturb}, the following holds:
\begin{Lemma}\label{lem-estimation-cocycle-difference}
Let $A=R_{\varrho}$ with $\varrho\in\T^1$ or $A=\left(\begin{array}{cc}1 & c_* \\ 0& 1\end{array}\right)$ with $|c_*|\leq 1$ and $\tilde A(\theta)=A(I+F(\theta))\in C(\T^1, SL(2,\R))$. 
Then for $0\leq m\leq \frac{1}{4\|F\|^{1/3}_{C^0}}$, we have 
\[d( T_{(\alpha,A)^m}(\theta,\varphi), T_{(\alpha, \tilde A)^m}(\theta,\varphi))\leq 16C^*m^3\|A^m\|^{4}\|F\|_{C^0}, \ \forall\ (\theta,\varphi)\in\T^1\times\R\mathbb P^1.\]

\end{Lemma}
\begin{pf}
By Lemma \ref{claim-afk},
we get
$\tilde A_m(\theta)=A^m(I+F^{(m)}(\theta))$ with 
$\|F^{(m)}(\theta)\|\leq e^{\sum_{k=1}^m\|A^k\|^2\|F\|_{C^0}}-1<8m^3\|F\|_{C^0}<1/6$ for any $\theta\in\T^1$. Then by Lemma \ref{lem-distance-perturb}, we obtain that  
\begin{eqnarray*}
d(T_{(\alpha, A)^m}(\theta,\varphi), T_{(\alpha, \tilde A)^m}(\theta,\varphi))&=&d( T_{(m\alpha, A^m)}(\theta, \varphi), T_{(m\alpha, A^m(I+F^{(m)}))}(\theta,\varphi) )\\
&\leq&2C^*\|A^m\|^{4}\|F^{(m)}\|_{C^0} 
\leq 16C^*m^3\|A^m\|^{4}\|F\|_{C^0}.
\end{eqnarray*}

\end{pf}

\noindent\textit{Proof of Proposition \ref{prop-measure-complexity-almost}.}
Let $\tau, \epsilon>0$. For any $(\theta,\varphi), (\tilde \theta,\tilde\varphi)\in \T^1\times\R\mathbb P^1$, we denote $\varphi_n=\pi_2\circ T_{(\alpha, A)^n}(\theta, \varphi) ,  \tilde\varphi_n=\pi_2\circ T_{(\alpha, A)^n}(\tilde \theta,\tilde\varphi)$. For any $\varepsilon>0$,   there exists $j_*\in\N$ such that $ \forall \ j\geq j_*$,
\[ (6C^*)^{4+\tau}\|W_j\|_{C^1}^{16}\|G_j\|_{C^0}^{\frac{\tau}{4}}<\epsilon^2\varepsilon,\ \ \textrm{and} \ \ \|W_j\|_{C^1}^4\|G_j\|_{C^0}^{\frac{1}{4}}<\epsilon,\] 
where $ C^*\geq 1$ is the global constant in Lemma \ref{lem-conjugate-estimate}. 
Then for $j\geq j_*$, by Lemma \ref{lem-conjugate-estimate}, \ref{lem-distance-perturb}, and \ref{lem-estimation-cocycle-difference}, we have 
\begin{eqnarray*}
\lefteqn{\|\varphi_n-\tilde\varphi_n\|=\|\pi_2\circ T_{(0,W_j)\circ (\alpha, R_{\varrho_j}(I+G_j))^n\circ (0, W_j)^{-1}}(\theta,\varphi)}\\
&&\ \ \ \ \ \ \ \ \ \ -\pi_2\circ T_{(0,W_j)\circ (\alpha, R_{\varrho_j}(I+G_j))^n\circ (0, W_j)^{-1}}(\tilde \theta,\tilde\varphi)\|\\
&\leq & C^*\|W_j\|_{C^1}^{4}d\left(\ T_{(\alpha, R_{\varrho_j}(I+G_j))^n\circ (0, W_j)^{-1}}(\theta,\varphi),T_{(\alpha, R_{\varrho_j}(I+G_j))^n\circ (0, W_j)^{-1}}(\tilde \theta, \tilde\varphi) \ \right)\\
&\leq & C^*\|W_j\|_{C^1}^{4}\left\{ d (\ T_{(\alpha, R_{\varrho_j}(I +G_j))^n}( T_{ (0,W_j)^{-1}}(\theta,\varphi) ),  T_{(\alpha, R_{\varrho_j})^n}(T_{(0, W_j)^{-1}}(\theta, \varphi))\ )\right.\\
&& \ \ \ \ \ \ \ \ \ \ \ \ \ \ + d(\ T_{(\alpha, R_{\varrho_j} )^n\circ (0,W_j)^{-1}}(\theta,\varphi),  T_{(\alpha, R_{\varrho_j})^n\circ (0, W_j)^{-1}}(\tilde \theta, \tilde\varphi)\ )\\
&& \ \ \ \ \ \ \ \ \ \ \ \ \ \ +\left.d( \ T_{(\alpha, R_{\varrho_j}(I +G_j))^n}(T_{(0,W_j)^{-1}}(\tilde \theta,\tilde \varphi)),  T_{(\alpha, R_{\varrho_j})^n}(T_{(0, W_j)^{-1}}(\tilde \theta, \tilde\varphi))\ ) \right\}\\
&\leq& C^*\|W_j\|_{C^1}^{4}\left( 64C^*n^3\|G_j\|_{C^0}+C^*\|W_j\|_{C^1}^{4}d((\theta,\varphi), (\tilde \theta, \tilde\varphi)) \right)
\end{eqnarray*}
for $0\leq n\leq \frac{1}{4\|G_j\|^{1/3}_{C^0}}$. Now let $M_j=\left\lceil\frac{1}{6C^{*}\|G_j\|_{C^0}^{1/4}}\right\rceil$. Then for $d((\theta,\varphi), (\tilde \theta, \tilde\varphi) )<\frac{\epsilon}{2 C^{*2}\|W_j\|_{C^1}^{8}}$, we have
\[\bar d_{M_j}( (\theta, \varphi), (\tilde \theta,\tilde\varphi) )<\epsilon.\]
If we let $L_j:=\left\{(k_1\cdot \frac{\epsilon}{2C^{*2}\|W_j\|_{C^1}^{8}}, k_2\cdot \frac{\epsilon}{2C^{*2}\|W_j\|_{C^1}^{8}})\ :\ 0\leq k_1, k_2\leq \left[\frac{2 C^{*2}\|W_j\|_{C^1}^{8}}{\epsilon}\right]+1\right\}$, then for any $(\theta, \varphi)\in \T^1\times \R\mathbb P^1$, there exists $(\tilde \theta,\tilde \varphi)\in L_j$ such that $d((\theta, \varphi), (\tilde \theta, \tilde\varphi))<\frac{\epsilon}{2 C^{*2} \|W_j\|_{C^1}^{8}}$, which implies that $S_{M_j}(d,\rho,\epsilon)< \frac{5 C^{*4}\|W_j\|_{C^1}^{16}}{\epsilon^2}$ for any $\rho\in\mathcal M(\T^1\times\R\mathbb P^1)$ . Hence, we have
\[\frac{S_{M_j}(d,\rho,\epsilon)}{M_j^\tau}<\frac{6^{\tau+1} C^{*(4+\tau)}}{\epsilon^2}\|G_j\|_{C^0}^{\frac{\tau}{4}}\|W_j\|_{C^1}^{16}<\varepsilon,\ \ \textrm{as}\ j\geq j_*.\]
Therefore, for any $ \rho\in\cal M(\T^1\times\R\mathbb P^1)$ and $\tau>0$, we have 
\[\liminf_{n\rightarrow\infty}\frac{S_n(d,\rho,\epsilon)}{n^\tau}=0,\ \ \forall \epsilon>0,\]
which means that the complexity of $(\T^1\times\R\mathbb P^1, T_{(\alpha, A)}, \rho)$ is sub-polynomial.
\qed

\section{M\"obius disjointness for the parabolic case: proof of Proposition \ref{prop-parabolic}}\label{sec-parabolic}

The main ideas of the proof are developed from  \cite{W17}. In \cite{W17}, within a reasonably long interval, any orbit of $T$ will come back to a neighborhood of the starting point in an almost periodic manner. However, in our case, the conjugated cocycle might be the perturbation of a parabolic matrix, and consequently  those points with the second variable $\varphi$ not close to 0 will leave the neighborhood of the starting 
point without returning in the future, which is the main difference compared to \cite{W17}. The key observation here is that in this situation the 
orbit of the second variable will approach 0. Inspired by this, we approximate the orbit of $(\theta,\varphi)$ by the orbit of $(\theta,0)$ by a periodic sequence with period $q_{k_j}$.
Then we decompose the periodic sequence into short average of Dirichlet characters\footnote{We say a function $\chi: \Z\rightarrow \C$ is a Dirichlet character modulus $k$ if it has the following properties: (1) $\chi(n)=\chi(n+k)$ for all $n$; (2) If gcd(n,k)>1 then $\chi(n)=0$, if gcd(n,k)=1 then $\chi(n)\neq 0$; (3) $\chi(mn)=\chi(m)\chi(n)$ for all integers $m$ and $n$.}, reducing the problem to control the average of multiplicative function on a typical interval, which is the content of the following lemma.

%

\begin{Lemma}[\cite{W17}]\label{lem-W}
For all $L, Q, M\in\N$ and any periodic function $D: \N\rightarrow \C$ of period $Q$ with $|D|\leq 1$,
\[\left|\E_{L\leq n<L+MQ}\mu(n)D(n)\right|^2\leq Q \E_{d | Q, \chi \mod^*\frac{Q}{d}} \left|\E_{\frac{L}{d}\leq r<\frac{L}{d}+M\frac{Q}{d}}\mu(r)\chi(r)\right|^2,\]
where the first average on the right hand side is taken over all pairs $(d,\chi)$ such that $d|Q$ and $\chi$ is a Dirichlet character of conductor $\frac{Q}{d}$.
\end{Lemma}

So as to control $\left|\E_{\frac{L}{d}\leq r<\frac{L}{d}+M\frac{Q}{d}}\mu(r)\chi(r)\right|$, we will use a recent result in \cite{MRT15} that allows to bound averages of a non-pretentious multiplicative function in random short intervals. We first recall the definition of pretentiousness.

For multiplicative functions $\nu, \nu' : \N \rightarrow \C$ with $|\nu|, |\nu'|\leq 1$, define 
$\D(\nu, \nu', X)=\left(\sum_{p\leq X}\frac{1-\Re(\nu(p)\overline{\nu'(p)})}{p}\right)^{\frac{1}{2}}$, and the function below measures how closely $\nu$ pretends to be $n^{it}$:
\[M(\nu,X)=\inf_{|t|\leq X}\D(\nu, n^{it}, X)^2.\]
By the same discussions in \cite{W17}, we know that 
\[\lim_{X\rightarrow \infty} M(\mu\chi, X)=\infty.\]
We will apply the following proposition to $\nu(n)=\mu(n)\chi(n)$.
\begin{Proposition}[\cite{MRT15}]\label{prop-MRT}
Let $\nu$ be a multiplicative function with $|\nu|\leq 1$ and $X\geq l\geq 10$. Then
\[\E_{X\leq L<2X}\Big|\E_{L\leq n<L+l}\ \nu(n)\Big|^2\lesssim e^{-M(\nu, X)}M(\nu,X)+(\log X)^{-\frac{1}{50}}+\left(\frac{\log\log l}{\log l}\right)^2.\]
\end{Proposition}

\noindent\textit{Proof of Proposition \ref{prop-parabolic}.}
Since trigonometric polynomials are dense in $C^0(\T^1\times\R \mathbb P^1)$, it suffices to prove the orthogonality for every $f(\theta,\varphi)=e^{2\pi i(\iota_1\theta+2\iota_2\varphi)}, \iota_1, \iota_2\in\Z$. Fix a function $f(\theta,\varphi)=e^{2\pi i(\iota_1\theta+2\iota_2\varphi)}$ and $ 0<\eta\ll1$. Suppose $k_j$ is sufficiently large and $q_{k_j+1}\geq e^{\beta q_{k_j}/2}$.  
Let $\tilde C$ be the constant in the proposition. Then for $j$ large enough, we have 
$6\tilde C^2e^{-(\frac{\beta}{4}-2\tilde C\eta_j-\tau_j)q_{k_j}}<1,$ $\frac{16\tilde C\eta_j}{\tau_j}<\xi$, and $4\tilde Cq_{k_j}<e^{\frac{1}{7}\xi q_{k_j}^{1/2}}.$
For any $(\theta,\varphi)\in \T^1\times\R\mathbb P^1$, denote 
\[x_0=(\theta,\varphi), \ \ \tilde x_0= T_{(0,W_j)^{-1}}x_0=:(\theta,\tilde\varphi), \ \  T=T_{(\alpha, A)}, \ \ \tilde T=T_{(\alpha, A_j(I+G_j))}, \  \  \tilde x_n=\tilde T^n\tilde x_0. \]
 Let $N_0<\frac{N}{2}$ be an integer. Then,
\begin{eqnarray*}
\lefteqn{\E_{n<N}\mu(n)f(T^n x_0)=\E_{n=N_0}^{N-1}\mu(n)f(T^nx_0)+O(\frac{N_0}{N})}
\\&=&\E_{L=N_0}^{N-1}\E_{n=0}^{M_jq_{k_j}-1}\mu(L+n)f(T^{L+n}x_0)+O(\frac{M_jq_{k_j}}{N})+O(\frac{N_0}{N})\\
&=& \E_{L=N_0}^{N-1}\E_{n=0}^{q_{k_j}-1}\E_{m=0}^{M_j-1}\mu(L+n+mq_{k_j})f(T^{L+n+mq_{k_j}}x_0)+O(\frac{M_jq_{k_j}}{N})+O(\frac{N_0}{N}).\end{eqnarray*}

Now we approximate $f(T^{L+n+mq_{k_j}}x_0)$ by some periodic sequence. We first approximate it by $T_{(0,W_j)}(T_{(0,A_j)^{mq_{k_j}}}\tilde x_{L+n})$. Then use the nature of parabolic matrix $A_j$ to count the number of $m$ where $0\leq m\leq M_j-1$ that the second variable of $T_{(0,A_j)^{mq_{k_j}}}\tilde x_{L+n}$ is away from $0$, and show that it is quite small comparing to $M_j$.

\begin{Lemma}
For any $L,m,n\in\N$ with $m<M_j\leq e^{\frac{1}{7}-\frac{\xi}{2}\tau_jq_{k_j}}$, we have 
\begin{eqnarray*}
\lefteqn{d\left(T^{L+n+mq_{k_j}}x_0, \ T_{(0,W_j)}(T_{ (0,A_j)^{mq_{k_j}}}\tilde x_{L+n}) \right)}
\\ &\lesssim M_j^3q_{k_j}^3(1+M_jq_{k_j}|c_j|)^4\|W_j\|_{C^1}^4\|G_j\|_{C^0}
+\frac{M_j\|W_j\|_{C^1}^6}{q_{k_j+1}},
\end{eqnarray*}
where $A_j=\left(\begin{array}{cc}1 & c_j\\ 0 & 1\end{array}\right)$ with $j$ large enough.
\end{Lemma}
\begin{pf}
By the fact that $(0,W_j)^{-1}\circ (\alpha, A)\circ (0,W_j)=(\alpha, A_j(I+G_j))$, we have $T^{L+n+mq_{k_j}}x_0=T_{(0,W_j)}(T_{(\alpha, A_j(I+G_j))^{mq_{k_j}}}(\tilde x_{L+n}))$. Then,  
\begin{eqnarray*}
\lefteqn{d\left(T^{L+n+mq_{k_j}}x_0, \ T_{(0,W_j)}(T_{ (0,A_j)^{mq_{k_j}}}\tilde x_{L+n}) \right)}\\
&\leq& d\left(T_{(0,W_j)}(T_{(\alpha, A_j(I+G_j))^{mq_{k_j}}}\tilde x_{L+n}), \ T_{(0, W_j(\cdot-mq_{k_j}\alpha))}(T_{(\alpha, A_j)^{mq_{k_j}}}\tilde x_{L+n})\right)\\
&&+d\left( T_{(0, W_j(\cdot-mq_{k_j}\alpha))}(T_{(\alpha, A_j)^{mq_{k_j}}}\tilde x_{L+n}), \ T_{(0,W_j)}(T_{ (0,A_j)^{mq_{k_j}}}\tilde x_{L+n}) \right).
\end{eqnarray*}
Since $A_j$ is a constant matrix, then $(\alpha, A_j)^{mq_{k_j}}=(mq_{k_j}\alpha, A_j^{mq_{k_j}})$, $(0, A_j)^{mq_{k_j}}=(0, A_j^{mq_{k_j}})$, and thus $(0, W_j(\cdot-mq_{k_j}\alpha))\circ (\alpha, A_j)^{mq_{k_j}}=(mq_{k_j}\alpha, W_j(\cdot)A_j^{mq_{k_j}})$, $(0, W_j)\circ (0,A_j)^{mq_{k_j}}=(0, W_j(\cdot)A_j^{mq_{k_j}})$. Therefore, 
\[d\left( T_{(0, W_j(\cdot-mq_{k_j}\alpha))}(T_{(\alpha, A_j)^{mq_{k_j}}}\tilde x_{L+n}), \ T_{(0,W_j)}(T_{ (0,A_j)^{mq_{k_j}}}\tilde x_{L+n}) \right)\leq \|mq_{k_j}\alpha\|_{\T}\leq \frac{m}{q_{k_j+1}}.\]
Denoting $\Delta W_{j,m}(\cdot)=W_{j}(\cdot)^{-1}W_{j}(\cdot-mq_{k_j}\alpha)$, then
for $0\leq m< M_j$ with $j$ sufficiently large,
\begin{eqnarray*}
\lefteqn{\|\Delta W_{j,m}-I\|_{C^0}}\\&\leq& \|W_j\|_{C^0}\cdot \|W_j(\cdot-mq_{k_j}\alpha)-W_j(\cdot)\|_{C^0}\\
&\leq &\|W_j\|_{C^1}\cdot \|W_j\|_{C^1}\cdot \|mq_{k_j}\alpha\|_{\T}\leq\frac{m\|W_j\|_{C^1}^2}{q_{k_j+1}}\\
&\leq& \tilde C^2e^{-(\frac{\beta}{2}-2\tilde C\eta_j-\tau_j)q_{k_j}}<\frac{1}{6}.
\end{eqnarray*} 
By Lemma \ref{lem-distance-perturb}, 
\begin{eqnarray*}
&&  d\left( T_{(0, W_j)}(T_{(\alpha, A_j)^{mq_{k_j}}}\tilde x_{L+n}), \  T_{(0, W_j(\cdot)\Delta W_{j,m}(\cdot))}(T_{(\alpha, A_j)^{mq_{k_j}}}\tilde x_{L+n})\right)\\
& \leq& 2C^*\|W_j\|_{C^1}^4\|\Delta W_{j,m}-I\|_{C^0}\lesssim \frac{m\|W_j\|_{C^1}^6}{q_{k_j+1}}.
\end{eqnarray*}
Moreover, by Lemma \ref{lem-conjugate-estimate} and \ref{lem-estimation-cocycle-difference}, we have 
\begin{eqnarray*}
&&d\left( T_{(0,W_j)}(T_{(\alpha, A_j(I+G_j))^{mq_{k_j}}}\tilde x_{L+n}), \  T_{(0, W_j)}(T_{(\alpha, A_j)^{mq_{k_j}}}\tilde x_{L+n})\right)\\
&\leq&C^*\|W_j\|_{C^1}^4d\left(T_{(\alpha, A_j(I+G_j))^{mq_{k_j}}}\tilde x_{L+n}, \ T_{(\alpha, A_j)^{mq_{k_j}}}\tilde x_{L+n}\right)\\
&\leq& 16{C^*}^2(mq_{k_j})^3\|W_j\|_{C^1}^4\|A_j^{mq_{k_j}}\|^4\|G_j\|_{C^0}  \\ 
&\lesssim&    m^3q_{k_j}^3(1+mq_{k_j}|c_j|)^4\|W_j\|_{C^1}^4\|G_j\|_{C^0},  
\end{eqnarray*}
since $4m\|G_j\|_{C^0}^{\frac{1}{3}}\leq 4e^{\frac{1}{7}-\frac{\xi}{2}\tau_jq_{k_j}}\tilde Ce^{-\tau_jq_{k_j}}<1$ for $j$ large enough. 

In conclusion, we get
\begin{eqnarray*}
&&d\left(T^{L+n+mq_{k_j}}x_0, \ T_{(0,W_j)}(T_{ (0,A_j)^{mq_{k_j}}}\tilde x_{L+n}) \right) \\
&\leq& d\left( T_{(0,W_j)}(T_{(\alpha, A_j(I+G_j))^{mq_{k_j}}}\tilde x_{L+n}), \  T_{(0, W_j)}(T_{(\alpha, A_j)^{mq_{k_j}}}\tilde x_{L+n})\right)\\
&& + d\left( T_{(0, W_j)}(T_{(\alpha, A_j)^{mq_{k_j}}}\tilde x_{L+n}), \  T_{(0, W_j(\cdot)\Delta W_{j,m}(\cdot))}(T_{(\alpha, A_j)^{mq_{k_j}}}\tilde x_{L+n})\right) \\
&&+d\left( T_{(0, W_j(\cdot-mq_{k_j}\alpha))}(T_{(\alpha, A_j)^{mq_{k_j}}}\tilde x_{L+n}), \ T_{(0,W_j)}(T_{ (0,A_j)^{mq_{k_j}}}\tilde x_{L+n}) \right)\\
&\lesssim&M_j^3q_{k_j}^3(1+M_jq_{k_j}|c_j|)^4\|W_j\|_{C^1}^4\|G_j\|_{C^0}
+\frac{M_j\|W_j\|_{C^1}^6}{q_{k_j+1}}.
\end{eqnarray*}
\end{pf}

%
%

Moreover, for any $0<\tilde\eta\ll 1, (\bar\theta,\bar\varphi)\in\T^1\times\R\mathbb P^1$, we denote $$\bar\varphi^{(m)}:=\pi_2\circ T_{(0,A_j)^{mq_{k_j}}}(\bar\theta,\bar\varphi)=\pi_2\circ T_{(0,A_j^{q_{k_j}})^{m}}(\bar\theta,\bar\varphi)$$ and $I_0(\bar\theta,\bar\varphi,\tilde\eta):=\{m\in\N : |\bar\varphi^{(m)}|>\tilde\eta\}$. 
Let $v_*=(a_*,1)^T\in \R^2$ be any vector and its angle is $\varphi_*$. Then the length for $a_*$ such that $|\tan\varphi_*|>\tilde\eta$ is no more than $\frac{2}{\tilde\eta}$. Moreover, we have $A_j^{q_{k_j}}v_*=(a_*+q_{k_j}c_j, 1)^{T}$, implying that the effect of $A_j^{q_{k_j}}$ on $v_*$ is to add $q_{k_j}c_j$ on its first variable. Since $I_0(\bar\theta,\bar\varphi,\tilde\eta)\subseteq\{m\in\N : |\tan(\bar\varphi^{(m)})|>\tilde\eta\}$,
then we obtain that 
$$\# I_0(\bar\theta,\bar\varphi,\tilde\eta)\leq\frac{2}{\tilde\eta q_{k_j}|c_j|}+1.$$ Therefore, we have 
\begin{eqnarray*}
 &\lefteqn{\E_{L=N_0}^{N-1}\E_{n=0}^{q_{k_j}-1}\E_{m=0}^{M_j-1}\mu(L+n+mq_{k_j})\left(f\circ T_{(0,W_j)}(T_{(0,A_j)^{mq_{k_j}}}\tilde x_{L+n})-f\circ T_{(0,W_j)}(\theta_{ L+n},0)\right)}\\
&=&\E_{L=N_0}^{N-1}\E_{n=0}^{q_{k_j}-1}\frac{1}{M_j}\left(\sum_{\substack{0\leq m< M_j\\ m\in I_0(\theta_{L+n},\tilde\varphi_{L+n},\tilde\eta)}}+\sum_{\substack{0\leq m< M_j\\ m\notin I_0(\theta_{L+n},\tilde\varphi_{L+n},\tilde\eta)}}\right)\\
&&\mu(L+n+mq_{k_j})\left(f\circ T_{(0,W_j)}(T_{(0,A_j)^{mq_{k_j}}}\tilde x_{ L+n})-f\circ T_{(0,W_j)}(\theta_{ L+n},0)\right).\ \ \ \ \ \ \ \ \ \ \ \ \ \ \ \ \  \ \ \ \ \ \ 
\end{eqnarray*}
Since $\# I_0(\bar\theta,\bar\varphi,\tilde\eta)\leq\frac{2}{\tilde\eta q_{k_j}|c_j|}+1$ for any $0<\tilde\eta\ll 1, (\bar\theta,\bar\varphi)\in\T^1\times\R\mathbb P^1$, then 
\begin{eqnarray*}
\lefteqn{\frac{1}{M_j}\sum_{\substack{0\leq m< M_j\\ m\in I_0(\theta_{L+n},\tilde\varphi_{L+n},\tilde\eta)}}\mu(L+n+mq_{k_j})\left(f\circ T_{(0,W_j)}(T_{(0,A_j)^{mq_{k_j}}}\tilde x_{L+n})-f\circ T_{(0,W_j)}(\theta_{ L+n},0)\right)}\\
&\leq &\frac{2 }{M_j}(\frac{2}{\tilde\eta q_{k_j}|c_j|}+1)\lesssim\frac{1}{M_j\tilde\eta q_{k_j}|c_j|}.\ \ \ \ \ \ \ \ \ \ \ \ \ \ \ \ \ \ \ \ \ \ \ \ \ \ \ \ \ \ \ \ \ \ \ \ \ \ \ \ \ \ \ \ \ \ \ \ \ \ \ \ \ \ \ \ \ \ \ \ \ \ \ \ \ \ \ \ \ 
\end{eqnarray*}
Moreover, for those $m$ such that $|\tilde\varphi^{(m)}_{L+n}|\leq\tilde\eta$, we have 
\begin{eqnarray*}
\lefteqn{|f\circ T_{(0,W_j)}(T_{(0,A_j)^{mq_{k_j}}}\tilde x_{L+n })-f\circ T_{(0,W_j)}(\theta_{L+n})|}\\
&\leq& \|f\|_{C^1}d\left(T_{(0,W_j)}(T_{(0,A_j)^{mq_{k_j}}}\tilde x_{L+n}), T_{(0,W_j)}(\theta_{L+n},0)\right)\\
&\overset{Lemma ~\ref{lem-conjugate-estimate}}{\leq}  &
C_*\|f\|_{C^1}\|W_j\|_{C^1}^4d\left( (\theta_{L+n},\tilde\varphi_{L+n}^{(m)}), (\theta_{L+n},0)\right)\\
&\lesssim& \|W_j\|_{C^1}^4\tilde\eta,
\end{eqnarray*}
and thus 
\begin{eqnarray*}
&&\frac{1}{M_j}\sum_{\substack{0\leq m< M_j\\ m\notin I_0(\theta_{L+n},\tilde\varphi_{L+n},\tilde\eta)}}\mu(L+n+mq_{k_j})\left(f\circ T_{(0,W_j)}(T_{(0,A_j)^{mq_{k_j}}}\tilde x_{L+n})-f\circ T_{(0,W_j)}(\theta_{L+n},0)\right)\\
&\lesssim& \|W_j\|_{C^1}^4\tilde\eta.
\end{eqnarray*}
In conclusion, we get
\begin{eqnarray*}
&&\lefteqn{\E_{n<N}\mu(n)f(T^n x_0)}\\
&=&\E_{L=N_0}^{N-1}\E_{n=0}^{q_{k_j}-1}\E_{m=0}^{M_j-1}\mu(L+n+mq_{k_j})f\circ T_{(0,W_j)}(\theta_{L+n},0)\\
& +&O(\|W_j\|_{C^1}^{4}\tilde\eta+\frac{1}{M_j\tilde\eta q_{k_j}|c_j|}+
M_j^7q_{k_j}^7|c_j|^4\|W_j\|_{C^1}^4\|G_j\|_{C^0}+\frac{M_j\|W_j\|_{C^1}^6}{q_{k_j+1}}+\frac{M_jq_{k_j}}{N}+\frac{N_0}{N}).
\end{eqnarray*}
Let $D_l:=f\circ T_{(0, W_j)}(\theta_l,0)$, and for each $L\in \N$ construct a function $D_L : \N\rightarrow \C$ by $D_L(n)=D_l$ where $l$ is the unique integer in $[L, L+q_{k_j})$ such that $l=n (\mod q_{k_j})$. Then $D_L$ is periodic with period $q_{k_j}$ and $|D_L|=1$.

 Now we let $\tilde \eta=e^{-4\tilde Cq_{k_j}\eta_j}\eta\leq \eta$, $M_j=[e^{(\frac{1}{7}-\frac{\xi}{2})\tau_jq_{k_j}}]$ and $J=\lceil \log_2\frac{1}{\eta} \rceil$. Choose $N_0=\lfloor 2^{-J}N\rfloor$, and then $\frac{N_0}{N}\leq \eta$. We have for $j$ sufficiently large,
\begin{eqnarray*}
\lefteqn{| \E_{n<N}\mu(n)f(T^n x_0)|}\\
&=&|\E_{L=N_0}^{N-1}\E_{n=L}^{L+M_jq_{k_j}-1}\mu(n)D_L(n)|+O(\|W_j\|_{C^1}^{4}\tilde\eta+\frac{N_0}{N})\\
&&+O(\frac{1}{M_j\tilde\eta q_{k_j}|c_j|}+
M_j^7q_{k_j}^7|c_j|^4\|W_j\|_{C^1}^4\|G_j\|_{C^0}+\frac{M_j\|W_j\|_{C^1}^6}{q_{k_j+1}})+O(\frac{M_jq_{k_j}}{N})\\
&\overset{Cauchy-Schwarz}{\lesssim} & \left(\E_{L=N_0}^{N-1}\left|\E_{n=L}^{L+M_jq_{k_j}-1}\mu(n)D_L(n)\right|^2\right)^{\frac{1}{2}} +\eta+\frac{e^{-q_{k_j}(\frac{\xi}{2}\tau_j-4\tilde C\eta_j)}}{\eta q_{k_j}}\\
&&+q_{k_j}^7e^{-q_{k_j}(\frac{7}{2}\xi\tau_j-4\tilde C\eta_j)}+
e^{-q_{k_j}(\frac{\beta}{2}-(\frac{1}{7}-\frac{\xi}{2})\tau_j-6\tilde C\eta_j)}+
\frac{q_{k_j} e^{(\frac{1}{7}-\frac{\xi}{2})\tau_jq_{k_j}}}{N}\\
&\overset{Lemma ~\ref{lem-W}}{\lesssim}  &\left( q_{k_j}\cdot \E_{L=N_0}^{N-1}
\E_{d | q_{k_j}, \chi\mod^*\frac{q_{k_j}}{d}}\left|\E_{n=\frac{L}{d}}^{\frac{L}{d}+\frac{M_jq_{k_j}}{d}-1}\mu(n)\chi(n)\right|^2\right)^{\frac{1}{2}}\\
&&+\eta+\frac{e^{-\frac{\xi}{4}\tau_jq_{k_j}}}{\eta}+e^{-2\xi\tau_jq_{k_j}}+e^{-\frac{\beta}{4}q_{k_j}}+\frac{ e^{(\frac{1}{7}-\frac{5}{14}\xi)\tau_jq_{k_j}}}{N}.
\end{eqnarray*}

Cut $[N_0, N)$ into J dyadic interval: $[2^{-\iota}N, 2^{-\iota+1}N)$ for $\iota=1,\ldots, J$. Let $\tilde\rho_{\chi}(X):=e^{-M(\mu\chi, X)}M(\mu\chi, X)+(\log X)^{-\frac{1}{50}}$, $\tilde\rho_{q_{k_j}}(X):=\max_{\chi \mod q_{k_j}}\tilde\rho_\chi(X)$, and  $\rho_{q_{k_j}}(X):=$\\
$\sup_{X'\geq X}\tilde\rho_{q_{k_j}}(X)$. Then $\rho_{q_{k_j}}$ is a positive function independent of $\chi$ that decreases to 0 as $X\rightarrow \infty$. 
By Proposition \ref{prop-MRT}, for each pair $(d,\chi)$ that $d | q_{k_j} $ and $\chi$ is a Dirichlet character of conductor $\frac{q_{k_j}}{d}$, and every $\iota\in [1, J]$,
\begin{eqnarray*}
\lefteqn{\E_{2^{-\iota}N\leq L< 2^{-\iota+1}N}\left|\E_{\frac{L}{d}\leq r<\frac{L}{d}+\frac{M_jq_{k_j}}{d}}\mu(r)\chi(r)\right|^2}\\
\nonumber&\lesssim&  \rho_{q_{k_j}}\left(\frac{2^{-\iota}N}{d}\right)+\left(\frac{\log\log\frac{M_jq_{k_j}}{d}}{\log \frac{M_jq_{k_j}}{d}}\right)^2\\
\nonumber&\lesssim& \rho_{q_{k_j}}\left(\frac{\eta N}{2q_{k_j}}\right)+ \left(\frac{\log\log M_j}{\log M_j}\right)^2\\
\nonumber&\lesssim& \rho_{q_{k_j}}\left(\frac{\eta N}{2q_{k_j}}\right)+\frac{\log^2( \tau_jq_{k_j})}{\tau_j^2q_{k_j}^2 }.
\end{eqnarray*}
Therefore, we have
\begin{eqnarray*}
\lefteqn{| \E_{n<N}\mu(n)f(T^n x_0) |}\\
&\lesssim& \left( q_{k_j}\cdot \rho_{q_{k_j}}(\frac{\eta N}{2q_{k_j}}) + \tau_j^{-2}q_{k_j}^{-1}\log^2(\tau_jq_{k_j}) \right)^{\frac{1}{2}}+ \eta+\frac{e^{-\frac{\xi}{4}\tau_jq_{k_j}}}{\eta} +\frac{ e^{(\frac{1}{7}-\frac{5}{14}\xi)\tau_jq_{k_j}}}{N}.
\end{eqnarray*}
Once $\eta$ is fixed, then for sufficiently large $q_{k_j}$ with $q_{k_j+1}\geq e^{\beta q_{k_j}/2}$, 
\[ \tau_j^{-2}q_{k_j}^{-1}\log^2(\tau_jq_{k_j})<q_{k_j}^{-2\epsilon}\log^2q_{k_j}\lesssim\eta^2, \ \ e^{-\frac{\xi}{4}\tau_jq_{k_j}}\lesssim\eta^2.\]
Fix such a $q_{k_j}$, and then
\[|\E_{n<N}\mu(n)f(T^n x_0)|\lesssim (q_{k_j}\cdot \rho_{q_{k_j}}(\frac{\eta N}{2q_{k_j}})+\eta^2)^{1/2}+\eta+\frac{ e^{(\frac{1}{7}-\frac{5}{14}\xi)\tau_jq_{k_j}}}{N}.\]
Since $\eta, q_{k_j}, M_j$ are now all fixed and $\rho_{q_{k_j}}$ is a function that decays to 0, then for sufficiently large $N$, we have $q_{k_j}\cdot \rho_{q_{k_j}}(\frac{\eta N}{2q_{k_j}})\lesssim\eta^2$ and $\frac{ e^{(\frac{1}{7}-\frac{5}{14}\xi)\tau_jq_{k_j}}}{N}\lesssim\eta$, which implies the result. 

\qed

\section*{Acknowledgments}
The authors would like to thank Tianyuan Mathematical Center in Southwest (No. 11826102).
W. Huang was supported by NSFC grant (12090012, 12031019, 11731003). J. Wang was supported by NSFC grant (11971233), the Outstanding Youth Foundation of Jiangsu Province (No. BK20200074), and Qing Lan Project of Jiangsu province. Z. Wang was supported by NSF grant (DMS-1753042). Q. Zhou was partially supported by National Key R\&D Program of China (2020YFA0713300), NSFC grant (12071232), The Science Fund for Distinguished Young Scholars of Tianjin (No. 19JCJQJC61300) and Nankai Zhide Foundation.

\end{document}